\documentclass[11pt,a4paper, reqno]{amsart}%{article}
\usepackage[a4paper, left = 2cm,  hmargin={25mm,25mm},vmargin={25mm,25mm}]{geometry}
\usepackage[utf8]{inputenc}
\usepackage{geometry}
\usepackage[english]{babel}
\usepackage{graphicx}
\usepackage{color}
\usepackage{xcolor}
\usepackage{pdfpages}
\usepackage{amsbsy}
\usepackage{amssymb}
\usepackage{amsmath}
\usepackage{bm}
\usepackage{mathtools}
\usepackage{matlab-prettifier}
\usepackage{amsthm}
\usepackage{tikz-cd}
\usepackage{amsfonts}
\usepackage{array}
\usepackage{hyperref}
\usepackage[OT2,T1]{fontenc}
\usepackage{csquotes}
\usepackage{mathabx}
\usepackage{bm}
\usepackage{url}
\usepackage{mathdots}
\usepackage{mathrsfs}

\numberwithin{equation}{section}

\newcommand\m[1]{\begin{pmatrix}#1\end{pmatrix}} 

\newcommand{\ol}{\overline}
\newcommand{\vol}{\textup{vol}}

\DeclareSymbolFont{cyrletters}{OT2}{wncyr}{m}{n}
\DeclareMathSymbol{\Sha}{\mathalpha}{cyrletters}{"58}
\newcolumntype{P}[1]{>{\centering\arraybackslash}p{#1}}
\everymath{\displaystyle}
\newtheorem{theorem}{Theorem}[section]

\newtheorem{lemma}[theorem]{Lemma}

\newtheorem{defn}[theorem]{Definition}

\newtheorem{proposition}[theorem]{Proposition}

\newtheorem{remark}[theorem]{Remark}

\newtheorem*{notation}{Notation}
\newtheorem{corollary}[theorem]{Corollary}

\newcommand{\su}{\textup{SU}}
\newcommand{\un}{\textup{U}}

\newcommand{\gl}{\textup{GL}}

\newcommand{\gu}{\textup{GU}}
\newcommand{\mf}{\mathfrak}
\newcommand{\wh}{\widehat}
\newcommand{\mc}{\mathcal}
\newcommand{\vphi}{\varphi}
\newcommand{\mb}{\mathbb}
\newcommand{\mbf}{\mathbf}
\newcommand{\msc}{\mathscr}

\definecolor{lightgray}{gray}{0.5}
\usepackage[backend=biber, style=alphabetic,]{biblatex}
\begin{document}
\title{The action of Hecke operators on Hermitian theta series}
\author{Thanasis Bouganis and Rafail Psyroukis}
\address{Department of Mathematical Sciences\\ Durham University\\
South. Rd.\\ Durham, DH1 3LE, U.K.}
\email{athanasios.bouganis@durham.ac.uk}
\address{University of Primorska\\ Faculty of Mathematics, Natural Sciences and Information Technologies \\ Glagoljaška ulica 8 \\ 6000\\ Koper \\Slovenia}
\email{rafail.psyroukis@famnit.upr.si}
.
\begin{abstract}
We prove $\mathrm{GU}$-Hecke stability results for spaces spanned by Hermitian theta series associated to positive definite lattices over an arbitrary quadratic imaginary number field. We also prove the stability of the space spanned by theta series attached to lattices in a fixed genus under a certain subalgebra of $\mathrm{GU}$-Hecke operators, the even subalgebra. Finally, we explicitly compute the theta matrices for certain Hecke operators in terms of representation numbers in this Hermitian setting.
% of spaces spanned by Hermitian theta series associated to lattices of fixed level, divisor and volume over an arbitrary quadratic imaginary number field. We also prove the stability of the space spanned by theta series attached to a fixed genus of similar lattices, under a certain subalgebra of $\textup{GU}$ Hecke operators, the even subalgebra. Finally, we calculate explicitly the Eichler matrices in our setting for some Hecke operators.
\end{abstract}
\maketitle
\vspace{-1cm}
\section{Introduction}
One of the central themes in the theory of quadratic forms is the study of representation numbers. Questions concerning their multiplicative properties lead naturally to the study of the action of Hecke operators on spaces spanned by theta series. Eichler \cite{eichler_book} was one of the first to undertake a systematic study of the subject for degree-one theta series (classical modular forms), and to develop what became known as ``Eichler Commutation Relations''. It was not until the work of Freitag \cite{freitag1} that such questions were investigated in the setting of Siegel theta series of arbitrary degree. Freitag's method involved singular modular forms and allowed him to find closed spaces of theta series under the Hecke action that were generally quite large. Andrianov \cite{andrianov1}, using a different approach involving parabolic Hecke algebras, then deduced more precise and explicit results. In this paper, we study the action of $\textup{GU}$-Hecke operators on Hermitian theta series over arbitrary quadratic imaginary number fields and give explicit formulas for the action of certain operators in terms of representation numbers.\newline

Let $K$ be any imaginary quadratic number field of discriminant $-d_K$ and $(V,\tau)$ be a finite-dimensional positive definite Hermitian vector space over $K$, of dimension $m \geq 1$. For $\ell, \delta \in \mb{Z}_{>0}$ and $d \in \mb{Q}_{>0}$, we define the set 
\begin{multline*}
\Lambda(\ell, \delta,d):= \left\{[(\mathcal{L'}, \rho)] \mid (\mc{L}', \rho) \,\, \mathbb{Z}\textup{-integral } \mathcal{O}_K \textup{-lattice in }\, V, \right.\\\left. \rho \in \mathbb{Q}_{>0}\tau,\,\, \ell((\mathcal{L'}, \rho))= \ell,\,\, \delta((\mathcal{L'}, \rho)) = \delta,\,\, \textup{vol}((\mathcal{L'}, \rho)) = d\right\},
\end{multline*}
where $\ell(\mathcal{L}',\rho)$, $\delta(\mathcal{L}',\rho)$ and $\textup{vol}(\mathcal{L}',\rho)$ denote the level, divisor and volume of the Hermitian lattice $(\mathcal{L}',\rho)$ (see Section \ref{theta and epsilon series}), and the brackets denote isometry classes of lattices (see \eqref{isometric lattices}). This is a finite set (see Lemma \ref{lambda finite}). For a lattice $(\mathcal{L}',\rho)$ in the set above, we write $\theta_{\mb{A},n}(\cdot; (\mc{L',\rho}))$ for the degree $n \geq 1$ adelic theta series of weight $m$ and level $\textup{lcm}(\ell, \delta d_K)$ on $\un(n,n)_{\mb{A}}$ (see \eqref{adelic theta}).\newline

We introduce the complex vector space,
\[
\Theta^n(\ell,\delta,d) := \textup{span}_{\mb{C}}\{ \theta_{\mathbb{A},n}(\cdot;(\mathcal{L}',\rho))\mid[(\mathcal{L}',\rho)] \in \Lambda(\ell, \delta,d) \},
\]
and for a fixed lattice $(\mathcal{L},\tau) \in \Lambda(\ell,\delta,d)$ the subspace
\[
\Theta^n(\textup{Gen}(\mathcal{L},\tau)) := \textup{span}_{\mb{C}}\{\theta_{\mathbb{A},n}(\cdot;(\mathcal{L}',\rho))\mid [(\mathcal{L}',\rho)] \in \textup{Gen}(\mathcal{L},\tau) \},
\]
where $\textup{Gen}(\mc{L}, \tau) \subseteq \Lambda(\ell, \delta,d)$ consists of classes of lattices which are locally isometric to $(\mc{L}, \tau)$ for every finite prime $p$.\\

By $\mc{H}^{(n)}$ we denote the Hecke algebra of $\gu(n,n)_{\mb{A}}$, taken trivial on primes dividing the level $\textup{lcm}(\ell, \delta d_K)$. Let also $\mc{H}^{(n)+}$ denote the subalgebra of even Hecke operators, whose rational similitude factor belongs in $N_{K/\mb{Q}}(K^{\times})$. Then our main Theorem \ref{Theorem m neq n} gives
\begin{theorem}\label{theorem intro}
For any positive definite $\mb{Z}$-integral Hermitian $\mc{O}_K$-lattice $(\mc{L},\tau)$, the space $\Theta^n(\textup{Gen}(\mathcal{L},\tau))$ is $\mathcal{H}^{(n)+}$-stable. Moreover, if $m$ is even or $m < 2n +1$, then the space $\Theta^n(\ell,\delta,d)$ is $\mathcal{H}^{(n)}$-stable. If $m$ is odd and $m \geq 2n+1$, $\Theta^n(\ell,\delta,d)$ is stable under the subalgebra of $\mc{H}^{(n)}$ generated by $\mc{H}_p^{(n)}$ for split primes $p$ in $K$ and $\mc{H}_p^{(n)+}$ for inert primes $p$ in K, all away from the level.
\end{theorem}
This theorem is the Hermitian analogue of Andrianov's \cite[Theorem 5.1.15]{andrianov_book}. In particular, Andrianov \cite[Chapter 5]{andrianov_book} studies the action of Hecke operators of $\textup{GSp}_n(\mb{Q})$ on Siegel theta series using a technique of Hecke algebra embeddings in the Siegel parabolic subgroup. Here we develop his approach in the Hermitian setting, working adelically with Hecke algebras of $\textup{GU}(n,n)$ over an arbitrary quadratic imaginary number field $K$. A non-trivial class number makes the arithmetic of Hermitian lattices much more intricate; for example, they are not always free. This is the main reason an adelic setting is essential. In the class number one case for unimodular lattices, we mention the result \cite[Satz 2.2.3]{thesis}, which follows from Freitag's approach in the Siegel case \cite{freitag_book}. Our results here cover arbitrary level and class number. An adelic setting also allows us to work locally, where the Hecke algebras of unitary groups differ depending on the splitting behaviour of each rational prime $p$ in $K$. In particular, working away from the level, the case of inert primes $p$ is completely analogous, and in some sense even simpler (see Section \ref{Section:Inert primes}) than the symplectic setting of Andrianov. On the other hand, the case of good split primes is substantially different and requires new ideas. We develop the required machinery in terms of parabolic Hecke algebras in Sections \ref{Hecke rings of the general linear group} and \ref{Split-prime local unitary Hecke rings}.\\
% We mention here the result of \cite[Satz 2.2.3]{thesis}, following Freitag's method of singular modular forms in the hermitian setting. However, the results there are restricted to unimodular lattices and trivial class number.\\

%We remark that our proof of the above theorem allows us to handle also the case where $m<n$ (the singular theta series) but then we need the assumption that the class number of $K$ is one (see Remark \ref{the case m<n}). We hope to return to this question in the near future.\newline
There is also a different point of view on tackling questions of Hecke stability of spaces of theta series, first developed by Rallis in \cite{rallis1} for the classical case and \cite{rallis2} for general Siegel theta series. The main idea in these works is to establish maps between Hecke algebras of the pair of algebraic groups $(\textup{Sp}_n, \textup{O}_m)$ through their Satake transforms. Interpreting Siegel theta series as theta-lifts of a basis of algebraic modular forms on the definite orthogonal group yields stability results. An important distinction is that \cite{rallis2}, which deals with $n\geq 2$, considers only operators with trivial similitude. Analogous local results exist in the similitude-free unitary case, for the dual pair $(\un(n,n), \un(m))$, in the work of Liu \cite[Appendix]{liu} and Ichino \cite[Proposition 2.1]{ichino}. However, Liu assumes $m$ is even in the inert case, while Ichino's proposition covers only the case of inert primes and is restricted to the range $m\leq 2n+1$.\\

The method employed here yields explicit Hecke stability results in quite a general setting. Theorem \ref{theorem intro} covers all cases $m,n\geq 1$ with the Hecke subalgebras specified above and gives explicit formulas for operators with nontrivial similitude. In particular, we show that the space $\Theta^{n}(\textup{Gen}(\mc{L},\tau))$ is stable under the even Hecke subalgebra, which includes all the operators with similitude one. We also explicitly identify the sets $\Lambda(\ell,\delta, d)$ and $\textup{Gen}(\mc{L}, \tau)$ of classes of lattices whose associated theta series span these stable spaces.\\ %These sets give sharp enough bounds on which theta series are necessary. Notably, $\textup{Gen}(\mc{L}, \tau)$ now also includes similar lattices to the original $(\mc{L}, \tau)$.\\ 

% The adelic formulation allows us to work locally and initially prove the Hecke stability results for each prime away from the level. To do this, in Section \ref{theta and epsilon series} and, in particular, Proposition \ref{relation theta epsilon}, we establish a connection between our adelic theta series and a certain epsilon series, which essentially records primitive representation numbers and is, in fact, an automorphic form for the Siegel parabolic subgroup. In Section \ref{The Hecke action on epsilon series}, we then study the action of parabolic Hecke algebras on the epsilon series and prove a particular expansion in Proposition \ref{action of hecke on epsilon}. This expansion then allow us to specify the sets $\Lambda(\ell, \delta, d)$ and $\textup{Gen}(\mc{L},\tau)$ (see Proposition \ref{level divisor det}). 
Once the Theorem above is established, one may ask for an explicit description of the Hecke action in terms of representation numbers appearing in the Fourier coefficients of the theta series. We address this question in Section \ref{Eichler Matrices and Representation Numbers}. After fixing a spanning set for the space $\Theta^n(\ell,\delta,d)$ of theta series attached to lattices $(\mc{L}_i,\tau_i)$, $i=1,\ldots, H$, we associate to every $T \in \mc{H}^{(n)}$ a matrix $\Theta(T)$ describing its action (see \eqref{theta matrix}). This matrix is unique if $n \geq m$, since this spanning set is a basis (see Lemma \ref{linear independence}), but need not be unique if $n < m$. We then compute these Theta matrices for certain Hecke operators. In particular, for a prime $p$ away from the level, define
\begin{equation*}
    R_{ij}(p):= \#\left\{ u\in M_m(K)\mid \tau_i[u]=p\tau_j,\quad u\mc{L}_j\subseteq\mc{L}_i \right\}.
\end{equation*}
If $p\mc{O}_K = \mf{p}\mf{\ol{p}}$ splits, we also define, for $0\leq k \leq m$
\begin{equation*}
    R_{ij}^{(k)}(p) := \#\left\{u \in M_m(K)\mid  \tau_i[u] = p\tau_j, \,\, u\mc{L}_j \subseteq \mc{L}_i, \,\, \mc{L}_i/(u\mc{L}_j) \cong \left(\mc{O}_K/\mf{p}\right)^{k} \oplus \left(\mc{O}_K/\mf{\ol{p}}\right)^{m-k}  \right\},
\end{equation*}
We then define the \textit{Eichler matrices}
\begin{equation*}
    M(p):=(e_i^{-1}R_{ij}(p))_{i,j=1}^{H}, \qquad M_k(p):=(e_i^{-1}R_{ij}^{(k)}(p))_{i,j=1}^{H},
\end{equation*}
where $e_i:=|\textup{Aut}(\mc{L}_i,\tau_i)|$.

\begin{theorem}
Let $m,n\geq1$. For the standard split-prime operators $T_r(p)$, $0\leq r\leq2n$, and the inert-prime operator $T(p)$, we may take \begin{equation*}
    \Theta(T)=
    \begin{cases}
        \displaystyle \chi(\mf{p})^{n-r} \sum_{k=0}^{m}p^{-k(r-k+m-n)}{2n-m\brack r-k}_p M_k(p),&\quad
        \begin{array}{l}
            T=T_r(p),\\
            p\textup{ split},
        \end{array}
        \\[16pt]
        \displaystyle
        \frac{\prod_{j=0}^{n-1}(1+p^{2j+1-m})}
             {\prod_{j=1}^{m/2}(1+p^{2j-1})}M(p),
        &\quad
        \begin{array}{l}
            T=T(p),\\
            p\textup{ inert},\ m\textup{ even},
        \end{array}
        \\[16pt]
        0,
        &\quad
        \begin{array}{l}
            T=T(p),\\
            p\textup{ inert},\ m\textup{ odd},\ m<2n+1.
        \end{array}
    \end{cases}
\end{equation*}
Here $\chi$ is the character in \eqref{chi defn}, defined using the rank-$m$ lattices, and the Gaussian binomial coefficients are those defined in \eqref{extended binomial}. 
\end{theorem}

\begin{notation}
    We denote the space of $m\times n$ matrices with coefficients in a commutative ring $R$ by $M_{m,n}(R)$. If $n=m$, we often use the notation $M_{n}(R)$. 
    % By $1_n$, we denote the $n\times n$ identity matrix. For any matrix $M \in M_{n}(R)$, we denote by $\det(M), \textup{ tr}(M)$ the determinant and trace of $M$ respectively. By $\textup{GL}_n(R)$, we denote the matrices in $M_{n}(R)$ with determinant in $R^{\times}$, the units of $R$, and by $\textup{SL}_n(R)$ the matrices with determinant $1$. For any vector $v$, we denote by $v^{t}$ its transpose.
    We use the bracket notation $A[B] := \overline{B}^{t}AB$ for suitably sized complex matrices $A,B$. 
    % By $\textup{diag}(A_1, A_2, \cdots, A_n)$, we will denote the block diagonal matrix with the matrices $A_1, A_2,\cdots, A_n$ in the diagonal blocks. 
    For an algebraic group $G$ over $\mathbb{Q}$, and a field $F$ containing $\mathbb{Q}$, we will denote by $G_{F}$ the group of $F$-rational points of $G$. We will abbreviate $G_{\mathbb{Q}_p}$ by $G_{p}$. We also denote by $G_{\mathbb{A}}, G_{\mathbf{h}}, G_{\infty}$ the adelized group of $G$, the finite part of $G_\mathbb{A}$ and the infinite part of $G_{\mb{A}}$ respectively. For a $\mathbb{Z}$-module $L$, we denote by $L_p := L \otimes_{\mathbb{Z}} \mathbb{Z}_p$. For $z \in \mathbb{C}$, let $e(z) := \textup{exp}(2\pi i z)$.
\end{notation}
\textbf{Statement on the use of AI}. GPT-6.1 Sol was used in this work to assist with the computations of Proposition \ref{action of B} and Theorem \ref{eichler matrices split}, as well as with the adaptation of the arguments in \cite{andrianov_book} to the inert prime case presented in the Appendix. We also used it to identify relevant literature, check calculations, help detect errors, and proofread the manuscript. We independently verified all the calculations and arguments and take full responsibility for the contents of the paper.
\section{Preliminaries}\label{preliminaries}
Let $n \geq 1$. Let $K/\mathbb{Q}$ denote a quadratic imaginary number field. Write $K = \mathbb{Q}(\sqrt{-d_K})$, where we normalise $d_K$ such that $-d_K$ is the discriminant of $K$. For any matrix $A$ with entries in a ring with involution $\rho$, we employ the notation
\begin{equation*}
    A^{*} := (A^{\rho})^t, \qquad \widehat{A} := (A^{*})^{-1}.
\end{equation*}

Let
\begin{equation*}
    \mathbf{G} := \textup{GU}(n,n)_{\mathbb{Q}} = \{g \in \textup{GL}_{2n}(K) \mid g^{*} J_n g = \mu(g) J_n, \,\, \mu(g) \in \mathbb{Q}^{\times}\}, \qquad J_n := \m{0 & -1_{n} \\ 1_n & 0_n},
\end{equation*}
\begin{equation*}
     G = \un(n,n)_\mathbb{Q} = \{g \in \gu(n,n)_\mathbb{Q} \mid \mu(g)=1\},
\end{equation*}
\begin{equation*}
    \su(n,n)_\mathbb{Q} = \un(n,n)_\mathbb{Q} \cap \textup{SL}_{2n}(K),
\end{equation*}
denote the general unitary, unitary, and special unitary groups, respectively. For any $n \times n$ block decomposition of a $2n \times 2n$ matrix $x$, we write $$x = \m{A & B \\ C& D}$$ and set $a_{x}:=A, \, b_x := B,\, c_{x} := C,\, d_x := D$. Let
\begin{equation}\label{hermitian upper half plan}
    \mathbb{H}_n := \{Z \in M_n(\mathbb{C}) \mid -i(Z - Z^{*}) >0\}
\end{equation}
denote the Hermitian upper half-space. Then the group $\mathbf{G}_{\infty}^{+} := \{x \in \mathbf{G}_{\infty} \mid \mu(x)>0\}$ acts transitively on $\mathbb{H}_n$ via
\begin{equation*}
    g Z := (a_g Z + b_g)(c_g Z + d_g)^{-1}.
\end{equation*}
Now, for each $g \in \mathbf{G}_{\infty}^{+}$ and $Z \in \mathbb{H}_n$, we set 
\begin{equation*}
    j(g, Z) := \det(c_{g}Z+d),
\end{equation*}
and for $k\in \mathbb{N}$ and a function $f : \mathbb{H}_n \longrightarrow \mathbb{C}$, set
\begin{equation*}
    f \mid _{k} g := j(g, Z)^{-k}f(gZ).
\end{equation*}
We then have the following Definition.
\begin{defn}\label{hermitian modular form classically}
    Let $\Gamma \subset \mathbf{G}_{\mb Q}^{+}$ be a congruence subgroup in the sense of \cite[Section 4.10]{shimura2000arithmeticity} and $\chi : \Gamma \longrightarrow \mb{C}^{\times}$ a character. We say that a function $f: \mathbb{H}_{n} \longrightarrow \mathbb{C}$ is a Hermitian modular form of weight $k$, level $\Gamma$ and character $\chi$ if $f$ is holomorphic and $f\mid_{k} \gamma = \chi(\gamma)f$ for all $\gamma \in \Gamma$. If $n=1$, we also require that $f$ is bounded at infinity. The space of modular forms of weight $k$ and level $\Gamma$ will be denoted by $M(\Gamma, k, \chi)$. If $\chi$ is the trivial character, we omit it from the notation.
\end{defn}
\begin{remark}
    Note that actually $\Gamma \subset G_{\mathbb{Q}}$ as is explained right after \cite[Definition 3.1]{klosin}.
\end{remark}

For any two fractional ideal $\mf{b}, \mf{c}$ in $\mathbb{Q}$ such that $\mf{bc} \subseteq \mathbb{Z}$, and a finite prime $p$, we set
\begin{equation*}
    \mathbf{D}_p[\mf{b}, \mf{c}] := \left\{x \in \mathbf{G}_{p}\mid a_{x} \in M_n(\mathcal{O}_{K,p}), \ b_{x} \in M_n(\mf{b}_p), \, c_{x} \in M_n(\mf{c}_p), \, d_x \in M_n(\mathcal{O}_{K,p}), \,\,  \mu(x) \in \mathbb{Z}_p^{\times}\right\}.
\end{equation*}
We then define
\begin{equation*}
    \mathbf{D}[\mf{b}, \mf{c}] := \prod_{p \nmid \infty} \mbf{D}_{p}[\mf{b}, \mf{c}] \subseteq \mathbf{G}_{\mathbb{A}}, \qquad \mathbf{D} := \mathbf{D}[\mf{b}^{-1}, \mf{bc}].
\end{equation*}
We also define $\mathbf{D}_{\infty}$ to be the stabiliser of $\mathbf{i}_n := i 1_n$ in $\mathbf{G}^{+}_{\infty}$, given by 
\begin{equation*}
    \mathbf{D}_{\infty} = \{g \in \mathbf{G}_{\infty}^{+} \mid g \mathbf{i} = \mathbf{i}\}.
\end{equation*}
Let also $\chi$ denote a Hecke character of $K$ such that $\chi_p(a) = 1$ for all $a \in \mathcal{O}_{K,p}^{\times}$ such that $a-1 \in \mathcal{O}_{K,p}\mf{c}_p$. We then have the following definition for automorphic forms on $\mathbf{G}$. 
\begin{defn}\label{automorphic form on GU}
Let $k \in \mathbb{N}$. We denote by $\mathcal{N}(\mbf{D},k, \chi)$ the set of all functions 
$\mathbf{f}:\mathbf{G_{\mathbb{A}}} \longrightarrow \mathbb{C}$ which satisfy
\begin{enumerate}
\item $\mathbf{f}(\gamma gw) = \chi_{\mathfrak{c}}^{-1}(\det (a_w)) j(w_{\infty}, \mathbf{i}_n)^{-k}\mathbf{f}(g)$ for all $\gamma \in \mbf{G}_{\mb{Q}}$, $g \in \mbf{G}_{\mb A}$, $w \in \mathbf{D}_{\infty}\mbf{D}$,
%\item $\mathbf{f}(ag) = a^{-2n\nu -nk } \mathbf{f}(g)$ for $a \in \mathbb{C}^{\times}$,
\item For any $g \in \mathbf{G}_{\mathbf{h}}$, there exists a modular form $f_{g} \in M(\Gamma, k)$ for some congruence subgroup $\Gamma$, called the $g$-component of $\mathbf{f}$, such that 
\[
\mathbf{f}(gy) = (f_g \mid _{k}y)(\mathbf{i}_n),
\]
for all $y \in \mathbf{G}_{\infty}^{+}$. Note that $f_{g}$ is uniquely determined by $\mbf{f}$ and $g$.
\end{enumerate}
\end{defn}
Similarly, we define the space of automorphic forms on $G = \textup{U}(n,n)$. Let $D := \mathbf{D} \cap G$ and $D_{\infty} := \mbf{D}_{\infty} \cap G$. 
\begin{defn}\label{automorphic form on U}
Let $k \in \mathbb{N}$. We denote by $\mathcal{M}_{}(D, k,\chi)$ the set of all functions which satisfy conditions $(1), (2)$ of the above definition with $\mbf{G}, \mbf{D}, \mbf{D}_{\infty}, \mbf{G}_{\infty}^{+}$ replaced by $G, D, D_{\infty}, G_{\infty}$ respectively.
% $\mathbf{f}:G_{\mathbb{A}} \longrightarrow \mathbb{C}$ which satisfy
% \begin{enumerate}
% \item $\mathbf{f}(\gamma gw) = \textup{det}(w_{\infty})^{-\nu} \chi_{\mathfrak{c}}^{-1}(\det (a_w)) j(w_{\infty}, \mathbf{i}_n)^{-k}\mathbf{f}(g)$ for all $\gamma \in G_{\mb Q}$, $g \in G_{\mb A}$ and $w \in D_{\infty}D$,
% %\item $\mathbf{f}(ag) = a^{-2n\nu -nk } \mathbf{f}(g)$ for all $a \in S^{1} := \{z \in \mathbb{C}^{\times} \mid |z|=1\}$,
% \item For any $g \in G_{\mathbf{h}}$, there exists a modular form $f_{g} \in M(\Gamma, k, \nu)$ for some congruence subgroup $\Gamma$ (recall $\Gamma \subseteq G_{\mathbb{Q}}$), called the $g$-component of $\mathbf{f}$, such that 
% \[
% \mathbf{f}(gy) = (f_g \mid _{k,\nu}y)(\mathbf{i}_n),
% \]
% for all $y \in G_{\infty}$. Note that $f_{g}$ is uniquely determined by $\mbf{f}$ and $g$.
% \end{enumerate}
\end{defn}
It is now possible to identify the spaces $\mc{M}(D, k, \chi)$ and $\mc{N}(\mbf{D}, k, \chi)$.
\begin{lemma}
    The map $\mbf{f} \longmapsto \mbf{f} \mid _{G}$ gives an isomorphism $\mc{N}(\mbf{D},k, \chi) \cong \mc{M}(D, k, \chi)$. 
\end{lemma}
\begin{proof}
    The crucial condition that makes this Lemma true is that, with our choice of $\mbf{D}$, we have $\mu(\mbf{D}) = \widehat{\mb{Z}}^{\times}$. We first show that the map is injective. Let $g \in \mbf{G}_{\mb{A}}$ and let $\mu = \mu(g) \in \mathbb{Q}_{\mathbb{A}}^{\times}$. We may write uniquely $\mu = \alpha r \omega$, with $\alpha \in \mb{Q}^{\times}$, $r \in \mb{R}_{>0}$ and $\omega \in \wh{\mb{Z}}^{\times}$. We then let $\gamma := \textup{diag}(1_n, \alpha1_{n})$, $(z_{\infty}, z_{\mbf{h}}) = (\sqrt{r}1_{2n}, 1)$ and choose $d \in \mbf{D}$ such that $\mu(d) = \omega$. This is possible since $\mu(\mbf{D}) = \wh{\mb{Z}}^{\times}$. We then observe $\mu(\gamma^{-1}z^{-1}gd^{-1}) = 1$, so we may write $g = \gamma z u c$, with $u \in G_{\mb{A}}$. Hence, from Definition \ref{automorphic form on GU}, we have
    \begin{equation*}
        \mathbf{f}(g) = \mathbf{f}(\gamma z u d) = r^{-nk/2} \chi_{\mf{c}}^{-1}(\det(a_{d}))\mathbf{f}(u),
    \end{equation*}
    so $\mathbf{f} \mid_{G} = \mathbf{f}' \mid _{G} $ implies $\mbf{f} = \mbf{f'}$.\\
    
    We now show that the map is surjective. We need to show that every $\mbf{f} \in \mc{M}(D, k,\chi)$ has an extension to $\mbf{G}_{\mathbb{A}}$. For $g \in \mathbf{G}_{\mathbb{A}}$, we write it again as $g = \gamma z u d$ as above. We then define
    \begin{equation}\label{extension to GU}
        \mbf{f}(g) := r^{- nk/2} \chi_{\mf{c}}^{-1}(\det(a_{d}))\mathbf{f}(u).
    \end{equation}
    We will show this is well-defined. Indeed, if $g = \gamma' z' u'd'$ is another such decomposition, by considering the finite component and taking similitude, we obtain $\mu(\gamma d) = \mu(\gamma d')$, which implies $\mu(\gamma'\gamma^{-1}) \in \mathbb{Z}_p^{\times}$ for all primes $p$. By also considering the infinite component, we obtain $\mu(\gamma z) = \mu(\gamma' z')$, which also implies $\mu(\gamma'\gamma^{-1})>0$, and since $\mu(\gamma'\gamma^{-1}) \in \mb{Q}$, we obtain $\mu(\gamma) = \mu(\gamma')$. Hence, this also shows $r = r'$. Now, by looking at the finite component, we obtain $\mu(d) = \mu(d')$. Hence, we may write $(\gamma')^{-1}\gamma \in G$ and $d(d')^{-1} \in D$ with $D = \mbf{D} \cap G$. Hence
    \begin{equation*}
        u' = (\gamma')^{-1}\gamma u d (d')^{-1}
    \end{equation*}
    and therefore from Definition \ref{automorphic form on U},
    \begin{equation*}
        \mbf{f}(u') = \mbf{f}((\gamma')^{-1}\gamma u d (d')^{-1}) = \chi_{\mf{c}}(\det (a_d))\chi_{\mf{c}}^{-1}(\det (a_{d'}))\mbf{f}(u),
    \end{equation*}
    which gives $\mbf{f}(g) = \mbf{f}(g)$ using \eqref{extension to GU}.
\end{proof}
% \begin{remark}
%     The proposition is analogous to \cite[Proposition 3.10]{klosin}. There, it is assumed that $\wh{\mathcal{O}}^{\times}_{K} \subseteq \det(\mbf{D})$, which cannot be satisfied for any open compact subgroup $\mbf{D}$. The author needs this to assume the existence of a unitary basis, as in \cite[Proposition 3.8]{klosin}. In fact, this condition needs to be replaced with a similar one (see \cite[Lemma 8.14]{Shimura_Euler_Product}). In our case, and since $\mu(\mbf{D}) = \wh{\mb{Z}}^{\times}$, we can prove the correspondence without the need to assume a unitary basis. 
% \end{remark}
We now define the Siegel parabolic subgroups given by
\begin{equation*}
    \mbf{P} := \left\{x \in \mbf{G} \mid c_{x} = 0\right\}, \qquad P := \mbf{P} \cap G.
\end{equation*}
We also define the corresponding open compact subgroups
\begin{equation*}
    \mbf{D}_0 = \mbf{D} \cap \mbf{P}, \qquad D_0 = \mbf{D}_0 \cap G,
\end{equation*}
and the corresponding groups at infinity
\begin{equation*}
    \mbf{D}_{\infty,0} := \mbf{D}_{\infty} \cap \mbf{P}_{\infty}, \qquad D_{\infty,0} := D_{\infty} \cap P_{\infty}.
\end{equation*}
We need a notion of automorphic forms on the parabolic subgroups $P$ and $\mathbf{P}$. Similarly to Proposition \ref{hermitian modular form classically}, we have the following definition of ``triangular'' forms.
\begin{defn}\label{triangular}
    Let $\Gamma \subset \mathbf{P}_{\mathbb{Q}}^{+} := \mathbf{P}_{\mb Q} \cap \mathbf{G}_{\infty}^{+}$ be a congruence subgroup in the sense of \cite[Section 4.10]{shimura2000arithmeticity} and let $\chi: \Gamma \longrightarrow \mb{C}^{\times}$ be a character. We say that a function $f: \mathbb{H}_{n} \longrightarrow \mathbb{C}$ is a triangular form of weight $k$, level $\Gamma$ and character $\chi$ if $f$ is holomorphic and $f\mid_{k} \gamma = \chi(\gamma)f$ for all $\gamma \in \Gamma$. If $n=1$, we also require that $f$ is bounded at infinity. The space of triangular forms of weight $k$ and level $\Gamma$ will be denoted by $M_0(\Gamma, k, \chi)$. If $\chi$ is the trivial character, we omit it from the notation.
\end{defn}
We then have the following completely analogous Definitions as above.

\begin{defn} Let $k \in \mathbb{N}$. We denote by $\mathcal{N}_0(\mbf{D}_0,k,\chi)$ the set of all functions  $\mathbf{f}:\mathbf{P_{\mathbb{A}}} \longrightarrow \mathbb{C}$ which satisfy conditions $(1), (2)$ of Definition \ref{automorphic form on GU} with $\mbf{G}, \mbf{D}, \mbf{D}_{\infty}, \mbf{G}_{\infty}^{+}$ replaced by $\mbf{P}, \mbf{D}_0, \mbf{D}_{0,\infty}, \mbf{P}_{\infty}$ respectively.
\end{defn}
% \begin{enumerate}
% \item $\mathbf{f}(\gamma gw) = \textup{det}(w_{\infty})^{-\nu} \chi_{\mathfrak{c}}^{-1}(\det (a_w)) j(w_{\infty}, \mathbf{i}_n)^{-k}\mathbf{f}(g)$ for all $\gamma \in \mbf{P}_{\mb{Q}}$, $g \in \mbf{P}_{\mb A}$, $w \in \mathbf{D}_{0,\infty}\mbf{D}_0$,
% %\item $\mathbf{f}(ag) = a^{-2n\nu -nk } \mathbf{f}(g)$ for $a \in \mathbb{C}^{\times}$,
% \item For any $g \in \mathbf{P}_{\mathbf{h}}$, there exists a modular form $f_{g} \in M_0(\Gamma, k, \nu)$ for some congruence subgroup $\Gamma$, called the $g$-component of $\mathbf{f}$, such that 
% \[
% \mathbf{f}(gy) = (f_g \mid _{k,\nu}y)(\mathbf{i}_n),
% \]
% for all $y \in \mathbf{P}_{\infty}^{+}$. Note that $f_{g}$ is uniquely determined by $\mbf{f}$ and $g$.
% \end{enumerate}
% \end{defn}
\begin{defn} Let $k \in \mathbb{N}$. We denote by $\mathcal{M}_0(D_0,k, \chi)$ the set of all functions  $\mathbf{f}:P_{\mb A}\longrightarrow \mathbb{C}$ which satisfy the conditions of the above Definition with $\mbf{P}, \mbf{D}_0, \mbf{D}_{0,\infty}, \mbf{P}_{\infty}^{+}$ replaced by $P, D_0, D_{0,\infty}, P_{\infty}^{+}$ respectively.
% \begin{enumerate}
% \item $\mathbf{f}(\gamma gw) = \textup{det}(w_{\infty})^{-\nu} \chi_{\mathfrak{c}}^{-1}(\det (a_w)) j(w_{\infty}, \mathbf{i}_n)^{-k}\mathbf{f}(g)$ for all $\gamma \in  P_{\mb{Q}}$, $g \in P_{\mb A}$, $w \in D_{0,\infty}D_0$,
% %\item $\mathbf{f}(ag) = a^{-2n\nu -nk } \mathbf{f}(g)$ for $a \in \mathbb{C}^{\times}$,
% \item For any $g \in P_{\mathbf{h}}$, there exists a modular form $f_{g} \in M_0(\Gamma, k, \nu)$ for some congruence subgroup $\Gamma$, called the $g$-component of $\mathbf{f}$, such that 
% \[
% \mathbf{f}(gy) = (f_g \mid _{k,\nu}y)(\mathbf{i}_n),
% \]
% for all $y \in P_{\infty}^{+}$. Note that $f_{g}$ is uniquely determined by $\mbf{f}$ and $g$.
% \end{enumerate}
\end{defn}
Next, we will consider certain Hecke algebras and their actions on automorphic triangular forms. For each finite rational prime $p$, we define the Hecke rings
\begin{equation*}
    \mathcal{H}_p:=\mathcal{H}_p(\mathbf{D}_p, \mathbf{G}_p) := \textup{span}_{\mb Q}\{\mathbf{D}_pM \mathbf{D}_p \mid M \in \mathbf{G}_{p}\}.
\end{equation*}
\begin{equation*}
    \mathcal{H}_{0,p}:=\mathcal{H}_{0,p}(\mathbf{D}_{0,p}, \mathbf{P}_p) := \textup{span}_{\mb Q}\{\mathbf{D}_pM \mathbf{D}_p \mid M \in \mathbf{P}_{p}\}.
\end{equation*}
We then define
\begin{equation*}
    \mathcal{H} := \bigotimes'_{p\nmid \mf{bc}} \mathcal{H}_{p}, \qquad \mathcal{H}_0 := \bigotimes'_{p\nmid \mf{bc}} \mathcal{H}_{0,p}.
\end{equation*}
Following \cite[Proposition 3.1.6]{andrianov_book}, we have injective ring homomorphisms given by 
\begin{align}\label{imbedding}
    \nonumber\iota_p:\mc{H}_p &\hookrightarrow{} \mc{H}_{0,p}\\
    \sum_{}a_i \mbf{D}_pg_i &\longmapsto \sum a_i \mbf{D}_{0,p} g_i.
\end{align}
for all primes $p \nmid \mf{bc}$. Indeed, to apply \cite[Proposition 3.1.6]{andrianov_book}, it suffices
\begin{equation*}
\mbf{D}_{0,p} \subseteq \mbf{D}_p, \qquad
\mbf{G}_p \subseteq \mbf{D}_p \mbf{P}_p, \qquad 
\mbf{D}_p \cap \mbf{P}_p \mbf{P}_p^{-1} \subset \mbf{D}_{0,p}.
\end{equation*}
The first condition follows from $\mathbf D_{0,p}=\mathbf D_p\cap\mathbf P_p$. For the second condition, we will show that $\mathbf{G}_p = \mathbf{P}_p \mathbf{D}_p$ for all primes $p \nmid \mf{bc}$. Recall that we have $G_p = D_p P_p$ for all $p \nmid \mf{bc}$ (\cite[Proposition 7.2 (2) and Subsection 9.7]{Shimura_Euler_Product}). Now, for an element $g \in \mathbf{G}_p$, we may set $$g'_p :=  \textup{diag}(1_n, \mu(g_p)^{-1}1_n)g_p \in G_p.$$ We may now decompose $g_p' = h_p' d_p'$ with $h_p' \in P_p$ and $d_p' \in D_p$. We then set $$h_p := \textup{diag}(1_n, \mu(g_v)1_n)h'_p$$ and obtain $g_p = h_p d_p' \in \mathbf{P}_p \mathbf{D}_p$, since $D = \mathbf{D} \cap G$. Finally, since \(\mathbf P_p\) is a group,
\[
\mathbf D_p\cap\mathbf P_p\mathbf P_p^{-1}
=\mathbf D_p\cap\mathbf P_p
=\mathbf D_{0,p},
\]
so the third condition also holds.\\

We now have the following action of $\mc{H}_0$ on $\mathcal{N}_0(\mathbf{D}, k, \chi)$ (cf. \cite[(20.10)]{shimura2000arithmeticity}).
\begin{defn}\label{hecke action}
    Let $\mathbf{f} \in \mathcal{N}_0(\mathbf{D}_0,k,\chi)$. Assume $\mathbf{D}_0 \tau \mathbf{D}_0= \bigsqcup_{\eta \in Y} \mathbf{D}_0\eta$, where $Y$ is a finite subset of $\mathbf{P}_{\mathbf{h}}$. We then define $\mathbf{f} \mid \mathbf{D}_0\tau \mathbf{D}_0 : \mathbf{P}_{\mathbb{A}} \longrightarrow \mathbb{C}$ as follows:
    \begin{equation*}
        (\mathbf{f} \mid \mathbf{D}_0\tau \mathbf{D}_0)(x) = \sum_{\eta \in Y}\chi_{\mf{c}}^{-1}(\det(a_{\eta}))\mathbf{f}(x\eta^{-1}), \,\,\, x \in \mathbf{P}_{\mathbb{A}}.
    \end{equation*}
\end{defn}
\begin{lemma} \label{unitary basis}
Let $\mc{B} \subset  \textup{GL}_n(K)_{\mbf{h}}$ be a finite set such that $$\textup{GL}_{n}(K)_{\mb{A}} = \bigsqcup_{b \in \mc{B} }\textup{GL}_n(K) \,b \, \textup{GL}_n(\widehat{\mathcal{O}}_K) \textup{GL}_n(\mathbb{C}).$$ 

For each $b \in \mc{B}$, put $p_{b}:= \textup{diag}(b, \wh{b})$. We then have, for our choice of level, that 
\[
\mathbf{G}_{\mathbb{A}} =\bigsqcup_{b \in \mc{B}} \mathbf{G}_{\mathbb{Q}}\, p_b\, \textbf{G}_{\mathbb{R}}\mathbf{D}.
\]
\end{lemma}

\begin{proof}
First of all, observe that with our choice of $D$, we have $\det (D) \subseteq \{u/\overline{u} \mid u \in \mc{O}_K^{\times}\}$. By the proof of \cite[Lemma 13.1]{Ikeda_2008}, we know that 
\[
G_{\mathbb{A}} =\bigsqcup_{b \in \mc{B}} G_{\mathbb{Q}}\, p_b\, G_{\mathbb{R}}D.
\]
Let $\mathbf{g} \in \mathbf{G}_{\mathbb{A}}$, and write $\mu:=\mu(\mathbf{g}) = \gamma_1 \alpha_1 \omega_1$, with $\gamma_1 \in \mathbb{R}^{\times}$, $\alpha_1 \in \mathbb{Q}^{\times}$ and $\omega_1 \in \widehat{\mathbb{Z}}^{\times}$. We then consider the element $$g:=\m{\alpha_1^{-1} \,1_n & 0 \\ 0 & 1_n} \mathbf{g} \m{\gamma_1^{-1} 1_n & 0 \\ 0 & 1_n} \m{\omega_1^{-1} 1_n & 0\\ 0 & 1_n}.$$ This element belongs to $G_{\mathbb{A}}$ and so we may write
\[
g = \alpha p_b \gamma \omega, 
\]
for some $b$ and with $\gamma \in G_{\mathbb{R}}$, $\alpha \in G_{\mathbb{Q}}$ and $\omega \in D$. But then we have
\[
\mathbf{g} = \m{\alpha_1 \,1_n & 0 \\ 0 & 1_n}\alpha p_b \gamma \m{\gamma_1 1_n & 0 \\ 0 & 1_n} \omega\m{\omega_1 1_n & 0\\ 0 & 1_n},
\]
which establishes the Lemma. The disjointness of the decomposition follows after considering similitudes.
\end{proof}
\section{Theta and epsilon series}\label{theta and epsilon series}
In this Section, we define theta and epsilon series, which are the main objects of study in this paper, and prove a connection between them.\\

Let $m \geq 1$ and take $V$ an $m$-dimensional vector space over $K$. Note that $V$ can be viewed as a $2m$-dimensional vector space over $\mathbb{Q}$. By a $\mathbb{Z}$-lattice $\mc{L}$ in $V$, we mean a finitely generated $\mathbb{Z}$-submodule of $V$ which spans $V$ over $\mathbb{Q}$. By a positive definite Hermitian $\mathbb{Z}$-lattice, we mean a pair $(\mc{L}, \tau)$ where $\mc{L}$ is a $\mathbb{Z}$-lattice in $V$ and $\tau$ is a positive definite Hermitian form on $V$. We call $(\mc{L}, \tau)$ a Hermitian $\mathcal{O}_K$-lattice, if $\mathcal{O}_K\mc{L} \subseteq \mc{L}$ (\cite[Section 2.10]{shimura_unitary}). We also denote by $\mf{n}(\mc{L}, \tau)$ the $\mathbb{Z}$-ideal generated by $\{\tau[x] \mid x \in \mc{L}\}$.  We call a $\mathbb{Z}$-lattice integral if $\mf{n}(\mc{L}, \tau) \subseteq \mathbb{Z}$. If it is clear what the Hermitian form is, we will just write $\mc{L}$ instead of $(\mc{L}, \tau)$. We have the following definitions.
\begin{defn}[Duals]
    We define the trace dual of a Hermitian $\mc{O}_K$-lattice $(\mc{L},\tau)$ as
    \begin{equation*}
        \widehat{\mc{L}} := \{x \in V \mid \tau(x, \mc{L}) \subseteq \mathfrak{d}^{-1}\},
    \end{equation*}
    where $\mathfrak{d}^{-1} = \mathfrak{d}^{-1}_{K/\mathbb{Q}} := \{x \in K \mid \textup{tr}_{K/\mathbb{Q}}(xy) \in \mathbb{Z}, \,\, \forall y\in \mathcal{O}_K\}$. We also define the usual Hermitian dual of $(\mc{L},\tau)$ as
    \begin{equation*}
        \mc{L}^{\#} := \{x \in V \mid \tau(x, \mc{L}) \subseteq \mathcal{O}_K\}.
    \end{equation*}
\end{defn}
\begin{defn}[Level]
    We define the level of $(\mc{L},\tau)$ to be
    \begin{equation*}
        \ell = \ell(\mc{L},\tau) := \textup{min}\left\{N \in \mathbb{Z}_{>0} \mid N\tau[x] \in \mathbb{Z},\,\, \forall x \in \widehat{\mc{L}}\right\}.
    \end{equation*}
    This definition agrees with the usual definition of the level for the associated quadratic form $q_{\tau}[x]:= \tau[x]$ for all $x \in V$, where now $V$ is viewed as a quadratic space of rank $2n$ over $\mathbb{Q}$.
\end{defn}
\begin{defn}[Divisor]
    For an integral lattice $(\mc{L},\tau)$, we define the divisor to be the positive integer
    \begin{equation*}
        \delta = \delta(\mc{L},\tau) := \gcd\{\tau[x] \mid x \in \mc{L}\}.
    \end{equation*}
\end{defn}

We note that the level and the divisor of $L$ are the same as those defined when $L$ is viewed as a lattice in the induced quadratic space. In particular we have that  $\delta | \ell$ (see \cite[page 349]{andrianov_book}).\\

For any Hermitian $\mc{O}_K$-lattice $(\mathcal{L},\tau)$ we know (see \cite[Theorem 81:3 and 81:5]{omeara}) that there exist fractional ideals $\{\mathfrak{a}_1,\ldots, \mathfrak{a}_n\}$, and a basis $\{x_1,\ldots, x_n\}$ of $V$ such that
\begin{equation*}
    \mathcal{L} = \sum_{i=1}^n\mathfrak{a}_i x_i.
\end{equation*}
The ideal class of $\mathfrak{a}:=\mathfrak{a}_1\cdots\mathfrak{a}_n$ depends only on the isomorphism class of $\mathcal{L}$ as an $\mc{O}_K$-module, and it's called its Steinitz class. We may select the basis $\{x_i\}$ such that $\mathcal{L} = \mathcal{O}_K x_1 + \ldots + \mathcal{O}_K x_{n-1} + \mathfrak{a} x_n$.  

\begin{defn}[Volume]\label{volume of lattice}
The volume of $(\mathcal{L},\tau)$ is defined as
\begin{equation}\label{volume defn}
    \textup{vol}(\mathcal{L},\tau):=\left(\det(\tau(x_i,x_j))\prod_{i=1}^{n}N(\mathfrak{a}_i)\right) \mathbb{Z} = \textup{det}(\tau(x_i,x_j))N_{K/\mathbb{Q}}(\mathfrak{a})\mb Z,
\end{equation}
a fractional $\mathbb{Z}$-ideal, which we will identify with its positive rational generator. 
\end{defn}

It is known that $\textup{vol}(\mathcal{L},\tau)$ is independent of the choice of the ideals $\{\mathfrak{a}_i\}$ and the choice of the basis $\{x_i\}$. Also, if $\mc{L} = \xi \mc{L}_0$, we know \cite[p. 348]{shimura_diophantine_2006} that $\mc{L}$ is isomorphic as a $\mathcal{O}_K$-module to the direct sum of $\mathcal{O}_K^{n-1}$ and $\det(\xi)\mathcal{O}_K$. Hence, we also have 
\begin{equation*}
     \textup{vol}(\mathcal{L},\tau) = \det (\tau)N_{K/\mb {Q}}(\det (\xi ))\mathbb{Z},
\end{equation*}
where by slight abuse of notation, we view $\det (\xi)$ as the corresponding fractional $\mathcal{O}_K$-ideal.\newline

We now fix, once and for all, a positive-definite $\mathbb{Z}$-integral Hermitian $\mathcal{O}_K$-lattice $(\mc{L},\tau)$ on $V$. By a slight abuse of notation, we will also denote by $\tau$ the Gram matrix of the Hermitian form $\tau$ with respect to a fixed basis $\{e_1, \ldots, e_{m}\}$ of $V$ over $K$. In particular, we will write $\tau^{*} = \tau \in \textup{GL}_m(K)$, but we will also view $\tau$ as a Hermitian form $V \times V \longrightarrow K$. It is also well-known that there exists some $\xi \in \textup{GL}_m(K)_{\mathbb{A}}$ such that $\mc{L} = \xi \mc{L}_0$, where $\mc{L}_0:= \mathcal{O}_K^m$ denotes the standard lattice in $V$. For $n \geq 1$, we now define
\begin{equation}\label{R and E}
    W := \prod_{p \nmid \infty}M_{m,n}(\mathcal{O}_{K,p}),\qquad E := \prod_{p \nmid \infty} \textup{GL}_{n}(\mathcal{O}_{K,p}),
\end{equation}
We define the functions $\lambda$ and $\lambda'$ on $W, E$ respectively, locally given by 
\begin{equation*}
    \lambda_{p}(x_p) = \begin{cases}1 \,\,\textup{ if } x_p \in W_p,\\0 \,\,\textup{ otherwise}\end{cases}, \,\,\, \lambda_{p}'(x_p) = \begin{cases}1 \,\,\textup{ if } x_p \in E_p,\\0 \,\,\textup{ otherwise}\end{cases}.
\end{equation*}
We now have the following definition of the classical theta series as in \cite[(A5.3)]{shimura2000arithmeticity}:
\begin{defn}
    We define the classical theta series attached to $\mc{L}$ as follows:
    \begin{equation}\label{classical theta}
        \theta(z; \mc{L}) := \sum_{u \in M_{m,n}(K)}\lambda (\xi^{-1}u)e(\textup{tr}(u^{*}\tau uz)), \,\, z \in \mathbb{H}_n.
    \end{equation}
    In the case $m=n$, we also define for $i=1, \cdots, n$, 
    \begin{equation*}
        \theta^{(i)}(z; \mc{L}) := \sum_{\substack{u \in M_{n}(K)\\\textup{rank } u = i}}\lambda (\xi^{-1}u)e(\textup{tr}(u^{*}\tau uz)), \,\, z \in \mathbb{H}_n.
    \end{equation*}

\end{defn}
Using now \cite[(A5.6)]{shimura2000arithmeticity}, we may define a $\theta_{\mathbb{A}}(x, \mc{L})$ on $G_{\mathbb{A}}$ via the rule
\begin{equation}\label{adelic theta}
    \theta_{\mb{A}}(x, \mc{L}) := j(x, \mbf{i})^{-m}\theta(x(\mbf{i}), \prescript{x}{}{\mc{L}}), \qquad x \in G_{\mb A}
\end{equation}
where the notation $\prescript{x}{}{\mc{L}}$ means we are taking the classical theta series \eqref{classical theta}, with $\prescript{x}{}{\lambda}$ in place of $\lambda$, where $\prescript{x}{}{\lambda}$ denotes the function of \cite[Theorem A5.4]{shimura2000arithmeticity}. We then have the following Lemma, due to Shimura.
\begin{lemma}\label{congruence subgroup theta}
    $\theta_{\mathbb{A}}(x, \mc{L})$ is an element of $\mathcal{M}(D,m, \varphi^{-m})$, where $D = D[\mathfrak{b}^{-1},\mathfrak{bc}]$, with 
        \begin{equation*}
            (\mathfrak{b},\mathfrak{bc}) = (\delta \mathcal{O}_K, \textup{lcm}(\ell, \delta d_K) \mathcal{O}_K).
        \end{equation*}
        the characteristics of the lattice $(\mc{L}, \tau)$ and $\varphi$ the Hecke character of $K$ which extends the quadratic Hecke character of $\mathbb{Q}$ corresponding to $K/\mathbb{Q}$ (see \cite[A5.3]{shimura2000arithmeticity}).
\end{lemma}
\begin{proof}
We apply \cite[Proposition A7.9]{Shimura_Euler_Product} with $H=2\tau$ and $\mc{M}:= \mc{L}^n \subset M_{m,n}(K)$. Then we consider 
\[
\mc{M}' := \{ \eta \in M_{m,n}(K) \mid \textup{tr}(\xi^* \tau \eta) \in \mathfrak{d}_{K}^{-1} \,\,\,\forall \xi \in \mc{M}\}. 
\]
Now we claim that, for our choice of $\mc{M}$, we have $\mc{M}' = \widehat{\mc{L}}^n \in M_{m,n}(K)$. Indeed we may write $\xi = (\xi_1,\ldots \xi_n)$ for $\xi_i \in \mc{L}$ (column vectors) and similarly $\eta = (\eta_1,\ldots \eta_n)$ with $\eta_i \in K^m$. Then we have that $\xi^* \tau \eta = (\xi_i \tau \eta_j)_{ij}$, and hence $\textup{tr}(\xi^* \tau \eta) = \sum_{i=1}^n \xi_i\tau \eta_i.$ That is, by selecting $\xi_j=0$ for all $j$ except an index $i$, we have that the condition on $\mc{M}'$ reads that $\eta_i \in \widehat{\mc{L}}$ for all $i$, which concludes the claim.\\

For $s = s^{*} \in M_n(\mc{O}_K)$, we have $\textup{tr}(\tau[u]s) \in \delta \mathbb{Z}$ for $u \in \mc{L}^n$ and $\textup{tr}(\tau[u]s) \in \ell^{-1} \mathbb{Z}$ for $u \in \mc{M}'$. Hence, with the notation of Shimura, we may take $\mf{x} = \delta \mathbb{Z}$, $\mf{y} = \ell^{-1}\mathbb{Z}$ and $\mathfrak{z} = \mathcal{O}_K$, since $\lambda$ is right-invariant under the action of $\textup{GL}_n(\wh{\mc{O}}_K)$. This gives the required level for the action $\textup{SU}(n,n)$. The extension to $\textup{U}(n,n)$, as well as the character $\vphi^{-m}$ can be deduced from \cite[Theorem A5.4]{shimura2000arithmeticity}.
% That is, the element $\widetilde{p}$ in the proof of \cite[Prop. A7.16]{Shimura_Euler_Product} is defined as
% $\widehat{\mc{L}}^n = \tau^{-1}\widetilde{p}M_{m,n}(\mathfrak{d}^{-1}_K)$, or equivalently $\widehat{\mc{L}} = \tau^{-1} \widetilde{p}(\mathfrak{d}^{-1}_K)^m$. In Shimura's notation the ideal $\mathfrak{t}$ is defined as $h^* \tau^{-1} h \in \mathfrak{t}^{-1}$ for all $h \in \widetilde{p}L_0$. We now claim that $N(\mathfrak{d}_K)\mathfrak{t} = \ell \mathcal{O}_K$. If we write $\mathfrak{d}_K = d \mathcal{O}_K$, then we have that $\widetilde{p}L_0 = d \tau \widehat{\mc{L}}$, and hence the condition reads $N_{K/\mathbb{Q}}(d) x^* \tau x \in \mathfrak{t}^{-1}$ for all $x \in \widehat{\mc{L}}$. That is, $\ell \mathcal{O}_K = \mathfrak{t} N(\mathfrak{d}_K)$.
\end{proof}
Let now
\begin{equation}\label{hermitian matrices in K}
    S := \{s \in M_n(K) \mid s^{*} = s\}
\end{equation}
For $q \in \textup{GL}_n(K)_{\mathbb{A}}$ and $s \in S_{\mathbb{A}}$, we use the notation $(q,s) := \m{q & s\widehat{q}\\ 0 & \widehat{q}}$. For such $q, s$ and $u \in M_{m,n}(K)$, we also set
\begin{equation}
    \Phi(u, q, s):= e(\textup{tr}(i\cdot q^*_{\infty}u_{\infty}^* \tau u_{\infty} q_{\infty} ))e_{\mathbb{A}}(\textup{tr}(u^* \tau u s)),
\end{equation}
\begin{equation}\label{betadef}
    \beta(q):=\overline{\textup{det}(q_{\infty})}^{m} |\textup{det}(q_{\mathbf{h}})|_K^{m/2} \varphi_{\mathbf{h}}(\det(q))^m,
\end{equation}
where for $x \in K_{\mathbb{A}}^{\times}$, we define
\begin{equation*}
    |x|_{K} := \prod_{v \in \mathbf{b} \cup \mathbf{k}} |x_v|_{v},
\end{equation*}
where $\mathbf{b}$ and $\mathbf{k}$ is the set of infinite and finite places of $K$, respectively.

\begin{proposition}\label{adelic expansion of theta}
    We have
    \begin{equation*}
        \theta_{\mathbb{A}}((q, s); (\mc{L},\tau)) = \beta(q)\sum_{u \in M_{m,n}(K)} \lambda(\xi^{-1}uq)\Phi(u, q, s).
    \end{equation*}
\end{proposition}

\begin{proof} This is stated in \cite[(A7.14.1)]{Shimura_Euler_Product} for $\textup{SU}(n,n)$. We now explain how this can be extended to $\textup{U}(n,n)$ and also compute the constant $\beta(q)$. We recall that 
\[
 \theta(z; \mc{L}) := \sum_{u \in M_{m,n}(K)}\lambda (\xi^{-1}u)e(\textup{tr}(u^{*}\tau uz)), \,\, z \in \mathbb{H}_n,
  \]  
  and 
  \[
\theta_{\mathbb{A}}(x;\mc{L})= j(x, \mathbf{i})^{-m} \theta(x(\mathbf{i}); \prescript{x}{}{\mc{L}}).
\]
We now set $x := \m{q & s\widehat{q} \\ 0 & \widehat{q}}$ in the above equation. We note that $$j_{x}(\mathbf{i})^{-m} = \overline{\textup{det}(q_{\infty})}^m, \qquad x(\mathbf{i}) = (q_{\infty}\mathbf{i} + s_{\infty}\widehat{q_{\infty}}) q_{\infty}^{*} = s_{\infty} + i q_{\infty}q_{\infty}^* =: z.$$
Moreover, by \cite[Theorem A5.4 (5)]{shimura2000arithmeticity} we have that $\prescript{x}{}{\lambda}$ depends only on the finite part of $x$, and there it is given as (we set $H= 2 \tau$ with the notation there),
\[
\prescript{x_{\mathbf{h}}}{}{\lambda}(v) = |\textup{det}(q_{\mathbf{h}})|_K^{m/2} \varphi_{\mathbf{h}}(\det(q))^m e_{\mathbf{h}}( \textup{tr}(v^* \tau v s )) \lambda(v q).
\]
That is,
\begin{multline*}
\theta_{\mathbb{A}}((q,s);\mc{L})= \overline{\textup{det}(q_{\infty})}^{m} |\textup{det}(q_{\mathbf{h}})|_K^{m/2} \varphi_{\mathbf{h}}(\det(q))^m \\\times\sum_{u \in M_{m,n}(K)} \lambda(\xi^{-1}u q) e_{\mathbf{h}}( \textup{tr}(u^*  \tau  u s )) e(\textup{tr}(u^{*}\tau u(s_{\infty} + i q_{\infty}q_{\infty}^*))).\qedhere
\end{multline*}
\end{proof}
From now on, and until Section \ref{m neq n} we are assuming that $m=n$.\newline

We now define, similar to \cite[(5.1.8)]{andrianov_book}, the associated $\epsilon$ series, which essentially amounts to considering ``primitive'' representations. We consider the series
\[
 \epsilon(z; \mc{L}) := \sum_{u \in \textup{GL}_n(K)}\lambda' (\xi^{-1}u)e(\textup{tr}(u^{*}\tau uz)), \qquad z \in \mathbb{H}_n.
  \]  
  We first establish that this is a triangular modular form for the group $\Gamma_0:= P_{\mathbb{Q}} \cap D[\mathfrak{b}^{-1},\mathfrak{bc}]$. First, the series converges on $\mb{H}_n$ by comparison with the theta series of $(\mc{L},\tau)$. For an element $\gamma = \m{ a & s\widehat{a} \\ 0 & \widehat{a}} \in \Gamma_0$ we have,  
\begin{equation*}
(\epsilon \mid_n \gamma)(z, \mc{L}) = \overline{\textup{det}(a)}^n \epsilon(aza^*+s,\mathcal{L}) = \overline{\textup{det}(a)}^n\sum_{u \in \textup{GL}_n(K)}\lambda' (\xi^{-1}u)e(\textup{tr}(a^*u^{*}\tau uaz))e(\textup{tr}(u^{*}\tau u s)).
\end{equation*}
By the definition of $\delta$ and the fact that $\mf{b} = \delta \mc{O}_K$, as well as the fact that $s = s^{*} \in M_n(\mf{b}^{-1})$, we have that the factor $e(\textup{tr}(u^{*}\tau u s))$ is $1$. Additionally, since $a \in \textup{GL}_n(\mc{O}_K)$, we have $\lambda'(\xi^{-1}ua) = \lambda'(\xi^{-1}u)$. Replacing $u$ by $ua$ in the sum gives
\begin{equation*}
    (\epsilon \mid_n \gamma)(z, \mc{L}) = \overline{\textup{det}(a)}^n\epsilon(z,\mathcal{L}).
\end{equation*}

But  $\det(a) \in \mc{O}_K^{\times}$, so the product formula for $\vphi$, together with the fact that $\vphi_{\infty}(t) = t^{-1}|t|$, gives $\overline{\textup{det}(a)}^n = \varphi_{\mf{c}}(\textup{det}(a))^{-n}$, and so $\epsilon(z,\mathcal{L})$ is triangular. \\

For $x \in P_{\mathbb{A}}$, we then define,
  \[
\epsilon_{\mathbb{A}}(x;(\mc{L},\tau)):= j(x, \mathbf{i})^{-n} \epsilon(x(\mathbf{i}); \prescript{x}{}{\mc{L}'}),
\]
where the notation this time means that we take the epsilon series with $\prescript{x}{}\lambda'$ in place of $\lambda'$, where $\prescript{x}{}\lambda'$ is as in \cite[Theorem A5.4]{shimura2000arithmeticity}. %We then, as we did with $\theta_{\mathbb{A}}$ we obtain
\begin{lemma}
    For the adelic epsilon series attached to $(\mc{L},\tau)$ we have:
    \begin{equation*}
        \epsilon_{\mathbb{A}}((q, s); (\mc{L},\tau)) = \beta(q)\sum_{u \in \textup{GL}_{n}(K)} \lambda'(\xi^{-1}uq)\Phi(u, q, s).
    \end{equation*}
Moreover $\epsilon_{\mathbb{A}}((q,s);(\mathcal{L},\tau)) \in \mathcal{M}_0(D_0,k,\varphi^{-n})$. 
    
\end{lemma}

\begin{proof}

This follows similarly to what we did for $\theta_{\mathbb{A}}$ in Proposition \ref{adelic expansion of theta}. We only remark that we can still apply \cite[Theorem A5.4]{shimura2000arithmeticity} even though the function $\lambda'$ is not a Schwartz function, since we take $x \in P_{\mathbb{A}}$ and hence the reciprocity law of the action of the finite adeles on functions of $M_n(K)_{\mathbb{A},\mathbf{h}}$ does not involve integration. The last statement follows from the fact that $D_0$ is the stabiliser of $\lambda'$ in $P_{\mathbb{A}}$.
%  We check the invariance properties. Let $(d,s_1) \in D_0$. Then $\epsilon_{\mathbb{A}}(q,s)(d,s_1) = \epsilon_{\mathbb{A}}(qd, qs_1q^* +s)$. That is
   % \[
 %  \epsilon_{\mathbb{A}}((q, s)(d,s_1); (\mc{L},\tau)) := \beta(qd)\sum_{u \in \textup{GL}_{n}(K)} \lambda'(\xi^{-1}uqd)\Phi(u, qd, qs_1q^*+s).
% \]
%However since $d \in E$ we have $\lambda'(\xi^{-1}uqd)=\lambda'(\xi^{-1}uq)$, and $\Phi(u,qd,qs_1q^*+s) = \Phi(u,q,s)$ as $d_{\infty}=1$ and $e_{\mathbb{A}}(tr((uq)^*\tau uq s_1)=1$ since $s_1 =s_1^*\in \prod_{p \nmid \infty}M_n(\mathfrak{b}^{-1}_{p})$ and hence in the dual of the matrix $(uq)^*\tau uq$.
\end{proof}
We now seek to obtain a relation between theta and epsilon series. To that end, let $0 \neq \mf{a} \subseteq \mathcal{O}_K$ be an ideal of $K$. We take an element $\widehat{\mf{a}} \in \prod_{v \textup{ finite}} K_{v}^{\times}$, such that $\widehat{\mf{a}}_v\mathcal{O}_{K_v} = \mf{a}_v$ for each finite place $v$ of $K$. Let also 
\begin{equation}\label{B set}
   B:= \textup{GL}_n(K)_{\mbf{h}} \cap W,
\end{equation}
where $W$ as in \eqref{R and E}. We define the following operator:
\begin{equation}\label{t minus original}
    t^{-}(\mf{a}) := \sum_{\substack{D \in E\backslash B/E\\ (\det D) \widehat{\mathcal{O}}_K=\widehat{\mf{a}}\widehat{\mathcal{O}}_K}} \mbf{D}_0 \m{D & \\ & \widehat{D}}\mbf{D}_0,
\end{equation}
For our later use we define the operator
\begin{equation}\label{t plus original}
    t^{+}(\mf{a}) := \sum_{\substack{D \in E\backslash B/E\\ (\det D) \widehat{\mathcal{O}}_K=\widehat{\mf{a}}\widehat{\mathcal{O}}_K}} \mbf{D}_0 \m{\widehat{D} & \\ & D}\mbf{D}_0.
\end{equation}
We then have the following Lemma.
\begin{lemma}\label{decomposition t minus}
    Let $q \in B$. We have the following decomposition.
    \begin{equation*}
        \mbf{D}_{0} \m{q & 0\\ 0 &\widehat{q}}\mbf{D}_{0} = \bigsqcup_{d \in E \backslash E qE } \mbf{D}_{0}\m{d & 0\\ 0 & \widehat{d}}. 
    \end{equation*}
    Therefore,
    \begin{equation*}
        t^{-}(\mf{a}) = \sum_{\substack{D \in E \backslash B\\(\det D) \wh{\mc{O}}_K=\wh{\mf{a}}\wh{\mc{O}}_K}} \mbf{D}_{0}\m{D & 0\\ 0 & \widehat{D}}. 
    \end{equation*}
\end{lemma}
\begin{proof}
    We first observe
    \begin{multline*}
        \m{\alpha_1 & \beta_1 \\ 0 & \mu_1\widehat{\alpha}_1} \m{q & 0\\ 0 &\widehat{q}} \m{\alpha_2 & \beta_2 \\ 0 & \mu_2\widehat{\alpha}_2} = \m{\alpha_1 q \alpha_2 & \alpha_1 q \beta_2 + \mu_2\beta_1 \widehat{q}\widehat{\alpha}_2\\ 0 & \mu_1\mu_2\widehat{\alpha_1q\alpha_2}} \\= \m{1 & x \\ 0 & \mu_1\mu_2} \m{\alpha_1 q \alpha_2 & 0\\ 0 & \widehat{\alpha_1q\alpha_2}},
    \end{multline*}
    where
    \begin{equation*}
        x = (\alpha_1 q \beta_2 + \mu_2\beta_1 \widehat{q}\widehat{\alpha}_2)(\widehat{\alpha_1q\alpha_2})^{-1} = \alpha_1q\beta_2\alpha_2^{*}q^{*}\alpha_1^{*} + \mu_2\beta_1\alpha_1^{*}.
    \end{equation*}
    Then $x^{*}=x$ and $x_p \in M_n(\mf{b}_p^{-1})$ for all primes $p$, since $\beta_{i,p} \in M_n(\mf{b^{-1}_p})$, $\mu_2 \in \mathbb{Z}_p^{\times}$ and $\alpha_{i,p}, q \in M_n(\mc{O}_{K,p})$ for $i=1,2$ and all primes $p$. Additionally, if 
    \begin{equation*}
        \m{d_1 & 0 \\ 0 & \widehat{d_1}} = \m{\alpha & \beta\\ 0 & \mu\widehat{\alpha}}\m{d_2 & 0 \\ 0 & \widehat{d_2}},
    \end{equation*}
    then $d_1 = \alpha d_2$, $\wh{d}_1 = \mu \wh{\alpha} \wh{d}_2$ and $\beta = 0$. This implies $\mu=1$ and therefore any two such left cosets are the same if and only if $d_i$ are the same modulo $E$.
\end{proof}
\begin{lemma}\label{action of t minus}
    If $F \in \mathcal{N}_{0}(\mbf{D}_0, k, \chi)$ is a triangular form, we have
    \begin{multline*}
        (F \mid t^{-}(\mathfrak{a})) \left((q,s)\right) = \sum_{\substack{D \in E \backslash B\\\det D = \mathfrak{a}}} \chi_{\mf{c}}^{-1}(\det (D))F\left(\m{qD^{-1}  & s\cdot  \widehat{qD^{-1}} \\ 0 & \widehat{qD^{-1}}}\right) \\= \sum_{\substack{D \in E \backslash B\\\det D = \mathfrak{a}}}\chi_{\mf{c}}^{-1}(\det (D)) F((qD^{-1},s)).
    \end{multline*}
\end{lemma}
\begin{proof}
    This follows directly from Lemma \ref{decomposition t minus} and Definition \ref{hecke action}. Note that the similitude of $\textup{diag}(D, \wh{D})$ is one, so the same Lemma is true if $F \in \mc{M}_0(D, k, \chi)$. 
\end{proof}

Our next step is to establish a relation between $\theta_{\mathbb{A}}^{(n)}$ and  $\epsilon_{\mathbb{A}}$. This may be seen as the Hermitian analogue of Andrianov's \cite[(5.1.12)]{andrianov_book}.
\begin{proposition}\label{relation theta epsilon}
    Let $\psi := \vphi_{\mbf{h}}\vphi_{\mf{c}}^{-1}$. We then have
    \begin{equation*}
        \theta_{\mathbb{A}}^{(n)}((q,s); (\mathcal{L}, \tau)) = \sum_{0 \neq \mathfrak{a} \subseteq \mathcal{O}_K}  \psi(\mathfrak{a})^{n}N_{K/\mathbb{Q}}(\mathfrak{a})^{-n/2}(\epsilon_{\mathbb{A}} \mid t^{-}(\mathfrak{a})) ((q,s); (\mathcal{L}, \tau)).
    \end{equation*}
\end{proposition}
\begin{proof}
Recall $\mathcal{L} = \xi \mathcal{L}_0$ for $\xi \in \textup{GL}_n(K)_{\mathbf{h}}$. Then we have
\begin{align*}
    \theta_{\mathbb{A}}^{(n)}((q,s); (\mathcal{L}, \tau)) &= \beta(q)\sum_{u \in \textup{GL}_n(K)} \lambda(\xi^{-1}uq) \Phi(u, q, s, \tau)\\
    &= \beta(q)\sum_{u \in \textup{GL}_n(K) \cap \xi Wq^{-1}} \Phi(u, q, s, \tau)\\
    &= \beta(q)\sum_{\mathfrak{a} \subseteq \mathcal{O}_K} \sum_{\substack{u \in \textup{GL}_n(K) \cap \xi Wq^{-1}\\\det ((\xi^{-1}uq)_{\mathbf{h}})=\mathfrak{a}}} \Phi(u, q, s, \tau).
 \end{align*}
Let now
\begin{equation*}
    B_{\mathfrak{a}}:= \{r \in B \mid \det(r)\widehat{\mathcal{O}}_K = \mathfrak{a}\widehat{\mathcal{O}}_K\}.
\end{equation*}
We then observe that in the above inner sum $(\xi^{-1}uq)_{\mathbf{h}} \in B_{\mathfrak{a}}$.
This is true because $\xi^{-1}uq \in W$ and $u \in \textup{GL}_n(K)$, so $u_{\mathbf{h}} \in  \textup{GL}_n(K)_{\mathbf{h}}$ and $\xi^{-1},q_{\mathbf{h}} \in \textup{GL}_n(K)_{\mathbf{h}}$. Hence, $(\xi^{-1}uq)_{\mathbf{h}} \in B$ and the determinant condition is true. By considering elementary divisors, we obtain
\begin{equation*}
    B_{\mathfrak{a}} = \bigsqcup_{\substack{D \in E \backslash B\\ \substack{\det D = \mathfrak{a}}}} ED.
\end{equation*}
Hence, the above sum can be written as
\begin{align*}
    \theta_{\mathbb{A}}^{(n)}((q,s); (\mathcal{L}, \tau)) &= \beta(q)\sum_{\mathfrak{a} \subseteq \mathcal{O}_K} \sum_{\substack{u \in \textup{GL}_n(K)\\ (\xi^{-1}uq)_{\mathbf{h}} \in B_{\mathfrak{a}}}}\Phi(u, q, s, \tau)\\
    &= \beta(q)\sum_{\mathfrak{a} \subseteq \mathcal{O}_K} \sum_{\substack{D \in E \backslash B\\ \substack{\det D = \mathfrak{a}}}} \sum_{\substack{u \in \textup{GL}_n(K)\\ (\xi^{-1}uq)_{\mathbf{h}} \in ED}}\Phi(u, q, s, \tau)\\
    &= \beta(q)\sum_{\mathfrak{a} \subseteq \mathcal{O}_K} \sum_{\substack{D \in E \backslash B\\ \substack{\det D = \mathfrak{a}}}} \sum_{\substack{u \in \textup{GL}_n(K)\\ \xi^{-1}u(qD^{-1}) \in E}}\Phi(u, q, s, \tau).
\end{align*}
But, for each $\mathfrak{a}$, we have from Lemma \ref{action of t minus},
\begin{align*}
    (\epsilon_{\mathbb{A}} \mid t^{-}(\mathfrak{a})) ((q, s); (\mathcal{L}, \tau)) &=  \sum_{\substack{D \in E \backslash B\\ \substack{\det D = \mathfrak{a}}}}\vphi_{\mf{c}}(\det D)^{n}\beta(qD^{-1})\sum_{u \in \textup{GL}_n(K)}\lambda'(\xi^{-1}uqD^{-1})\Phi(u, qD^{-1}, s, \tau) \\
    &=  \sum_{\substack{D \in E \backslash B\\ \substack{\det D = \mathfrak{a}}}}\vphi_{\mf{c}}(\det D)^{n}\beta(qD^{-1})\sum_{\substack{u \in \textup{GL}_n(K)\\ \xi^{-1}u(qD^{-1}) \in E}}\Phi(u, qD^{-1}, s, \tau) \\
    &= \psi(\mf{a})^{-n}N_{K/\mathbb{Q}}(\mathfrak{a})^{n/2}\beta(q)\sum_{\substack{D \in E \backslash B\\ \substack{\det D = \mathfrak{a}}}}\sum_{\substack{u \in \textup{GL}_n(K)\\ \xi^{-1}u(qD^{-1}) \in E}}\Phi(u, q, s, \tau),
\end{align*}
since $D$ is a finite adele and $q$ only appears in the first exponential, which only ``sees'' the infinite part. Also $\det (D) = \mathfrak{a}$, which gives
\begin{equation*}
    \vphi_{\mf{c}}(\det D)^{n}\beta(qD^{-1})= \beta(q) \psi(\mathfrak{a})^{-n}N_{K/\mathbb{Q}}(\mathfrak{a})^{n/2},
\end{equation*}
since $\psi = \vphi_{\mbf{h}}\vphi_{\mf{c}}^{-1}$. By combining these two relations, we obtain
\begin{equation*}
    \theta_{\mathbb{A}}^{(n)}((q,s);(\mathcal{L}, \tau)) = \sum_{0 \neq \mathfrak{a} \subseteq \mathcal{O}_K}   \psi(\mf{a})^{n}N_{K/\mathbb{Q}}(\mathfrak{a})^{-n/2}(\epsilon_{\mathbb{A}} \mid t^{-}(\mathfrak{a})) ((q, s);(\mathcal{L}, \tau)).\qedhere
\end{equation*}
\end{proof}
\section{The Hecke action on epsilon series}\label{The Hecke action on epsilon series}
In this Section, we will obtain an explicit formula for the action of Hecke operators from the parabolic Hecke algebra $\mathcal{H}_0$ on the epsilon series $\epsilon_{\mathbb{A}}$. Recall $\epsilon_{\mathbb{A}}$ is defined on $P_{\mb A}$. Using the rule \eqref{extension to GU}, we may extend $\epsilon_{\mathbb{A}}$ to $\mathbf{P}$.\\

Let $g = (\mu, q,s ):=  \m{\mu q & s\widehat{q} \\ 0 & \widehat{q}} \in \mathbf{P}$, with $\mu \in \mathbb{Q}_{\mathbb{A}}^{\times}$. We write uniquely $\mu = \alpha r_{\infty} \omega$ for some $\alpha \in \mathbb{Q}^{\times}$, $r_{\infty} \in \mathbb{R}_{>0}$ and $\omega \in \wh{\mathbb{Z}}^{\times}$. We then have
\begin{multline*}
g = \m{\mu q & s\widehat{q} \\ 0 & \widehat{q}} = \m{\alpha 1_n \\ & 1_n} \m{r_{\infty}^{1/2}\\& r_{\infty}^{1/2}} \m{r_{\infty}^{1/2}q & \alpha^{-1} sr_{\infty}^{-1/2}\widehat{q} \\ 0 & r_{\infty}^{-1/2}\widehat{q}}\m{\omega1_n \\ & 1_n}\\=\m{\alpha 1_n \\ & 1_n} z \m{r_{\infty}^{1/2}q & \alpha^{-1} sr_{\infty}^{-1/2}\widehat{q} \\ 0 & r_{\infty}^{-1/2}\widehat{q}}\m{\omega1_n \\ & 1_n},
\end{multline*}
where $z:= \textup{diag}(r_{\infty}^{1/2}, r_{\infty}^{1/2})$, using also that $r_{\infty}\mu^{-1}u = \alpha^{-1}$. Note that $\textup{diag}(\omega1_n, 1_n) \in \mathbf{D}_0$. \\

Using the rule \eqref{extension to GU}, and since $k=n$, we have
\begin{equation}\label{extension of epsilon to GU}
    \epsilon_{\mathbb{A}}((\mu,q,s);\mc{L}) = r_{\infty}^{-n^2/2} \phi_{\mf{c}}^{n^2}(\omega)\epsilon_{\mathbb{A}}((r_{\infty}^{1/2}q, \alpha^{-1}s); \mc{L}), 
\end{equation}
where $r_{\infty}$ is defined as the element in $\mathbb{R}^{\times}$ such that $\mu = \alpha r_{\infty}  \omega$ with $\omega \in \prod_{p} \mathbb{Z}_p^{\times}$, and $\alpha \in \mathbb{Q}^{\times}$.\\

Our next step is to study the elements of the parabolic Hecke algebras and, in particular, to obtain an analogue of \cite[Lemma 3.3.25]{andrianov_book}. In particular, we want to find a left coset decomposition of elements $\mbf{D}_0 M \mbf{D}_0$ of $\mathcal{H}_0$. Let 
\begin{equation}\label{symmetric hermitian matrices}
    \mathcal{S} := S \cap M_n(\mathcal{O}_K),
\end{equation}
where $S$ as in \eqref{hermitian matrices in K}. Analogously to Andrianov, we have:
\begin{lemma} \label{triangular decomposition} Let $p$ be a finite prime of $\mathbb{Q}$ with $p$ prime to $\mathfrak{bc}$ (inert or split). We then have:
\begin{enumerate}
\item The left coset $\mbf{D}_{0,p} M$ for a matrix $M = \m{A & B \\ 0 & D} \in \mathbf{P}$ is uniquely determined by the valuation of $\mu(M)$, the left coset $E_p D$ with $E_p = \gl_n(\mathcal{O}_{K_p})$ and, for a fixed $D$, the right residue class $B + \mathcal{S}_p D$ with $\mathcal{S}_p$ the integral ``symmetric'' matrices locally. Moreover, each choice above uniquely determines a left coset. 

\item For $M = \m{A & B \\ 0 & D} \in \mathbf{P}$ we have
\[
\mathbf{D}_{0,p} M \mathbf{D}_{0,p}= \bigsqcup_{L,N} \mathbf{D}_{0,p} \m{\mu \widehat{L} & N \\ 0 & L},
\]
where $L \in E_p \backslash E_p D E_p$ and $N \in B(L,M) \pmod{ \mathcal{S}L}$ with 
$$B_p(L,M) := \left\{ N; \m{\mu \widehat{L} & N \\ 0 & L}\in \mathbf{D}_{0,p}M\mathbf{D}_{0,p}  \right\}.$$
\end{enumerate}
\end{lemma}

That is, Andrianov's statement also holds in the unitary case when one works locally over $\mathbb{Q}_p$. In this case, $K_p = \mathbb{Q}_p \times \mathbb{Q}_p$ for $p$ split and $K_p$ a quadratic extension of $\mathbb{Q}_p$ when $p$ is inert. In particular in the statement above $E_p = \gl_n(\mathbb{Z}_p) \times \gl_n(\mathbb{Z}_p)$ and $\mathcal{S}_p = \{ (a,a^{t}) \mid a\in M_n(\mathbb{Z}_p)\}$ for $p$ split.\newline

Next, in analogy to \cite{andrianov_book}, we define certain sums, the ``Gauss sums''. We first define
\begin{equation}\label{T set}
    \widehat{\mc{S}} := \{x \in \mathcal{S} \mid \textup{tr}(yx) \in \mathbb{Z} \textup{ for all } y \in \mathcal{S} \cap M_n(\mathcal{O}_K)\},
\end{equation}
the dual lattice of $\mc{S}$ (see \eqref{symmetric hermitian matrices}). This is a $\mathbb{Z}$-lattice in $\mc{S}$, so $\widehat{\mc{S}}_p$ is meaningful for all rational primes $p$ (see \cite[p. 147]{Shimura_Euler_Product}).

\begin{defn}[Gauss sums]\label{gauss sum defn}
    Let $M = \m{A&B\\0& D} \in \mbf{P}_{\mbf h}$. For $Q \in \widehat{\mathcal{S}}_{\mbf{h}}$ and $L \in EDE$, we define the Gauss sum
    \begin{equation*}
        s(Q, L, M) := \sum_{N \in B(L, M)\pmod{L}}e_{\mathbf{h}}\{QNL^{-1}\}.
    \end{equation*}
\end{defn}
We now have the following Lemma.
\begin{lemma}\label{gauss sum lemma 1}
    Let $A, B \in E$. We then have 
    \begin{equation*}
        s(\tau[A], L, M) = s(\tau, \widehat{A}LB, M),
    \end{equation*}
    where $s(-)$ denote the Gauss sums of Definition \ref{gauss sum defn}. In particular, 
    \begin{equation*}
        s(\tau, L,M) = s(\tau, LB, M)
    \end{equation*}
    for all $B \in E$. 
\end{lemma}
\begin{proof}
    We have 
    \begin{align*}
        s(\tau[A], L, M) &= \sum_{N \in  B(L, M) \pmod{L}} e_{\mathbf{h}}\{\tau[A]NL^{-1}\} = \sum_{N \in  B(L, M) \pmod{L}} e_{\mathbf{h}}\{\tau (ANB)B^{-1}L^{-1}A^{*}\}\\ &= \sum_{N \in  B(L, M) \pmod{L}} e_{\mathbf{h}}\{\tau (ANB)(\widehat{A}L B)^{-1}\}.
    \end{align*}
    We now have the relation
    \begin{equation*}
        \widehat{V} \{B(L, M)\pmod{L}\} V' = B(VLV', M) \pmod{VLV'}
    \end{equation*}
    for $V,V' \in E$ (cf. with \cite[Eq. (5.1.19)]{andrianov_book}). Indeed, to see this, we recall that 
    \[
    B(L,M) = \left\{ N; \m{\mu \widehat{L} & N \\ 0 & L} \in \mbf{D}_0 M \mbf{D}_0\right\},
    \]
    and hence 
\[
    B(VLV',M) = \left\{ \widetilde{N}; \m{\mu \widehat{VLV'} & \widetilde{N} \\ 0 & VLV'} \in \mbf{D}_0 M \mbf{D}_0\right\},
    \]
    We now observe that $\m{\widehat{V} & 0 \\ 0 & V }, \m{\widehat{V'} & 0 \\ 0 & V' } \in \mbf{D}_0$ for $V,V' \in E$, and hence in particular
    \[
    \m{\widehat{V} & 0 \\ 0 & V }\m{\mu \widehat{L} & N \\ 0 & L}\m{\widehat{V'} & 0 \\ 0 & V' } = \m{\mu \widehat{VLV'} & \widehat{V}NV' \\ 0 & VLV'} \in \mbf{D}_0M\mbf{D}_0.
    \]
In particular, from the above we see that for a fixed set of representatives $\{N_i\}$ of $B(L,M) \pmod{L}$, taking $\widetilde{N_i} := \widehat{V}N_iV'$ gives a set of representatives of  $B(VLV',M) \pmod{VLV'}$ and hence establishes the claim. Using this with $V' = B$ and $V = \widehat{A}$, we obtain after setting $N' = ANB$:
    \begin{equation*}
        s(\tau[A], L, M) = \sum_{N' \in B(\widehat{A}LB, M) \pmod{\widehat{A}LB}} e_{\mathbf{h}}\{\tau N'(\widehat{A}LB)^{-1}\} = s(\tau,  \widehat{A}LB, M). \qedhere
    \end{equation*}
\end{proof}
% Following Andrianov we now introduce the notation $\epsilon(L_1,L_2,T)$ where $L_1,L_2$ are lattices in $(V,\tau)$ and $T \in \mathbf{H}_0$, a Hecke operator of the parabolic Hecke algebra. We start with the following observation. If we write $U(\tau)$ for the unitary group corresponding to the hermitian matrix $\tau$, then we have an equality $U(\tau) = U(\alpha \tau)$ for all $\alpha \in \mathbb{Q}^{\times}_{+}$, since for any matrix $M \in \textup{GL}_n(K)$ we have $M \tau M^* = \tau$ if and only if $M \alpha \tau M^* = \alpha \tau$. In particular the genus, and the class of a lattice $L \subset V$ depends only on $\tau$ and not on scaling by $\alpha$. \\
% Note 
% \begin{equation*}
%     E\widehat{D}Eq^{-1} \cap \textup{GL}_n(K) = \bigsqcup_{L \in E\backslash EDE} (E\widehat{L}q^{-1} \cap \textup{GL}_n(K)).
% \end{equation*}

%That is we obtain,

%\[
%(\epsilon_{\mathbb{A}} \mid T_{M})(q, s) = \sum_L G(\tau,L) \sum_{{\substack{u \in \textup{GL}_n(K)\cap Eq^{-1} }}}\Phi_{\infty}(u_{\infty}, q_{\infty}, s_{\infty}) e_{\mathbf{h}}(\textup{tr}(\tau[u_{\mathbf{h}}]s )) 
%\]
%with $G(\tau,L) := \sum_N e_{\mathbf{h}}(\textup{tr}(-\tau(\mu^{-1} L^* N)$, which is non-zero if and only if $(\mu^{-1} \tau[L^*])$ in the dual of $\mathcal{S}_v$. That is
%\[
%(\epsilon_{\mathbb{A}} \mid T_{M})(q, s) = \sum_L G(\tau,L) 
%\]
We also have the following Lemma.
\begin{lemma}\label{invariance of epsilon under U}
    Let $\gamma \in \textup{GL}_n(K)$. We then have for any $\mathcal{O}_K$-lattice $\mathcal{L}$, that
\begin{equation}
    \epsilon_{\mathbb{A}}(\gamma \mathcal{L}, \tau) = \epsilon_{\mathbb{A}}(\mathcal{L}, \tau[\gamma]).
\end{equation}
\end{lemma}
\begin{proof}
    Assume $\mathcal{L} = \xi\mathcal{L}_0$ for some $\xi \in \textup{GL}_n(K)_{\mathbf{h}}$. Then $\gamma \mathcal{L} = \gamma \xi \mathcal{L}_0$. Hence
    \begin{equation*}
        \epsilon_{\mathbb{A}}(\gamma \mathcal{L}, \tau) = \beta(q)\sum_{u \in \textup{GL}_n(K)} \ell'(\xi^{-1}\gamma^{-1}uq) \Phi(u,q,s,\tau).
    \end{equation*}
    Let now $v:= \gamma^{-1}u$. Then $v$ runs through $\textup{GL}_n(K)$ as $u \in \textup{GL}_n(K)$. Additionally,
    \begin{equation*}
        \Phi(u,q,s,\tau) = \Phi(\gamma v, q,s, \tau) = \Phi(v, q, s, \tau[\gamma]).
    \end{equation*}
    Hence, 
    \begin{equation*}
        \epsilon_{\mathbb{A}}(\gamma \mathcal{L}, \tau) = \beta(q)\sum_{v \in \textup{GL}_n(K)} \ell'(\xi^{-1}v q) \Phi(v,q,s,\tau[\gamma]) = \epsilon_{\mathbb{A}}(\mathcal{L}, \tau[\gamma]).\qedhere
    \end{equation*}
\end{proof}
We can now compute the action of $\mathcal{H}_0$ on $\epsilon_{\mathbb{A}}$ as follows. First, we will say that two lattices $(\mc{L}, \tau)$ and $(\mc{L'},\rho)$ are isometric if there exists some $\gamma \in \textup{GL}_n(K)$ such that 
\begin{equation}\label{isometric lattices}
    \mc{L}'= \gamma \mc{L}, \qquad \gamma^{*}\rho \gamma = \tau.
\end{equation}
We will denote by $[(\mc{L}, \tau)]$ the resulting isometry class. We can now extend Andrianov's \cite[Proposition 5.1.3]{andrianov_book} to our Hermitian setting.
\begin{proposition}\label{action of hecke on epsilon}
    Let 
    \begin{equation*}
        T_M := \mbf{D}_0 \m{\mu\widehat{D} & B\\ 0 &D}\mbf{D}_0 = \mbf{D}_0M\mbf{D}_0\in \mathcal{H}_0,
    \end{equation*}
    where $\mu \in \mathbb{Q}_{\mathbb{A}}^{\times}$ is written uniquely as $\mu  = \alpha r_{\infty}\omega$, with $\alpha \in \mathbb{Q}^{\times}$, $r_{\infty} \in \mb{R}_{>0}$ and $\omega \in \wh{\mb{Z}}^{\times}$. Then, for a $\mathbb{Z}$-integral $\mc{O}_K$-lattice $\mathcal{L} := \xi \mathcal{L}_0$, with $\xi \in \textup{GL}_n(K)_{\mathbf{h}}$, we have 
    \begin{equation}\label{epsilon under hecke}
         (\epsilon_{\mathbb{A}} \mid T_M )((q,s); (\mathcal{L}, \tau)) = \sum_{[(\mc{L'}, \rho)]} \epsilon((\mathcal{L},\tau), (\mathcal{L}',\rho), T_M) \epsilon_{\mathbb{A}}((q,s); (\mathcal{L}', \rho)),
    \end{equation}
    where the sum is taken over isometry classes, and
    \begin{equation}\label{coefficients}
        \epsilon((\mathcal{L},\tau), (\mathcal{L}',\rho), T_M) := \vphi_{\mf{c}}(\omega)^{-n^2}\beta(D^{*}) \sum_{\substack{R \in EDE/E\\ [(\xi\widehat{R}\mathcal{L}_0,\alpha \tau)] = [(\mc{L'}, \rho)]}} s((-\omega^{-1}\tau)[\xi], R, M),
    \end{equation}
    where $s(-)$ denote the Gauss sums of Definition \ref{gauss sum defn}.
\end{proposition}
\begin{proof}
Choose a left-coset decomposition
\[
 T_M=
 \bigsqcup_{\substack{L\in E\backslash EDE\\
                     N\in B(L,M) \pmod{L}}}
 \mathbf D_{0}
 \begin{pmatrix}\mu\widehat L&N\\0&L\end{pmatrix}.
\]
All representatives are taken to be the identity at infinity and at the primes dividing $\mathfrak c$.
Note that
\begin{equation*}
    \m{q & s\widehat{q} \\ & \widehat{q}} \m{\mu \widehat{L} & N\\ & L}^{-1} = \m{\mu^{-1}qL^{*} & (-\mu^{-1}qL^{*}Nq^{*} + s)\widehat{qL^*}\\ & \widehat{qL^{*}}}.
\end{equation*}
Substituting $\mu=r_{\infty}\alpha \omega$, we obtain that under the above multiplication,  $q \longmapsto r_{\infty}^{-1/2} qL^{*}$ and $s \longmapsto \alpha(-\mu^{-1}qL^{*}Nq^{*}+s)$. \\

We observe that for $L \in EDE$, we have $|\det(L)|_{K} = |\det(D)|_{K}$, since $|a|_{K}=1$ for all $a \in E$. Additionally, since $\varphi$ is unramified for $v \nmid \mathfrak{c}$ (i.e. $\varphi_v(x)=1 \, \forall x\in \mathcal{O}_{K,v}^{\times}$), and both $L,D$ are trivial at the places dividing $\mathfrak{c}$, we also obtain $\varphi_{\mathbf{h}}(\det(L)) = \varphi_{\mathbf{h}}(\det(D))$. Hence, $\beta(L) = \beta(D)$ for all $L \in EDE$. Finally, $\beta(r_{\infty}^{-1/2}) = r_{\infty}^{-n^2/2}$. Hence $\beta(r_{\infty}^{-1/2}qL^{*}) = r_{\infty}^{-n^2/2}\beta(qD^{*})$. Using \eqref{extension of epsilon to GU}, we obtain:
\begin{align*}
&(\epsilon_{\mathbb{A}} \mid T_M)((q,s); \mathcal{L}, \tau)\\
&= r_{\infty}^{n^2/2}\vphi_{\mf{c}}^{-n^2}(\omega)\sum_{L, N} \sum_{u \in \textup{GL}_n(K)}\beta(r_{\infty}^{-1/2}qL^{*})\lambda'(\xi^{-1}uqL^{*})\Phi(u, r_{\infty}^{-1/2} qL^{*}, \alpha(-\mu^{-1}qL^{*}Nq^{*}+s), \tau)\\
&= \!\begin{multlined}[t][10.5cm]\vphi_{\mf{c}}^{-n^2}(\omega)\beta(qD^{*})\sum_{L,N}\sum_{u \in \textup{GL}_n(K)}\lambda'(\xi^{-1}uqL^{*})\Phi_{\infty}(u_{\infty},r_{\infty}^{-1/2}q_{\infty}, \alpha s_{\infty}, \tau)\\ \times \Phi_{\mathbf{h}}(u, qL^{*}, \alpha(-\mu^{-1}qL^{*}Nq^{*}+s), \tau)\end{multlined}\\
&= \vphi_{\mf{c}}^{-n^2}(\omega)\beta(qD^{*})\sum_{L,N}\sum_{u \in \textup{GL}_n(K)} \lambda'(\xi^{-1}uqL^{*})\Phi(u,q,s, \alpha\tau)e_{\mathbf{h}}\{\alpha\tau[u_{\mathbf{h}}](-\mu^{-1}q L^* Nq^* )\}\\
&= \vphi_{\mf{c}}^{-n^2}(\omega) \beta(qD^{*})\sum_{L}\sum_{u \in \textup{GL}_n(K)}\lambda'(\xi^{-1}uqL^{*})\Phi(u,q,s, \alpha\tau)s((-\omega^{-1}\tau)[uqL^{*}], L, M)\\
&=  \vphi_{\mf{c}}^{-n^2}(\omega)\beta(qD^{*})\sum_{L}\sum_{u \in \textup{GL}_n(K)}\lambda'(\xi^{-1}uqL^{*})\Phi(u,q,s, \alpha\tau)s((-\omega^{-1}\tau)[\xi][\xi^{-1}uqL^{*}], L, M).
\end{align*}
Notice now that $\lambda'(\xi^{-1}uqL^{*}) \neq 0 \iff \xi^{-1}uqL^{*}\in E$, and in that case
\begin{equation*}
    s((-\omega^{-1}\tau)[\xi][\xi^{-1}uqL^{*}], L, M) = s((-\omega^{-1}\tau)[\xi], \xi^{*}\widehat{u}\widehat{q}, M)
\end{equation*}
using Lemma \ref{gauss sum lemma 1}. Therefore, the above can be written as
\begin{align*}
    &(\epsilon_{\mathbb{A}} \mid T_M)((q,s); \mathcal{L}, \tau)\\
    &=\vphi_{\mf{c}}^{-n^2}(\omega)\beta(qD^{*})\sum_{u \in \textup{GL}_n(K)}s((-\omega^{-1}\tau)[\xi], \xi^{*}\widehat{u}\widehat{q}, M)\Phi(u,q,s,\alpha\tau)\sum_{L \in E\backslash EDE}\lambda'(\xi^{-1}uqL^{*})\\
    &=\vphi_{\mf{c}}^{-n^2}(\omega)\beta(qD^{*})\sum_{u \in \textup{GL}_n(K)}s((-\omega^{-1}\tau)[\xi], \xi^{*}\widehat{u}\widehat{q}, M)\Phi(u,q,s,\alpha\tau)\sum_{i=1}^{t}\lambda'(\xi^{-1}uqL_{i}^{*}),
\end{align*}
where we have fixed a coset decomposition
\begin{equation}\label{coset decomposition}
    EDE = \bigsqcup_{i=1}^{t} E L_i.
\end{equation}
We now observe that for $u \in \textup{GL}_n(K)$,
\begin{equation*}
\sum_{i=1}^{t}\lambda'(\xi^{-1}uqL_{i}^{*})=
    \begin{cases}
        1 \textup{ if }  u \in \xi E\widehat{D}Eq^{-1}\\
        0 \textup{ otherwise}
    \end{cases}.
\end{equation*}
Indeed, if there is some $j$ such that $\xi^{-1}uqL_j^{*} \in E$, then $u \in \xi E\widehat{L}_jq^{-1} \subset \xi E\widehat{D}Eq^{-1}$, since $L_i \in EDE$ for all $i$. Additionally, if $u \in \xi E\widehat{D}Eq^{-1}$, then $u \in \xi E\widehat{L}_i q^{-1}$ for a unique $i$, since from \ref{coset decomposition}, we also have the decomposition
\begin{equation*}
    \xi E\widehat{D}Eq^{-1} = \bigsqcup_{i=1}^{t} \xi E \widehat{L}_iq^{-1}.
\end{equation*}
Therefore, we have
\begin{equation*}
    (\epsilon_{\mathbb{A}} \mid T_M)((q,s); \mathcal{L}, \tau)=\vphi_{\mf{c}}^{-n^2}(\omega)\beta(qD^{*})\sum_{u \in \textup{GL}_n(K) \cap \xi E\widehat{D}Eq^{-1}}s((-\omega^{-1}\tau)[\xi], \xi^{*}\widehat{u}\widehat{q}, M)\Phi(u,q,s,\alpha\tau).
\end{equation*}
Consider now a decomposition
\begin{equation*}
    EDE = \bigsqcup_{i=1}^{t}R_iE \implies \xi E\widehat{D}Eq^{-1} = \bigsqcup_{i=1}^{t}\xi \widehat{R}_iEq^{-1}.
\end{equation*}
Then, the above sum can be written as 
\begin{align*}
    &(\epsilon_{\mathbb{A}} \mid T_M)((q,s); \mathcal{L}, \tau)\\
    &= \vphi_{\mf{c}}^{-n^2}(\omega)\beta(qD^{*})
    \sum_{i=1}^{t}\sum_{u \in \textup{GL}_n(K)\cap \xi \widehat{R}_iEq^{-1}} s((-\omega^{-1}\tau)[\xi], \xi^{*}\widehat{u}\widehat{q}, M)\Phi(u,q,s,\alpha\tau)\\
    &=\vphi_{\mf{c}}^{-n^2}(\omega)\beta(qD^{*}) \sum_{u \in \textup{GL}_n(K)}s((-\omega^{-1}\tau)[\xi], \xi^{*}\widehat{u}\widehat{q}, M)\Phi(u,q,s,\alpha\tau) \sum_{i=1}^{t} \ell'(R_i^{*}\xi^{-1}uq)\\
    &= \vphi_{\mf{c}}^{-n^2}(\omega)\beta(qD^{*})
    \sum_{R \in EDE/E} \sum_{u \in \textup{GL}_n(K)}\lambda'(R^{*}\xi^{-1}uq)s((-\omega^{-1}\tau)[\xi], \xi^{*}\widehat{u}\widehat{q}, M)\Phi(u,q,s,\alpha\tau)
\end{align*}
using the same reasoning as above. But now, if $R^{*}\xi^{-1}uq\in E$, then
\begin{equation*}
    s((-\omega^{-1}\tau)[\xi], \xi^{*}\widehat{u}\widehat{q}, M) = s((-\omega^{-1}\tau)[\xi], (\xi^{*}\widehat{u}\widehat{q})(R^{*}\xi^{-1}uq)^{*}, M) = s((-\omega^{-1}\tau)[\xi], R, M),
\end{equation*}
using Lemma \ref{gauss sum lemma 1}. Therefore, we obtain
\begin{align*}
    &(\epsilon_{\mathbb{A}} \mid T_M)((q,s); \mathcal{L}, \tau)\\
    &= \vphi_{\mf{c}}^{-n^2}(\omega)\beta(D^{*})\sum_{R \in EDE/E} s((-\omega^{-1}\tau)[\xi], R, M)\times \beta(q)\sum_{u \in \textup{GL}_n(K)}\lambda'(R^{*}\xi^{-1}uq)\Phi(u,q,s,\alpha\tau)\\
    &= \vphi_{\mf{c}}^{-n^2}(\omega)\beta(D^{*})\sum_{R \in EDE/E}s((-\omega^{-1}\tau)[\xi], R, M)\epsilon_{\mathbb{A}}((q,s); (\xi\widehat{R}\mathcal{L}_0, \alpha\tau)).
\end{align*}
Finally, grouping the terms according to their isometry classes and using \ref{invariance of epsilon under U} gives \eqref{epsilon under hecke}.
\end{proof}
\begin{lemma}\label{multiplicative epsilon}
    The function $\epsilon((\mathcal{L},\tau), (\mathcal{L}',\rho), T_M)$ of \eqref{coefficients} has the following properties.
    \begin{enumerate}
        \item It depends only on the isometry classes $[(\mathcal{L},\tau)], [(\mc{L}', \rho)]$.
        \item If the class $[(\mc{L}, \tau)]$ (resp. $[(\mc{L}', \rho)]$) and $T_M$ are fixed, then $\epsilon((\mathcal{L},\tau), (\mathcal{L}', \rho), T_M)=0$ for all but a finite number of the classes $[(\mc{L}', \rho)]$ (resp. $[(\mc{L}, \tau)]$).
       \item For all lattices $(\mathcal{L},\tau), (\mathcal{L}',\rho)$ and $T,T' \in \mc{H}_0$, we have 
        \begin{equation*}
            \epsilon((\mathcal{L},\tau), (\mathcal{L}',\rho), TT') = \sum_{\{[(\mc{M},\pi)]\}} \epsilon((\mc{L},\tau), (\mc{M},\pi), T)\epsilon((\mc{M},\pi), (\mc{L}',\rho), T').
       \end{equation*} 
    \end{enumerate}
    \end{lemma}
    \begin{proof}
       \begin{enumerate}
            \item This is a straightforward calculation, using the notion of the isometry class as above.
            \item For fixed $\mc{L}=\xi\mc{L}_0$ and fixed $T_M$, $(\mc{L}', \rho)$ runs through the classes $(\xi\widehat{R}\mathcal{L}_0,\alpha\tau)$ for $R \in EDE/E$, a finite number.
            \item 
        The adelic epsilon series attached to distinct isometry classes are linearly independent. We prove this slightly more generally in Lemma \ref{linear independence} and we refer the reader there to understand the idea. Here, it follows similarly, with $\lambda'$ in place of $\lambda$ there. Indeed, writing $\mc{L}'=\xi'\mc{L}_0$ and $p_{\xi'}=\textup{diag}(\xi',\widehat{\xi'})$, the $\rho$-Fourier coefficient of the $p_{\xi'}$-component of $\epsilon_{\mathbb{A}}(\cdot;(\mc{M},\pi))$ is
        \begin{equation*}
            \beta(\xi')\#\left\{ u\in\textup{GL}_n(K) \mid \pi[u]=\rho,\quad u\mc{L}'=\mc{M}\right\}.
        \end{equation*}
        Here, the equality $u\mc{L}'=\mc{M}$ follows from the condition involving $\lambda'$. Thus this coefficient vanishes unless $(\mc{M},\pi)$ and $(\mc{L}',\rho)$ are isometric, in which case it equals
        $\beta(\xi')|\textup{Aut}(\mc{L}',\rho)|\neq0$. Hence, these epsilon series are linearly independent. Applying Proposition \ref{action of hecke on epsilon} to both sides of
        \begin{equation*}
            \epsilon_{\mathbb{A}}(\cdot;(\mc{L},\tau))\mid TT' = \left(\epsilon_{\mathbb{A}}(\cdot;(\mc{L},\tau))\mid T\right)\mid T',
        \end{equation*}
        and comparing the coefficients of these linearly independent epsilon series gives the required identity. All sums are finite by (2) above.\qedhere
        \end{enumerate}
    \end{proof}

Let us now consider when $s(\omega^{-1}\tau[p],R,M) = 0$. We have the following Lemma.
\begin{lemma}\label{non-vanishing lemma}
With the same notation as in Definition \ref{gauss sum defn}, the Gauss sum $s(Q,L,M)$ vanishes if $$\mu Q[\widehat{L}] \not \in \prod_{p \nmid \infty }\widehat{\mc{S}}_{p}.$$
% With notation as above (in particular $R \in EDE$),  $s(\omega^{-1}\tau[p],R,M) = 0$ when $\alpha\tau[p\widehat{R}] \not \in \prod_v T_v$ where $T_v$ is the dual with respect to the trace map of $\mathcal{S}(\mathbb{Z})$ (see \cite[page 147]{Shimura_Euler_Product}).

%$p\widehat{r}\mathcal{O}_K^n $ is not the dual lattice to $\mathcal{O}_K^n$ (wrt to $\alpha\tau$) for all finite places.
    
\end{lemma}
\begin{proof}
We recall that 
\begin{equation*}
s(Q,L,M) = \sum_{N \in  B(L, M) \pmod{L}} e_{\mathbf{h}}\{Q NL^{-1}\}.
\end{equation*}
If $\{N_i\}$ is a finite set of representatives of $B(L,M) \pmod{L}$, so is also $\{N_i+ \mu\widehat{L}A\}$ for any $A \in \mathcal{S}_{\mathbf{h}}$. We then have
\begin{equation*}
s(Q,L,M) = e_{\mathbf{h}}\{\mu Q \widehat{L}AL^{-1}\}\sum_i e_{\mathbf{h}}\{Q N_iL^{-1}\}=e_{\mathbf{h}}\{\mu Q[\widehat{L}]A\}s(Q,L,M).
\end{equation*}
In particular, if $e_{\mathbf{h}}\{\mu Q[\widehat{L}]A\} \neq 1$, then $s(Q,L,M) = 0$. Now, if $\mu Q[\widehat{L}] \not \in \prod_{p}\widehat{\mc{S}}_p$, there is some $p$ and $y_p \in \mathcal{S}_p \cap M_n(\mathcal{O}_{K,p})$ such that $\textup{tr}(\mu_p Q_p[\widehat{L}_p]y_p) \not\in \mathbb{Z}_p$. By then taking $A \in \mathcal{S}_{\mathbf{h}}$ given by $A_p = y_p$ and $A_{q} = 0$ for all $q \neq p$, we obtain 
\begin{equation*}
    e_{\mathbf{h}}\{\mu Q[\widehat{L}]A\} = e_{p}\{\mu_p Q_p[\widehat{L}_p]y_p\} \neq 1.\qedhere
\end{equation*}
%That is, if the lattice $p\hat{r}O_K^n$ is not the dual of $O_K^n$, with the respect to $\alpha\tau$ for all finite places, then $s(\omega^{-1}\tau[p],r,M)=0$.
\end{proof}
Let now
\begin{equation*}
    \widetilde{\mbf{P}} = \left\{M= \m{\mu \widehat{D} & B\\0 & D} \in \mbf{P}_{\mathbb{A}}  \mid N(\det D) \alpha^{-n}=1,\, \,\, M_p \in \mathbf{D}_{0,p} \,\  \forall p \mid \mf{bc}\right\}.
\end{equation*}
and its subgroup
\begin{equation*}
    \widetilde{\mbf{P}}^+ = \left\{\m{\mu \widehat{D} & B\\0 & D} \in \widetilde{\mbf{P}} \mid \alpha \in N_{K/\mathbb{Q}}(K^{\times})\right\},
\end{equation*}
where we write $\mu = \alpha r_{\infty}\omega$ as usual. We define the Hecke algebras 
\begin{equation}\label{even operators}
    \widetilde{\mc{H}}_0 := \mathcal{H}(\mbf{D}_0, \widetilde{\mbf{P}}), \qquad \widetilde{\mc{H}}_0^+ := \mathcal{H}(\mbf{D}_0, \widetilde{\mbf{P}}^+)
\end{equation}
We have the following Proposition.
\begin{proposition}\label{level divisor det}
    Let $\mc{L} = \xi \mc{L}_0$, where $\xi \in \textup{GL}_n(K)_{\mathbf{h}}$ be a lattice and let $T_M \in \widetilde{\mc{H}}_0$. Suppose $\epsilon((\mc{L},\tau), (\mc{L'},\rho), T_M) \neq 0$. Then the lattice $(\mc{L}', \rho)$ is integral, and its volume, level and divisor are equal to those of $(\mc{L}, \tau)$. If in addition $T_M \in \widetilde{\mc{H}}_0^+$, i.e. $\alpha \in N_{K/\mathbb{Q}}(K^{\times})$, we have that $(\mc{L'}, \rho)$ and $(\mc{L}, \tau)$ are locally isometric for all finite primes $p$.
\end{proposition}
\begin{proof}
   Recall we write $\mu = r_{\infty}\alpha\omega$. From \eqref{coefficients}, and since $\epsilon((\mc{L},\tau), (\mc{L}',\rho), T_M) \neq 0$, there exists $R \in EDE$ such that
    \begin{equation*}
        s((-\omega^{-1}\tau)[\xi],R,M)\neq 0,
        \qquad
        \xi\widehat{R}\mc{L}_0 =\gamma\mc{L}'
    \end{equation*}
    for some $\gamma$ such that $\rho = \alpha \tau[\gamma]$.\\

    We will first show that $(\mc{L}', \rho)$ is $\mathbb{Z}$-integral. It suffices to show $\rho[x] \in \mathbb{Z}$ for all $x \in \mc{L}'$, or equivalently that $\rho[x] \in \mathbb{Z}_p$ for all $x \in \mc{L}'_p$. Since $s((-\omega^{-1}\tau)[\xi],R,M)\neq 0$, from Lemma \ref{non-vanishing lemma} we obtain
    \begin{equation*}
        \alpha \tau[\xi \widehat{R}]_p \in \widehat{\mc{S}}_p 
    \end{equation*}
    for all finite primes $p$. Let now $x \in \mc{L}'_{p}$. We have that $\gamma x = \xi_p \wh{R}_py$ for some $y \in \mc{L}_{0,p}$. Hence
    \begin{equation*}
        \rho[x] = \alpha\tau[\gamma x] = \alpha\tau[ \xi_p \wh{R}_p y] = \textup{tr}(\alpha\tau[\xi \wh{R}]yy^{*}) \in \mathbb{Z}_p,
    \end{equation*}
    showing the integrality.\\
    
    Now, since $ \det(\rho)=\alpha^n\det(\tau)N_{K/\mathbb Q}(\det\gamma)$, we have
    \begin{multline*}
        \vol(\mc{L}', \rho)= \det(\rho) N_{K/\mb{Q}}(\det(\gamma^{-1} \xi \wh{R})) = \det(\rho)N_{K/\mb{Q}}(\det\gamma)^{-1}N_{K/\mb{Q}}(\det (\xi \wh{R})) \\= \det(\tau)N_{K/\mb{Q}}(\det\xi)\alpha^n N_{K/\mb{Q}}(\det(\wh{R})) = \vol(\mc{L}, \tau),
    \end{multline*}
    since $\alpha^n N_{K/\mb{Q}}(\det\wh{R}) = \alpha^n N_{K/\mb{Q}}(\det\wh{D})=1$ as $T_M \in \widetilde{H}_{0}$ and $R \in EDE$.\\

We now consider the level and divisor. Let first $p \nmid \ell(\mc{L},\tau)$. Let $\{e_1, \ldots, e_{2n}\}$ denote a $\mathbb{Z}$-basis for $\mc{L}$. We consider then the integral, symmetric Gram matrix $Q$ with even diagonal entries associated to $q_{\tau}$ with respect to this basis. From the proof of \cite[Lemma 5.1.5]{andrianov_book}, we have
    \begin{equation}\label{lev, det, div 1}
        \ell(\mc{L},\tau) \mid \det(Q),
        \qquad
        \det(Q) \mid \ell(\mc{L},\tau)^{2n}.
    \end{equation}
    Additionally, since $\delta(\mc{L},\tau)^{-1}Q$ is an integral matrix, we also obtain
    \begin{equation}\label{lev, det, div 2}
        \delta(\mc{L},\tau)^{2n} \mid \det(Q).
    \end{equation}
    Now, it is well-known that
    \begin{equation*}
        \det(Q) = |d_K|^n\vol(\mc{L},\tau)^2.
    \end{equation*}
    But $\vol(\mc{L},\tau)=\vol(\mc{L}',\rho)$, which shows $\det(Q)=\det(Q')$, where $Q'$ is the analogous Gram matrix for $q_{\rho}$ with respect to a $\mathbb Z$-basis for $\mc{L}'$.\\

    From \eqref{lev, det, div 1}, we see that $p\nmid\det(Q)$, which means $p\nmid\det(Q')$. Hence, applying \eqref{lev, det, div 2} to $(\mc{L}',\rho)$ gives $p\nmid\delta(\mc{L}',\rho)$. Additionally, applying \eqref{lev, det, div 1} to $(\mc{L}',\rho)$ gives $p\nmid\ell(\mc{L}',\rho)$. Hence, in this case
    \begin{equation*}
        v_p(\ell(\mc{L},\tau))=v_p(\ell(\mc{L}',\rho)),
        \qquad
        v_p(\delta(\mc{L},\tau))=v_p(\delta(\mc{L}',\rho)).
    \end{equation*}
    \medskip
Assume now $p\mid\ell(\mc{L},\tau)$. By our assumption on the bad primes, $R_p\in E_p$ and $\alpha\in\mathbb Z_p^\times$. We then have
    \begin{equation*}
        \gamma\mc{L}'_p=\xi_p\widehat{R}_p\mc{L}_{0,p}
        =\xi_p\mc{L}_{0,p}=\mc{L}_p.
    \end{equation*}
    Hence,
    \begin{equation*}
        \tau[\mc{L}_p]=\tau[\gamma\mc{L}'_p]
        =\alpha^{-1}\rho[\mc{L}'_p],
    \end{equation*}
    which shows $v_p(\delta(\mc{L},\tau))=v_p(\delta(\mc{L}',\rho))$.\\
    
    Let now $B$ and $B'$ denote the associated bilinear forms
    \begin{equation*}
        B(x,y):=\textup{tr}_{K/\mathbb Q}(\tau(x,y)),
        \qquad
        B'(x,y):=\textup{tr}_{K/\mathbb Q}(\rho(x,y))=\alpha B(\gamma x,\gamma y).
    \end{equation*}
    We denote by $\widehat{\mc{L}}$ and $\widehat{\mc{L}'}$ the duals with respect to $B$ and $B'$, respectively. We then have
    \begin{equation*}
        \widehat{\mc{L}}_p
        =\{x\in V_p:B(x,\mc{L}_p)\subseteq\mathbb Z_p\}.
    \end{equation*}
    Now
    \begin{equation*}
        x\in\widehat{\mc{L}}_p
        \iff B(x,\gamma\mc{L}'_p)\subseteq\mathbb Z_p
        \iff \alpha^{-1}B'(\gamma^{-1}x,\mc{L}'_p)
                    \subseteq\mathbb Z_p
        \iff x\in\gamma\alpha\widehat{\mc{L}}_p'.
    \end{equation*}
    Hence, we obtain
    \begin{equation*}
        \widehat{\mc{L}}_p=\gamma\alpha\widehat{\mc{L}}_p',
        \qquad
        \tau[\widehat{\mc{L}}_p]
          =\alpha\rho[\widehat{\mc{L}}_p'].
    \end{equation*}
    Since $v_p(\alpha)=0$, we obtain
    $v_p(\ell(\mc{L},\tau))=v_p(\ell(\mc{L}',\rho))$.\\
    
Assume now $\alpha \in N_{K/\mathbb{Q}}(K^{\times})$. Consider first the primes $p \mid |d_K|\ell(\mc{L},\tau)$. We then have that  
$R_p \in E_p$ and hence $\mc{L}_p = \gamma \mc{L}_p'$. But also $\alpha \in \mathbb{Z}_p^{\times}$, so there exists $\beta \in \mathcal{O}_{K,p}^{\times}$ such that $\alpha = N_{K_p/\mb{Q}_p}(\beta)$. Then $\mc{L}_p = \beta \gamma \mc{L}_p'$ since $\beta$ is a unit in $\mathcal{O}_{K,p}$, and $\tau[\beta \gamma] = \alpha \tau[\gamma] = \rho $. Hence, the lattices $(\mc{L}', \rho)$ and $(\mc{L},\tau)$ are locally isometric at these primes.\\

For the remaining primes, we claim that both lattices $\mc{L}_p$ and $ \mc{L}_p'$ are unimodular. First of all,
\begin{equation*}
    \mc{L}_p^{\#} = \mathfrak{d}_{K/\mathbb Q}^{-1}\widehat{\mc{L}}_p = \widehat{\mc{L}}_p ,
\end{equation*}
since $p \nmid |d_K|$. But $\widehat{\mc{L}}_p$ is the dual lattice of $\mc{L}_p$, with respect to $B$. But since $p \nmid \ell(\mc{L}, \tau)$, we obtain $p \nmid \det (Q)$ from \eqref{lev, det, div 1}. Hence, $\mc{L}_p$ is locally unimodular with respect to $B$, meaning $\widehat{\mc{L}}_p = \mc{L}_p$. Hence, $ \mc{L}_p^{\#} = \mc{L}_p$. Similarly for $(\mc{L}_p', \rho)$ since $\ell(\mc{L},\tau) = \ell(\mc{L}', \rho)$.\\

But from \cite[Proposition 3.3.5]{thesis_kirschmer}, we obtain that there is only one class of local unimodular lattices when $p$ is not a ramified prime. Hence, $(\mc{L}_p, \tau)$ and $(\mc{L}_p', \rho)$ are isometric. Hence, they are isometric at every finite prime, meaning they lie in the same $U(\tau)$-genus.
\end{proof}

\section{Hecke rings of the general linear group}\label{Hecke rings of the general linear group}
Our next goal is to compute the action of the global Hecke algebra $\mathcal{H}$ on the theta series $\theta_{\mathbb{A}}$. Our strategy is similar to that of \cite[Section 5.1.2]{andrianov_book}. In particular, we will first compute the local actions of $\mathcal{H}_p$, and then ``glue'' them together. The case of inert primes is quite similar to the situation encountered by Andrianov. Our main focus will therefore be on the case of split primes, where the unitary group is locally related to a general linear group. In this Section, we study Hecke rings of the general linear group and prove various useful results.\\

Let $n \geq 1$ and $p$ denote any prime. Let $G_n := \textup{GL}_{n}(\mathbb{Q}_p)$ and $\Gamma_n := \textup{GL}_{n}(\mathbb{Z}_p)$. Denote also 
\begin{equation*}
    G_{2n,n} = \left\{\m{a & b \\ 0 & d} \in G_{2n}\right\}, \qquad \Gamma_{2n,n} = G_{2n,n}\cap \Gamma_{2n}.
\end{equation*} 
We are interested in the Hecke rings $\mathcal{H}(\Gamma_n, G_{n})$ and $\mathcal{H}(\Gamma_{2n,n}, G_{2n,n})$. These are related by the following embedding:
\begin{align}\label{embedding gl}
    \iota : \mathcal{H}(\Gamma_{2n}, G_{2n}) &\longrightarrow \mathcal{H}(\Gamma_{2n,n}, G_{2n,n})\\\nonumber
    \sum a_i \Gamma_{2n} g_i &\longmapsto \sum a_i \Gamma_{2n,n}g_i,
\end{align}
This is a well-defined injective ring homomorphism, as seen, for example, in \cite{parabolic_gl_I}, \cite{parabolic_gl_II}.\\

Now, from \cite[Proposition 3.1.7]{andrianov_book}, if $\mathcal{H} = \mathcal{H}(\Gamma, S)$ is any Hecke ring, where $S$ is a group, the map
\begin{align}\label{j map}
    j : \mathcal{H} &\longrightarrow \mathcal{H}\\\nonumber
    \Gamma g \Gamma &\longmapsto \Gamma g^{-1} \Gamma,
\end{align}
is an anti-isomorphism, i.e. an isomorphism which satisisfies $j(ab) = j(b)j(a)$ for all $a, b \in \mathcal{H}$. \\

Moreover, the map 
\begin{align}\label{* map}
    * : \mathcal{H}(\Gamma_{2n,n}, G_{2n,n}) &\longrightarrow \mathcal{H}(\Gamma_{2n,n}, G_{2n,n})\\\nonumber
    \Gamma_{2n,n} g \Gamma_{2n,n} &\longmapsto \Gamma_{2n,n} g^{*} \Gamma_{2n,n}, 
\end{align}
with 
\begin{equation*}
    g^{*} = Ig^{t}I, \qquad I := \m{& & 1 \\ & \iddots & \\ 1 & &},
\end{equation*}
is also an anti-isomorphism (see \cite[Eq. (3)]{parabolic_gl_II}).\\

We now define the following standard Hecke elements of $\mathcal{H}(\Gamma_{n}, G_{n})$:

\begin{equation}\label{generators of gl}
    \pi_{i} = \pi_{i}(p) =\pi_{i,n}(p) = \Gamma_{n}\m{1_{n-i} & \\ & p1_i}\Gamma_{n}, \qquad i=0, \cdots, n.
\end{equation}
For convenience, we shall denote
\begin{equation}\label{A matrices}
A_{i,n} := \textup{diag}(1_{n-i}, p1_i), \qquad i=0,\ldots, n.
\end{equation}

We now define, for $i=0,\ldots, n$ the elements
\begin{equation*}
    \pi_{i,n}^{-} = \Gamma_{2n,n}\textup{diag}(1_{n-i}, p1_{i}, 1_{n})\Gamma_{2n,n}, \qquad \pi_{i,n}^{+} = \Gamma_{2n,n}\textup{diag}(1_{n}, 1_{n-i}, p1_{i})\Gamma_{2n,n}.
\end{equation*}

We now record the following Lemma of \cite[Lemma 2]{parabolic_gl_I}, regarding right coset decomposition of elements in $\mathcal{H}(\Gamma_{2n,n}, G_{2n,n})$.  Let $n \geq 1$. For a given matrix $R \in \textup{GL}_{n}(\mathbb{Q}_p)$, we define $$\Gamma_{n}^{R} := \Gamma_{n} \cap R^{-1}\Gamma_{n}R.$$ Also, if $A,D \in \textup{GL}_n(\mathbb{Q}_p)$, we define
\begin{equation*}
     V(A,D) := \{ AY \,\,|\,\, Y \in M_{n}(\mathbb{Z}_p) \pmod{D^{\times}} \} / \sim,
\end{equation*}
where $AY_1 \sim AY_2 \pmod{D^{\times}}$ if and only if  $AY_1D^{-1}-AY_2D^{-1} \in M_{n}(\mathbb{Z}_p)$. We then have
\begin{lemma}[Gritsenko]\label{gritsenko decomposition}
We have 
\begin{equation*}
\Gamma_{2n,n}\m{A & \\ & D}\Gamma_{2n,n} = \sum_{\mu,\nu, B} \Gamma_{2n,n}\m{A & B\\ &D}\m{\mu & \\ & \nu}
\end{equation*}
where $\mu \in \Gamma_{n}^{A} \backslash \Gamma_{n}$, $\nu \in \Gamma_{n}^{D} \backslash \Gamma_{n}$ and $B \in V(A,D)$.
\end{lemma}
\begin{remark}
    In \cite[Lemma 2]{parabolic_gl_I}, the Lemma is given for the localised rings. However, it extends in the obvious way in our completed setting, as can be seen, for example, from \cite[Proposition 3.1]{heyer2}.
\end{remark}

We can now prove the following Lemma regarding the image of the standard elements of $\mathcal{H}(\Gamma_{2n}, G_{2n})$ under the map $\iota$ of \eqref{embedding gl}.
\begin{lemma}[Embeddings in the parabolic Hecke algebra]\label{Embeddings in the parabolic Hecke algebra}
    We have 
    \begin{equation*}
        \iota(\pi_{r,2n}) = \sum_{\substack{a+b=r\\0 \leq a,b \leq n}}p^{-ab}\pi_{a,n}^{-}\pi_{b,n}^{+}, \qquad 0 \leq r \leq 2n.
    \end{equation*}
\end{lemma}
\begin{proof}
    For the proof, we suppress the subscript denoting the size of the general linear group, as it follows from context. We will first show
    \begin{equation*}
        \iota(\pi_{r}) = \sum_{\substack{a+b=r\\0 \leq a, b \leq n}} \Gamma_{2n,n} \m{A_a & \\ & A_b}\Gamma_{2n,n}.
    \end{equation*}
    Recall, 
    \begin{equation*}
        \pi_r = \Gamma_{2n} \m{1_{2n-r} & \\ & p1_{r}} \Gamma_{2n}.
    \end{equation*}
    We will first show that as sets, 
    \begin{equation}\label{intermediate step}
        \pi_r \cap G_{2n,n} = \bigsqcup_{\substack{a+b=r\\a, b \leq n}} \Gamma_{2n,n} \m{A_a & \\ & A_b}\Gamma_{2n,n}.
    \end{equation}
    Take any
    \begin{equation*}
        g = \m{A & X \\ 0 & D}
    \end{equation*}
    on the left-hand side. Since $g \in \pi_{r}$, $g$ and $pg^{-1}$ have integral entries. Hence $A, D, pA^{-1}, pD^{-1}$ also have integral entries, which shows that the elementary divisors of both $A,D$ are either $1$ or $p$. Writing $a = v_p (\det A)$ and $b = v_p (\det D)$, we see that $0 \leq a, b \leq n $ and $a+b = v_{p}(\det g) = r$.   \\

    Now, using row/column operations, we may find $\gamma_1, \gamma_2 \in \Gamma_{2n,n}$ such that $$\gamma_1 g \gamma_2 = \m{A_a & X' \\ & A_b}.$$ Using now again the fact that $pg^{-1}$ has integral entries, this means $-pA_a^{-1}X'A_b^{-1}$ has integral entries. This forces the bottom right $a \times b $ submatrix of $X'$ to be divisible by $p$. Now,
    \begin{equation*}
        \m{1_n & Y \\ & 1_n} \m{A_a & X' \\ & A_b} \m{1_n & Z \\ & 1_n} = \m{A_a & A_a Z+ X' + YA_b \\ & A_b}
    \end{equation*}
    and we can choose integral $Y,Z$ so that the top right corner becomes zero. Therefore, $$g \in \Gamma_{2n,n}\textup{diag}(A_a, A_b) \Gamma_{2n,n},$$ hence contained in the right-hand side of \eqref{intermediate step}. On the other hand, if $g \in \Gamma_{2n,n}\textup{diag}(A_a, A_b) \Gamma_{2n,n}$, for some $0 \leq a,b \leq n$ with $a+b = r$, then clearly $g \in \pi_{r}$, and also $g \in G_{2n,n}$, since $\textup{diag}(A_a, A_b) \in G_{2n,n}$. Moreover, the valuations of the determinants of the two diagonal blocks are invariant by multiplication on either side by $\Gamma_{2n,n}$. Thus, distinct pairs $(a,b)$ give disjoint double cosets, and \eqref{intermediate step} is established. Let now
    \begin{equation*}
        \pi_r = \sum_{i} \Gamma_{2n} g_i,
    \end{equation*}
    for some fixed representatives $g_i \in \pi_{r} \cap G_{2n,n}$. This is possible since $G_{2n}=\Gamma_{2n}G_{2n,n}$. Then 
    \begin{equation*}
        \iota(\pi_r) = \sum_{i} \Gamma_{2n,n} g_i.
    \end{equation*}
    Note that if $g, h$ come from the same $\Gamma_{2n}$ coset, then $g = \gamma h$, which implies $\gamma \in \Gamma_{2n} \cap G_{2n,n} = \Gamma_{2n,n}$. Conversely, equality of the $\Gamma_{2n,n}$ cosets implies equality of the $\Gamma_{2n}$ cosets. Thus, the cosets in $\iota(\pi_r)$ are all distinct and exhaust the cosets in $\pi_r \cap G_{2n,n}$. Using \eqref{intermediate step}, we may therefore rewrite
    \begin{multline*}
        \iota(\pi_r) = \sum_{\substack{\textup{distinct left cosets } \Gamma_{2n,n}g\\g \in \pi_r \cap G_{2n,n}}} \Gamma_{2n,n}g = \sum_{\substack{a+b=r\\0 \leq a,b \leq n}} \sum_{\substack{\textup{distinct left cosets } \Gamma_{2n,n}g\\g \in \Gamma_{2n,n}\textup{diag}(A_a, A_b)\Gamma_{2n,n}}}\Gamma_{2n,n}g 
        \\ =\sum_{\substack{a+b=r\\0 \leq a,b \leq n}} \Gamma_{2n,n}\textup{diag}(A_a, A_b)\Gamma_{2n,n},
    \end{multline*}
    with each double coset on the right-hand side appearing exactly once.\\
    
    Fix now $a,b$ with $0 \leq a,b \leq n$. We choose fixed representative sets for $\Gamma_n^{A_a}\backslash\Gamma_n$ and $\Gamma_n^{A_b}\backslash\Gamma_n$. In what follows, $X,L$ range over the first set and $Z,M$ over the second. We also take $K \in V(A_a,A_b)$ and $Y \in V(1_n,A_b)$. Using Lemma \ref{gritsenko decomposition}, we have
    \begin{equation*}
        \Gamma_{2n,n} \m{A_a & \\ & A_b}\Gamma_{2n,n} = \sum_{K,L,M} \Gamma_{2n,n} \m{A_aL & KM\\ & A_b M}.
    \end{equation*}
    \begin{equation*}
        \pi_{a}^{-} = \Gamma_{2n,n} \m{A_a & \\ & 1_n}\Gamma_{2n,n} = \sum_{X} \Gamma_{2n,n} \m{A_aX & \\ & 1_n}.
    \end{equation*}
    \begin{equation*}
        \pi_{b}^{+} = \Gamma_{2n,n} \m{1_n & \\ & A_b}\Gamma_{2n,n} = \sum_{Y,Z} \Gamma_{2n,n} \m{1_n & YZ\\ & A_bZ}.
    \end{equation*}
    Hence,
    \begin{equation*}
        \pi_{a}^{-} \pi_{b}^{+} = \sum_{X,Y,Z} \Gamma_{2n,n} \m{A_aX & (A_a X)YZ\\ & A_bZ}.
    \end{equation*}
    We now check when 
    \begin{equation*}
        \Gamma_{2n,n} \m{A_aX & (A_a X)YZ\\ & A_bZ} = \Gamma_{2n,n} \m{A_aL & KM\\ & A_b M}.
    \end{equation*}
    Comparing the diagonal blocks shows that
    \begin{equation*}
        \Gamma_n^{A_a}X=\Gamma_n^{A_a}L,
        \qquad
        \Gamma_n^{A_b}Z=\Gamma_n^{A_b}M.
    \end{equation*}
    By our choice of representatives, this means $X=L$ and $Z=M$. Looking at the top right corner, we obtain
    \begin{equation*}
        (A_a XY-K)A_b^{-1} \in M_n(\mathbb{Z}_p).
    \end{equation*}
    Now, since $K \in V(A_a, A_b)$, it has the form $K = A_a \widetilde{K}$, and hence the relation is equivalent to 
    \begin{equation*}
        [A_a XY] = [K] \textup{ in } V(A_a, A_b).
    \end{equation*}
    Now, for fixed $K, X$, we need to find how many $Y \in V(1_n, A_b)$ exist, so that $[A_a X Y] = [K]$ in $V(A_a, A_b)$. This is possible by considering the map
    \begin{align*}
        V(1_n, A_b) &\longrightarrow V(A_a, A_b)\\
        [Y] &\longmapsto [A_a X Y].
    \end{align*}
    This is a well-defined homomorphism, since $A_aX$ is integral. It is also surjective, since $X \in \Gamma_n$ and therefore $XM_n(\mathbb Z_p)=M_n(\mathbb Z_p)$. Modulo $M_n(\mathbb Z_p)A_b$, the set $V(1_n,A_b)$ is determined by the last $b$ columns modulo $p$, while $V(A_a,A_b)$ is determined by the first $n-a$ entries in each of these columns. Thus,
    \begin{equation*}
        |V(1_n,A_b)|=p^{nb},
        \qquad
        |V(A_a,A_b)|=p^{(n-a)b}.
    \end{equation*}
    Each fibre therefore has size
    \begin{equation*}
        \frac{|V(1_n,A_b)|}{|V(A_a,A_b)|}=p^{ab}.
    \end{equation*}
    It follows that
    \begin{equation*}
        \pi_a^{-}\pi_b^{+}
        = p^{ab}\Gamma_{2n,n}\m{A_a & \\ & A_b}\Gamma_{2n,n}.
    \end{equation*}
    Combining this with the first part of the proof gives
    \begin{equation*}
        \iota(\pi_r)
        = \sum_{\substack{a+b=r\\0 \leq a,b \leq n}}
        p^{-ab}\pi_a^{-}\pi_b^{+}. \qedhere
    \end{equation*}
\end{proof}

Using this Lemma, it is possible to obtain a factorisation of the standard Hecke polynomials of $\mathcal{H}( \Gamma_{2n}, G_{2n})$. We use the notation $\langle i\rangle := i(i+1)/2$.
\begin{defn}\label{standard polynomials gl}
    We define the following standard Hecke polynomials (cf. \cite[Example 2]{parabolic_gl_II})
    \begin{equation*}
        R_{1,p}^{(2n)}(v) = \sum_{i=0}^{2n} (-1)^{i}p^{\langle i-1\rangle } \pi_{i,2n} v^{i}, \qquad  R_{2,p}^{(2n)}(v) = \sum_{i=0}^{2n} (-1)^{i}p^{\langle i-1\rangle } j(\pi_{i,2n}) v^{i},
    \end{equation*}
    where $j$ is the anti-isomorphism of \eqref{j map}.
\end{defn}
We then have the following Proposition.
\begin{proposition}\label{factorisation gl}
    We have the following factorisations in the ring $\mathcal{H}(\Gamma_{2n,n}, G_{2n,n})[v]$ (i.e. we identify the polynomials $R$ with their images in the parabolic Hecke ring):
    \begin{equation*}
        R_{1,p}^{(2n)}(v) = \left(\sum_{i=0}^n (-1)^i p^{\langle i-1 \rangle} \pi_{i,n}^{-} v^i\right)\left(\sum_{i=0}^n (-1)^i p^{\langle i-1 \rangle} \pi_{i,n}^{+} v^i\right).
    \end{equation*}
    \begin{equation*}
        R_{2,p}^{(2n)}(v) = \left(\sum_{i=0}^n (-1)^i p^{\langle i-1 \rangle} j(\pi_{i,n}^{+}) v^i\right)\left(\sum_{i=0}^n (-1)^i p^{\langle i-1 \rangle} j(\pi_{i,n}^{-}) v^i\right).
    \end{equation*}
\end{proposition}
\begin{proof}
    The first factorisation follows directly from Lemma \ref{Embeddings in the parabolic Hecke algebra}. The second factorisation follows after applying the anti-isomorphism $j$ to the first, since the maps $\iota$ and $j$ commute from \cite[Lemma 3.1.9]{andrianov_book}.
\end{proof}
We end this Section, by recalling the spherical (Satake) mapping. First, it is well-known that any element of $g \in \mathcal{H}(\Gamma_{2n}, G_{2n})$ can be written as $X = \sum_{i} \Gamma_{2n}g_i$ with
\begin{equation*}
    g_i = \m{p^{d_{i1}}&*&* \\& \ddots & * \\ & &  p^{d_{i2n}}}.
\end{equation*}
\begin{defn}[Spherical Mapping, \cite{parabolic_gl_I}]\label{spherical map}
For $X$ as above, the map
        \begin{equation}\label{spherical mapping global ring}
            \Phi_{2n}(X) := \sum_{i}a_i  \Phi_{2n}(\Gamma_{2n}g_i) = \sum_{i} a_i \prod_{j=1}^{2n} (x_jp^{-j})^{d_{ij}},
        \end{equation}
        defines an isomorphism of $\mathcal{H}(\Gamma_{2n}, G_{2n})$ with the ring of symmetric polynomials $\mathbb{Q}^{\Sigma}[x_1^{\pm}, \cdots, x_{2n}^{\pm}]$.\\
        
Additionally, if $X = \sum_{i}a_i \Gamma_{2n,n}\m{A_i & B_i \\ 0 & D_i}$, we use $x_i$ for the Satake variables associated with the $A$-block and $y_i$ for the ones associated with the $D$-block. We then set
\begin{equation}\label{omega mapping}
    \Omega(X) := \sum_{i} a_i (\Gamma_n A_i)\Phi_{n}^{(y)}(\Gamma_nD_i) \in \mathcal{H}(\Gamma_n, G_n)[y_1^{\pm}, \ldots, y_n^{\pm}].
\end{equation}
Then, we set $\Phi_{2n,n}(X) := \Phi_{n}^{(x)}(\Omega(X))$, where $\Phi^{(n)}$ is extended $\mb Q[y_1^{\pm}, \ldots, y_n^{\pm}]$- linearly. This is a homomorphism from $\mathcal{H}(\Gamma_{2n,n}, G_{2n,n})$ to $\mathbb{Q}[x_1^{\pm}, \ldots, x_{n}^{\pm}, y_1^{\pm}, \ldots, y_n^{\pm}]$. 
\end{defn}
We have the following property of $\Omega$. Let
\begin{equation}
    \mathcal{C}^{-} := \{X \in \mathcal{H}(\Gamma_{2n,n}, G_{2n,n}) \mid X\pi_{n,n}^{-} = \pi_{n,n}^{-} X\},
\end{equation}
denote the centraliser of the element $\pi_{n,n}^{-}$. We then have the following.
\begin{lemma}
    The restriction of $\Phi_{2n,n}$ to $\mathcal{C}^{-}$ is injective.
\end{lemma}
\begin{proof}
    This follows from \cite[Lemma 3]{parabolic_gl_II}. 
\end{proof}
\section{Split-prime local unitary Hecke rings}\label{Split-prime local unitary Hecke rings}
Assume $p\mathcal{O}_K = \mathfrak{p} \overline{\mathfrak{p}}$, with $\mathfrak{p} \neq \ol{\mathfrak{p}}$ and assume $p \nmid \mf{bc}=\ell d_K$. We fix the ordering of the two primes as well as an isomorphism
\begin{equation*}
    K_p := K \otimes_{\mathbb Q}\mathbb Q_p
\xrightarrow{\ \sim\ }\mathbb Q_p\times\mathbb Q_p,
\end{equation*}
where the first factor corresponds to $\mf{p}$ and the other to $\ol{\mf{p}}$.  We choose the uniformisers 
\begin{equation}\label{uniformisers}
    \varpi_{\mathfrak p}=(p,1), \qquad \varpi_{\overline{\mathfrak p}}=(1,p),
\end{equation}
under the identification $K \otimes_{\mathbb Q}\mathbb Q_p
\xrightarrow{\ \sim\ }\mathbb Q_p\times\mathbb Q_p$. We recall
\begin{equation*}
    \textup{GU}(n,n)(\mathbb{Q}_p) = \{g \in \textup{GL}_{2n}(K_p) \mid \ol{g}^{t}J_n g = \mu(g)J_n, \,\ \mu(g) \in \mathbb{Q}_p^{\times}\}.
\end{equation*}
Writing $g=(g_1,g_2)$, and using the fact that conjugation interchanges the two factors, the defining relation becomes
\[
g_2^tJ_ng_1=\mu(g)J_n.
\]
Projection to the first factor gives an isomorphism
\begin{equation*}
\varphi: \textup{GU}(n,n)(\mathbb Q_p)\xrightarrow{\ \sim\ }\mathrm{GL}_{2n}(\mathbb Q_p) \times \mathbb{Q}_p^{\times},
\qquad
(g_1, g_2) \longmapsto (g_1, \mu(g)),
\end{equation*}
whose inverse is
\[
(h, \nu)\longmapsto \bigl(h,-\nu J_n(g^{-1})^tJ_n\bigr).
\]
Under $\varphi$, we have that $\textup{GU}(n,n)(\mathbb{Z}_p)$ gets mapped to $\textup{GL}_{2n}(\mathbb{Z}_p) \times \mathbb{Z}_p^{\times}$, and moreover, the standard Siegel parabolic gets mapped to the Siegel parabolic of $\textup{GL}_{2n}$. Thus, with the notation of Section \ref{preliminaries},  $\varphi$ identifies
$$\mathcal{H}(\mbf{D}_{p},  \mbf{G}_p) \cong \mathcal{H}(\textup{GL}_{2n}(\mathbb{Z}_p), \textup{GL}_{2n}(\mathbb{Q}_p))[x^{\pm}],$$
as well as the corresponding parabolic Hecke rings $$\mathcal{H}(\mbf{D}_{0,p}, \mbf{P}_p) \cong \mathcal{H}(\textup{GL}_{2n,n}(\mathbb{Z}_p), \textup{GL}_{2n,n}(\mathbb{Q}_p))[x^{\pm}].$$

Here, $x$ corresponds to the central double coset of \(\varpi_{\overline{\mathfrak p}}1_{2n}\), whose image under \(\varphi\) is \((1_{2n},p)\) and whose similitude multiplier is \(p\).
\medskip

Let now $\mathfrak{a}$ denote a fractional ideal of $K$. We take an element $\widehat{\mf{a}} \in \prod_{v \textup{ finite}} K_{v}^{\times}$, such that $\widehat{\mf{a}}_v\mathcal{O}_{K_v} = \mf{a}_v$ for each finite place $v$ of $K$. We have defined (see \eqref{t minus original})
\begin{equation*}
    t^{-}(\mf{a}) = \sum_{\substack{D \in E\backslash B/E\\ (\det D) \widehat{\mathcal{O}}_K=\widehat{\mf{a}}\widehat{\mathcal{O}}_K}} \mbf{D}_0 \m{D & \\ & \widehat{D}}\mbf{D}_0,
\end{equation*}
where $B$ as in \eqref{t minus original} and $E$ as in \eqref{R and E}. We now consider the element $t^{-}(\mf{p}^r)$ for $r \geq 0$. We observe that 
\begin{equation*}
    t^{-}(\mf{p}^r)_p = \sum_{\substack{D \in E_p\backslash B_{p}/E_p\\(\det D)\mathcal{O}_{K,p} = \varpi^{r}\mc{O}_{K,p}}}\mbf{D}_{0,p}\m{D & \\ & \widehat{D}} \mbf{D}_{0,p},
\end{equation*}
\begin{equation*}
    t^{+}(\mf{p}^r)_p = \sum_{\substack{D \in E_p\backslash B_{p}/E_p\\(\det D)\mathcal{O}_{K,p} = \varpi^{r}\mc{O}_{K,p}}}\mbf{D}_{0,p}\m{\widehat{D} & \\ & D} \mbf{D}_{0,p},
\end{equation*}
where now $E_p = \textup{GL}_n(\mathcal{O}_{K,p})$ and $B_p = \textup{GL}_n(K_p) \cap M_n(\mathcal{O}_{K,p})$. Note that since $p \nmid \mf{bc}$, $\mbf{D}_{0,p} = \mbf{P}(\mathbb{Z}_p)$.
From the elementary divisor theorem, we may choose representatives for $D$ to be
\begin{equation*}
\left\{\textup{diag}(\varpi_{\mf{p}}^{r_1}, \ldots, \varpi_{\mf{p}}^{r_n}), \,\,\, 0 \leq r_1 \leq \ldots \leq r_n, \,\,\,\,\, \sum_{i=1}^{n}r_i = r,\right\}
\end{equation*}
with the notation of \eqref{uniformisers}. For convenience, let
\[
\mathcal A_r
:=
\left\{
\operatorname{diag}(p^{r_1},\ldots,p^{r_n})
\;\middle|\;
0\leq r_1\leq\cdots\leq r_n,\quad
\sum_{i=1}^n r_i=r
\right\}.
\]
\(\mathcal A_r\) is a complete set of representatives for the integral \(\Gamma_n\)-double cosets whose determinant generates \(p^r\mathbb Z_p\). All the matrices occurring in \(t^{\pm}(\mathfrak p^r)_p\) have similitude multiplier \(1\). Hence, under \(\varphi\), their images lie in
\[
\mathrm{GL}_{2n,n}(\mathbb Q_p)\times\{1\}.
\]
We shall suppress the second component in what follows. Recalling that conjugation exchanges the components, under $\vphi$, the elements $t^{\pm}(\mathfrak{p}^r)$ therefore get mapped to 
\begin{equation*}
    t^{-}_{1}(p^r) := \sum_{A \in \mathcal A_r}\Gamma_{2n,n}\m{A & 0 \\ 0& 1} \Gamma_{2n,n},
    \qquad
    t^{+}_{1}(p^r) := \sum_{A \in \mathcal A_r}\Gamma_{2n,n}\m{1 & 0 \\ 0& A} \Gamma_{2n,n}.
\end{equation*}
Analogously, and using the element $\varpi_{\ol{\mf{p}}}$, we have that $t^{\pm}(\overline{\mf{p}}^r)$ get mapped to
\begin{equation*}
    t^{-}_{2}(p^r) := \sum_{A \in \mathcal A_r}\Gamma_{2n,n}\m{1 & 0 \\ 0& A^{-1}} \Gamma_{2n,n}, \qquad t^{+}_{2}(p^r) := \sum_{A \in \mathcal A_r}\Gamma_{2n,n}\m{A^{-1} & 0 \\ 0& 1}\Gamma_{2n,n}.
\end{equation*}
Note that
\begin{equation}\label{relations 1 2}
    t_{2}^{-}(p^r) = j(t_{1}^{+}(p^r)), \qquad t_{2}^{+}(p^{r}) = j(t_1^{-}(p^{r})).
\end{equation}
We now have the following Lemma.
\begin{lemma}\label{series gl}
    We have 
    \[
        \sum_{r=0}^{\infty} t^{\pm}_1(p^r) v^r = \left(\sum_{i=0}^n (-1)^i p^{\langle i-1 \rangle} \pi_{i}^{\pm} v^i\right)^{-1}.
    \]
\end{lemma}
\begin{proof}
    Note that it suffices to show the statement only for the minus operators, as then the other statement will follow after applying the anti-isomorphism $*$ of \eqref{* map}, since
    \[
\bigl(t_1^{-}(p^r)\bigr)^*=t_1^{+}(p^r),
\qquad
(\pi_i^{-})^*=\pi_i^{+}.
\]
    
    Using Lemma \ref{gritsenko decomposition}, we obtain
    \begin{equation}\label{decomposition t1 minus}
    t_{1}^{-}(p^r) = \sum_{A \in \mathcal A_r}\sum_{\mu \in \Gamma_n^{A} \backslash \Gamma_n}\Gamma_{2n,n} \m{A & 0 \\ 0 & 1}\m{\mu & 0 \\ 0 & 1}.
\end{equation}
Hence, 
\begin{equation*}
    \Omega(t_{1}^{-}(p^r)) = \sum_{A \in \mathcal A_r}\sum_{\mu \in \Gamma_n^{A} \backslash \Gamma_n} \Gamma_{n}(A\mu) = \sum_{A \in \mathcal A_r} \Gamma_{n}A\Gamma_{n} := t(p^{r}),
\end{equation*}
where $\Omega$ is the map of \eqref{omega mapping}. Recall now $\pi_{i}^{-} = \Gamma_{2n,n}\textup{diag}(A_{i},1_n)\Gamma_{2n,n}$, where $A_i$ are the matrices of \eqref{A matrices}. Hence, using again Lemma \ref{gritsenko decomposition}, we obtain 
\begin{equation*}
    \Omega(\pi_{i}^{-}) = \sum_{\mu \in \Gamma_{n}^{A_i}\backslash \Gamma_{n}}\Gamma_{n}(A_i\mu) = \pi_{i}.
\end{equation*}
We now use the standard Hecke-series identity
\[
\sum_{r=0}^{\infty}t(p^r)v^r
=
\left(
\sum_{i=0}^n
(-1)^ip^{\langle i-1\rangle}\pi_iv^i
\right)^{-1},
\]
as in \cite[Proposition 3.2.22]{andrianov_book}. Since, $\pi_{i}^{-}$ and $t_{1}^{-}(p^r)$ are both in $\mathcal{C}^{-}$, and $\Omega$ is injective there, the Lemma follows.
\end{proof}
\begin{corollary}
    We also have 
    \begin{equation*}
        \sum_{r=0}^{\infty} t^{\pm}_2(p^r) v^r = \left(\sum_{i=0}^n (-1)^i p^{\langle i-1 \rangle} j(\pi_{i}^{\mp}) v^i\right)^{-1},
    \end{equation*}
    which follows after applying the $j$ map in Lemma \ref{series gl} as well as \eqref{relations 1 2}. Note that the plus/minus are exchanged on the left and right-hand sides.
\end{corollary}
We now let $x_{\pm i}(\mf{p}) := \varphi^{-1}(\pi_{i}^{\pm})$ and $x_{\pm i}(\ol{\mf{p}}) := \varphi^{-1}(j(\pi_{i}^{\mp}))$ denote the corresponding operators in $\textup{U}(n,n)$. Let also $Z^{(2n)}(\mf{p};v) := \varphi^{-1}(R_{1,p}^{(2n)}(v))$ and $Z^{(2n)}(\ol{\mf{p}};v) := \varphi^{-1}(R_{2,p}^{(2n)}(v))$ We then have:
\begin{proposition}\label{factorisation un}
    We have 
    \begin{equation*}
       Z_1(v) = Z^{(2n)}(\mf{p};v) = \left(\sum_{i=0}^n (-1)^i p^{\langle i-1 \rangle} x_{-i}(\mf{p}) v^i\right)\left(\sum_{i=0}^n (-1)^i p^{\langle i-1 \rangle} x_{i}(\mf{p}) v^i\right).
    \end{equation*}
    \begin{equation*}
        Z_2(v) = Z^{(2n)}(\ol{\mf{p}};v) = \left(\sum_{i=0}^n (-1)^i p^{\langle i-1 \rangle} x_{-i}(\ol{\mf{p}}) v^i\right)\left(\sum_{i=0}^n (-1)^i p^{\langle i-1 \rangle} x_{i}(\ol{\mf{p}}) v^i\right).
    \end{equation*}
\end{proposition}
\begin{proof}
    This is straightforward from Proposition \ref{factorisation gl}.
\end{proof}
If we now let $X_1^{-}(v) := X_1^{-}(\mf{p}; v)$ and $X_{1}^{+}(v) := X_1^{+}(\mf{p}; v)$ denote the $\pm$ series appearing as factors above, from Lemma \ref{series gl}, we also get
\begin{equation*}
    \sum_{r=0}^{\infty} t^{\pm}(\mf{p}^r) v^r = X_{1}^{\pm}(v)^{-1},
\end{equation*}
and the analogous relation for $\mf{p}$, where now $X_2^{-}(v) := X_2^{-}(\mf{\ol{p}}; v)$ and $X_{2}^{+}(v) =: X_2^{+}(\mf{\ol{p}}; v)$.
\section{Split-prime Hecke action on theta series}
% \begin{defn}
%     The image of the standard $L$-function at $p$ under the Satake mapping is given by (see \cite[Section 20.6]{shimura2000arithmeticity})
%     \begin{equation*}
%         Z_{p}(T) = 
%         \begin{cases}
%             \prod_{i=1}^{n}\left\{(1 -N(\mf{p})^{n-1}x_{p,i}[\mathfrak{p}])(1-N(\mf{p})^nx_{p,i}^{-1}[\mf{p}])\right\}^{-1}, & p\mathcal{O}_K = \mathfrak{p}^{e}\\
%             \prod_{i=1}^{2n}\left\{(1 -N(\mf{p})^{2n}x_{p,i}^{-1}[\mathfrak{p}])(1-N(\overline{\mf{p}})^{-1}x_{p,i}[\overline{\mf{p}}])\right\}^{-1}, & p\mathcal{O}_K = \mf{p}\ol{\mf{p}}.\\
%         \end{cases}
%     \end{equation*}
% \end{defn}

We are now ready to compute the explicit Hecke action on the theta series $\theta_{\mathbb{A}}$ attached to a lattice $\mc{L} = \xi \mc{L}_0$ for the case of split primes. We first have the following Lemma, which is true for all primes $p \nmid \mf{bc}$ (i.e. both split and inert).
\begin{lemma}
    If $(\mathfrak{a}, \mathfrak{bc} \cdot p\mathcal{O}_K)=1$ and $X \in \iota_p(\mathcal{H}_{p}) \subset \mathcal{H}_{0,p}$, where $\iota$ is the embedding \eqref{imbedding}, then 
    \begin{equation*}
        t^{-}(\mathfrak{a}) X = X t^{-}(\mathfrak{a}).
    \end{equation*}
\end{lemma}
\begin{proof}
    Since \((\mathfrak a,p\mathcal O_K)=1\) and \(p\nmid\mathfrak{bc}\), the local component of \(t^{-}(\mathfrak a)\) at \(p\) is the identity. On the other hand, \(X\in\iota(\mc{H}_p)\) has identity components away from \(p\). Thus, \(t^{-}(\mathfrak a)\) and \(X\) belong to different factors of the restricted tensor product defining the global (parabolic) Hecke algebra, and hence commute.
\end{proof}
Let now $p\mathcal{O}_K = \mathfrak{p} \overline{\mathfrak{p}}$, with $\mathfrak{p} \neq \ol{\mathfrak{p}}$ and assume $p \nmid \mf{bc}=\ell d_K$. We first have the following:
\begin{lemma}
    For any $r_1, r_2 \geq 0$, we have 
    \begin{equation*}
        t^{-}(\mf{p}^{r_1}\ol{\mf{p}}^{r_2}) = t^{-}(\mf{p}^{r_1})t^{-}(\ol{\mf{p}}^{r_2}).
    \end{equation*}
\end{lemma}
\begin{proof}
    By passing into the Hecke algebra of $\textup{GL}_{2n}$, we obtain
    \begin{equation*}
        \varphi(t^{-}(\mf{p}^{r_1}\ol{\mf{p}}^{r_2})) = \sum_{A,B} \Gamma_{2n,n} \m{A & \\ & B^{-1}}\Gamma_{2n,n},
    \end{equation*}
    where $A,B$ are elementary divisor matrices with $\det A = r_1$ and $\det B = r_2$. Using Lemma \ref{gritsenko decomposition}, we can write this as
    \begin{equation*}
        \varphi(t^{-}(\mf{p}^{r_1}\ol{\mf{p}}^{r_2})) = \sum_{A,B} \sum_{\mu, \nu} \Gamma_{2n,n} \m{A & \\ & B^{-1}}\m{\mu & \\ & \nu} =  \sum_{A,B} \sum_{\mu, \nu} \Gamma_{2n,n} \m{A\mu & \\ & B^{-1}\nu}.
    \end{equation*}
    Using now Lemma \ref{gritsenko decomposition} for each of $\varphi(t^{-}(\mf{p}^{r_1}))$ and $\varphi(t^{-}(\mf{\ol{p}}^{r_2}))$, the Lemma follows (note that $\mu, \nu$ run over the correct representatives).
\end{proof}
Define 
\begin{equation*}
    I_{p} := \{\mathfrak{a \subseteq \mc{O}}_K \mid \gcd(\mathfrak{a}, p\mathcal{O}_K)=1\}.
\end{equation*}
Let also
\begin{equation}\label{chi defn}
    \chi(\mathfrak{a}) := \psi(\mathfrak{a})^{n}N_{K/\mathbb{Q}}(\mathfrak{a})^{-n/2},
\end{equation}
where $\psi$ is the character of Lemma \ref{relation theta epsilon}. Let $T = T_M \in \mc{H}_p$. If $\mathbf{T} = \iota_{p}(T) \in \mathcal{H}_{0,p}$, from Proposition \ref{relation theta epsilon}, we have 
\begin{align}\label{action on theta}
    (\theta_{\mathbb{A}}^{(n)} \mid \mathbf{T})((q,s); (\mathcal{L}, \tau)) &=\nonumber \sum_{\mathfrak{a} \subseteq \mathcal{O}_K} \chi(\mathfrak{a})(\epsilon_{\mathbb{A}} \mid t^{-}(\mathfrak{a}) \mathbf{T}) ((q,s); (\mathcal{L}, \tau))\\\nonumber
    &= \sum_{r=0}^{\infty}\sum_{s=0}^{\infty} \sum_{\mathfrak{a} \subseteq I_{p}} \chi(\mathfrak{a}\mathfrak{p}^r\ol{\mf{p}}^s)  (\epsilon_{\mathbb{A}} \mid t^{-}(\mathfrak{a}\mathfrak{p}^r\ol{\mf{p}}^s) \mathbf{T}) ((q,s); (\mathcal{L}, \tau))\\
    &= \sum_{\mathfrak{a} \subseteq I_{p}}\sum_{r=0}^{\infty}\sum_{s=0}^{\infty} \epsilon_{\mathbb{A}}\mid t^{-}(\mathfrak{p}^r)t^{-}(\ol{\mf{p}}^{s})\mathbf{T}\chi(\mathfrak{p})^{r}\chi(\ol{\mf{p}}^{s}) \mid \chi(\mf{a})t^{-}(\mf{a}).
\end{align}
In analogy to \cite[(3.5.83)]{andrianov_book}, we have the following definition of the $p$-signature.
\begin{defn}
    Let $$X = \mbf{D}_{0,p}\m{A & B \\ 0 & D}\mbf{D}_{0,p} \in \mc{H}_{0,p}.$$ Then, $A,D \in \textup{GL}_n(K_p)$, and we write $A = (A_1, A_2)$ and $D = (D_1, D_2)$ under the identification $K_p = \mathbb{Q}_p \times \mathbb Q_p$. We then define the $p$-signature of $X$ as
    \begin{equation*}
        \sigma_p(X) := v_{p}(\det D_1) - v_{p}(\det A_1) = v_{p}(\det D_2) - v_{p}(\det A_2).
    \end{equation*}
\end{defn}
This is well defined since left or right multiplication by an element of $\mbf{D}_{0,p}$ multiplies the determinants of the diagonal blocks by $p$-adic units. Moreover, $A = \mu(D^{*})^{-1}$ which gives that the two differences are equal.\\

For an element $X = \sum_{i}a_i\mbf{D}_{0,p}M_i \mbf{D}_{0,p}\in \mc{H}_{0,p}$, we set
\begin{equation}\label{p homoegeneous}
    X(\sigma) := \sum_{i, \sigma_p(M_i)=\sigma} a_i  \mbf{D}_{0,p}M_i \mbf{D}_{0,p},
\end{equation}
and we call $X(\sigma)$ the $p$-homogeneous component of $X$ of $p$-signature $\sigma$.

\begin{lemma}\label{pos signature}
    Assume $X$ is $p$-homogeneous. Then if $\sigma_p(X) >0$ and $v_{p}(\vol(\mc{L}, \tau))= 0$, we have
    \begin{equation*}
        \epsilon_{\mathbb{A}}((q,s); (\mathcal{L}, \tau)) \mid X  = 0.
    \end{equation*}
\end{lemma}
\begin{proof}
    Take a double coset appearing in $X$ and let $D$ denote its lower-right part. We claim that under the assumptions of the Lemma, 
    \begin{equation*}
        \mu_p\tau[\xi \widehat{R}] \not \in \widehat{\mc{S}}_p,
    \end{equation*}
    if $R \in EDE$ and $\widehat{\mc{S}}$ is as in \eqref{T set}. The result then follows from \eqref{coefficients} because the relevant Gauss sums will vanish from Lemma \ref{non-vanishing lemma}. First, note that since $p$ is a split prime, we have that 
    \begin{equation*}
        T_p = \{(A,A^{t}) \mid A \in M_n(\mathbb{Z}_p)\}.
    \end{equation*}
    Writing $R = (R_1, R_2)$, $\xi = (\xi_1, \xi_2)$, we have that
    \begin{equation*}
        \mu_p \tau[\xi\widehat{R}] = (\mu R_1^{-1}\xi_2^{t}\tau_{\mf{p}} \xi_{1} (R_2^{-1})^{t}, \mu R_2^{-1}\xi_1^{t}\tau_{\mf{p}}^{t} \xi_2 (R_1^{-1})^{t}).
    \end{equation*}
    Now, the assumption on the signature means that
    \begin{equation*}
        v_{p}(\det(R_1)) + v_p(\det(R_2)) - nv_p(\mu_p) > 0.
    \end{equation*}
    Additionally, since $v_{p}(\vol(\mc{L}))= 0$, this means 
    \begin{equation*}
        v_{p} (\det (\tau)\det(\xi_1)\det(\xi_2)) \leq 0.
    \end{equation*}
    Hence, by looking at the $p$-valuation of the determinant of the first component of $\mu_p \tau[\xi \widehat{R}]$, we see that it is negative, which proves the Lemma.
\end{proof}
Going back to \eqref{action on theta}, we will now consider the series
\begin{equation}\label{start of epsilon action}
    \epsilon_{\mathbb{A}} \mid \sum_{r=0}^{\infty}\sum_{s=0}^{\infty} t^{-}(\mathfrak{p}^r)t^{-}(\ol{\mf{p}}^{s})\mathbf{T}v^{r}w^{s}.
\end{equation}
We let 
\begin{equation}\label{Y polynomials}
    Y^{-}(v,w) := X_{1}^{-}(v)X^{-}_{2}(w), \qquad Y^{+}(v,w) := X_1^{+}(v)X_2^{+}(w),
\end{equation}
where $X_i^{\pm}$ are the factors of Proposition \ref{factorisation un}. Note that $X_i^{-}, i=1,2$ commute and the same is true for $X_i^{+}, i=1,2$. We then have
\begin{equation}\label{next expression}
    \sum_{r=0}^{\infty}\sum_{s=0}^{\infty} t^{-}(\mathfrak{p}^r)t^{-}(\ol{\mf{p}}^{s})v^{r}w^{s}\mathbf{T} = Y^{-}(v,w)^{-1}\mbf{T}  = Q(v,w) Y^{-}(v,w)^{-1},
\end{equation}
where
\begin{equation*}
    Q(v,w) := Y^{-}(v,w)^{-1} \mbf{T} Y^{-}(v,w).
\end{equation*}

Using the factorisation of the standard polynomials $Z_i = X_i^{-}X_i^{+}, i=1,2$ from Proposition \ref{factorisation un}, we obtain
%Assume $R_1R_2 = X_1^-X_2^- %X_2^+X_1^+$. Then $W$ reads
%\[
%(X_1^-X_2^-)^{-1} T X_{1}^- X_2^- = %X_1^+X_2^+ (R_1 R_2)^{-1}T R_1R_2 %(X_1^+ X_2^+)^{-1}
%\]
\[
Z_1(v)Z_2(w) = X_1^-(v)X_1^+(v) X_2^-(w) X_2^+(w) = X_1^-(v) X_2^-(w) B(v,w) X_1^+(v)X_2^+(w), 
\]
where
\begin{equation}\label{B polymomial}
    B(v,w) := X_2^-(w)^{-1} X_1^{+}(v) X_2^{-}(w) X_1^+(v)^{-1} \in \mc{H}_{0,p}[[v,w]]
\end{equation}

\begin{lemma}
   $B(v,w)$ is a polynomial in $vw$, whose coefficients all have zero signature.
\end{lemma}
\begin{proof}
For convenience, we set
\begin{equation*}
    K = K(v) := X_1^{-}(v), \,\,\, L = L(v) := X_1^{+}(v), \,\,\, M = M(w) := X_2^{-}(w), \,\,\, N = N(w) := X_2^{+}(w).
\end{equation*}
Then, we have defined 
\begin{equation}\label{B first}
    B(v,w) := M^{-1}(w)L(v)M(w)L^{-1}(v).
\end{equation}
But $$K(v)L(v)M(w)N(w) = Z_1(v) Z_2(w)\implies L(v) M(w) = K^{-1}(v)Z_1(v)Z_2(w)N^{-1}(w)$$ and $Z_1(v), Z_2(w)$ commute since they belong in $\iota(\mc{H}_p)$. Hence (we omit the variables $v,w$)
\begin{multline}\label{B second}
    B(v,w) = M^{-1}K^{-1}R_1R_2N^{-1}L^{-1} = M^{-1}K^{-1}R_2R_1 N^{-1}L^{-1}\\ = M^{-1}K^{-1}MNKLN^{-1}L^{-1} = K^{-1}NKN^{-1},
\end{multline}
using also the fact that $M^{-1}, K^{-1}$ and $N^{-1}, L^{-1}$ commute (since they belong in the minus and plus spaces respectively). Let now 
\begin{equation*}
    B(v,w) = \sum_{r, s \geq 0}b_{r,s}v^r w^s.
\end{equation*}
We now consider the signature of the coefficient of $v^r w^s$. We first look at \eqref{B first}. We observe that we need to take the coefficients of $w^{a}, v^{b}, w^{c}, v^{d}$ from $M^{-1}, L, M, L^{-1}$ respectively such that $a+c = s$ and $b+d = r$. Note that since $L,M$ are polynomials of degree $n$, we will have $0 \leq b,c \leq n$. Moreover, we observe that the coefficient of $w^{a}$ has signature $-a$, the coefficient of $v^{b}$ signature $b$, of $w^{c}$ signature $-c$ and of $v^{d}$ signature $d$. Hence, the signature of $b_{r,s}$ will be 
\begin{equation*}
    -a+b-c+d = (b+d) - (a+c) = r-s.
\end{equation*}
We now look at \eqref{B second}.  We now take the coefficients of $v^{a}, w^{b}, v^{c}, w^{d}$ from $K^{-1}, N, K, N^{-1}$ respectively such that $a+c = r$ and $b+d = s$. Note that since $K,N$ are polynomials of degree $n$, we will have $0 \leq b,c \leq n$. Moreover, we observe that the coefficient of $v^{a}$ has signature $-a$, the coefficient of $w^{b}$ signature $b$, of $v^{c}$ signature $-c$ and of $w^{d}$ signature $d$. Hence, the signature of $b_{r,s}$ will be 
\begin{equation*}
    -a+b-c+d = (b+d)-(a+c) = s - r.
\end{equation*}
This shows that if $r \neq s$, then $b_{r,s}=0$. Hence, we may write
\begin{equation*}
    B(v,w) = \sum_{t=0}^{\infty} B_{t} (vw)^{t},
\end{equation*}
where the coefficients $B_{t}$ all have signature $0$. We will now show that this series, in fact, terminates. Observe that from \eqref{B first}
\begin{equation}\label{expression for B}
    B(v,w) L(v) = M^{-1}(w)L(v) M(w).
\end{equation}
We now look at the coefficient of $(vw)^{m}$ for $m > n$ in both sides. On the right-hand side, this is zero, since the maximum power of $v$ appearing is at most $n$ (since $L$ is a polynomial of degree $n$ in $v$ and $M$ is a polynomial in $w$). On the left-hand side, 
\begin{equation*}
    B(v,w)L(v) = \sum_{t=0}^{\infty}\sum_{i=0} ^{n}B_t\ell_iv^{t+i}w^{t}.
\end{equation*}
So we need $t=t+i=m$, which implies $i=0$ and $t=m$. Hence, the coefficient is exactly $B_{m}$. This shows $B_m = 0$ for all $m > n$ and hence
\begin{equation*}
    B(v,w) = \sum_{t=0}^{n} B_t (vw)^{t},
\end{equation*}
with all the $B_{t}$ having signature equal to zero.
\end{proof}
\begin{lemma}\label{coefficients of B}
    The coefficients $B_t$ are given by
    \begin{equation*}
        B_{t} = (-1)^{t}p^{\langle t-1\rangle}\sum_{j=0}^{t} (-1)^{j}p^{\langle j-1\rangle} t^{-}(\overline{\mf{p}}^{t-j})x_{t}(\mf{p}) x_{-j}(\ol{\mf{p}}), \qquad 0 \leq t \leq n.
    \end{equation*}
\end{lemma}
\begin{proof}
    We use \eqref{expression for B}. As in the previous Lemma, the coefficient of $(vw)^{t}$ on the left-hand side is given by $B_t$. Hence, by equating coefficients,
    \begin{equation*}
        B_{t} = (-1)^{t}p^{\langle t-1\rangle}\sum_{j=0}^{t} (-1)^{j}p^{\langle j-1\rangle} t^{-}(\overline{\mf{p}}^{t-j})x_{t}(\mf{p}) x_{-j}(\ol{\mf{p}}).\qedhere
    \end{equation*}
\end{proof}
The main property of the polynomial $B$ is that it has a scalar action on $\epsilon_{\mb A}$, in analogy with \cite[Theorem 4.2.27]{andrianov_book}.

\begin{proposition}\label{action of B}
    The action of $B(v,w)$ on $\epsilon_{\mathbb{A}}$ is given 
    \[
    (\epsilon_{\mathbb{A}}\mid B(v,w))((q,s);(\mathcal{L};\tau)) = B(v,w,(\mathcal{L}, \tau))\epsilon_{\mathbb{A}}((q,s);(\mathcal{L};\tau)),
    \]
    where 
    \[
    B(v,w,(\mathcal{L},\tau)) = \prod_{i=0}^{r-1}(1-p^{2n-r+i}vw),
\]
with $r:=\textup{rank}_{\mathbb{F}_p}(\overline{\xi^*\tau \xi})_{\mathfrak{p}}$. Here $\xi^*\tau \xi$ is the adelic symmetric matrix corresponding to the lattice $\mathcal{L} = \xi \mathcal{L}_0$, and for a split prime $p\mathcal{O}_K =\mathfrak{p} \overline{\mathfrak{p}}$ we write $(\xi^*\tau \xi)_\mathfrak{p}$ for the projection to the first component at the prime $p$. This is well-defined since the lattice $\mc{L}$ is taken $\mathbb{Z}$-integral.
\end{proposition}
\begin{proof}
    In the following, we omit the notation $(q,s)$ from
    $\epsilon_{\mathbb{A}}$ whenever no confusion can arise.
    We write
    \begin{equation}\label{initial defns}
        c_i:=(-1)^ip^{\langle i-1 \rangle},
        \qquad
        H_i:=t^{-}(\ol{\mf{p}}^{i}),
        \qquad
        V_i:=x_i(\mf{p}),
        \qquad
        W_i:=x_{-i}(\ol{\mf{p}}).
    \end{equation}
    Then, from Lemma \ref{coefficients of B}, we have
    \begin{equation}\label{B three factor expression}
        B_t=c_t\sum_{j=0}^{t}c_jH_{t-j}V_tW_j.
    \end{equation}

    We first compute this expression in the parabolic Hecke
    algebra of $\textup{GL}_{2n}(\mathbb{Q}_p)$, using the
    identification $\varphi$ fixed in Section \ref{Split-prime local unitary Hecke rings}. Put
    \begin{equation}\label{proof definitions}
        \Gamma:=\Gamma_{2n,n}(\mathbb{Z}_p),
        \qquad
        \Lambda_0:=\mathbb{Z}_p^n,
        \qquad
        \mathbf{n}(X):=
        \begin{pmatrix}
            1_n&X\\
            0&1_n
        \end{pmatrix}.
    \end{equation}
    For $0\leq j\leq n$, define
    \begin{equation*}
        \mathcal{U}_j
        :=
        \sum_{\substack{U\subseteq\mathbb{F}_p^n\\
                        \dim U=j}}
        \ \sum_{\substack{Z\in M_n(\mathbb{F}_p)\\
                           \operatorname{im}Z\subseteq U}}
        \Gamma\mathbf{n}(p^{-1}\widetilde Z),
    \end{equation*}
    where $\widetilde Z \in M_n(\mathbb{Z}_p)$ is a lift of $Z$.
    This definition is independent of the lift. 
    % These sums
    % are both left and right $\Gamma$-invariant, since the multiplicity of a
    % matrix $Z$ depends only on its rank.
    In particular,
    $\mathcal{U}_0=1$. We claim that
    \begin{equation}\label{B finite unipotent expansion}
        B_t
        =
        c_t\sum_{j=0}^{t}
        c_jp^{n(t-j)}
        {n-j\brack t-j}_p\,\mathcal{U}_j,
    \end{equation}
    where 
    \begin{equation}\label{gaussian binomial coefficient}
        {n\brack j}_p
        =\prod_{i=0}^{j-1}\frac{p^{n-i}-1}{p^{j-i}-1}.
    \end{equation}
    We prove this by expanding the products in
    \eqref{B three factor expression} into left cosets. A left coset represented by
    \begin{equation*}
        \begin{pmatrix}
            1_n&C\\
            0&D
        \end{pmatrix}
    \end{equation*}
    is determined by the row lattice $L:=\Lambda_0D$ and
    the rows of $C$ modulo $L$. Indeed, left multiplication
    by an integral parabolic matrix changes $D$ by an
    integral change of basis and changes $C$ by a matrix
    whose rows belong to $L$. From Lemma \ref{gritsenko decomposition}, consider representatives
    \begin{equation*}
        h=
        \begin{pmatrix}1_n&0\\0&H\end{pmatrix},
        \qquad
        v=
        \begin{pmatrix}1_n&N\\0&D\end{pmatrix},
        \qquad
        w=
        \begin{pmatrix}1_n&0\\0&K\end{pmatrix}
    \end{equation*}
    occurring in $H_{t-j}$, $V_t$, and $W_j$, respectively.
    Their product is
    \begin{equation*}
        hvw=
        \begin{pmatrix}
            1_n&NK\\
            0&HDK
        \end{pmatrix}.
    \end{equation*}
    Thus, and because of the form of the representatives occurring from Lemma \ref{gritsenko decomposition}, the left-coset expansion is parametrised by
    lattices
    \begin{equation*}
        M=\Lambda_0K,\qquad
        J=\Lambda_0DK,\qquad
        L=\Lambda_0HDK,
    \end{equation*}
    satisfying
    \begin{equation}\label{inclusion 1}
        \Lambda_0\subseteq M\subseteq p^{-1}\Lambda_0,
        \qquad [M:\Lambda_0]=p^j,
    \end{equation}
    \begin{equation}\label{inclusion 2}
        pM\subseteq J\subseteq M,
        \qquad [M:J]=p^t,
    \end{equation}
    \begin{equation}\label{inclusion 3}
        J\subseteq L,
        \qquad [L:J]=p^{t-j}.
    \end{equation}
    For example, $K$ runs over a set of representatives for $\textup{GL}_n(\mb Z_p) \backslash \textup{GL}_n(\mb Z_p) \m{1_{n-j} \\ & p^{-1}1_j}\textup{GL}_n(\mb Z_p)$, and similarly the others can be read off Lemma \ref{gritsenko decomposition}. Note also from the same Lemma that $N$ runs through $M_n(\mb Z_p)/M_n(\mb Z_p)D$. This means that each row of \(N\) is chosen independently modulo the row lattice \(\Lambda_0D\). Multiplication on the right by \(K\) gives a bijection
\[
\Lambda_0/\Lambda_0D
\xrightarrow{\ \sim\ }
\Lambda_0K/\Lambda_0DK=M/J,
\qquad
z+\Lambda_0D\longmapsto zK+J.
\]
Therefore, each row of the upper-right matrix \(NK\) runs independently through \(M/J\).\\

We now fix one coset of the form
\[
\Gamma
\begin{pmatrix}1_n&C\\0&D_L\end{pmatrix},
\qquad L=\Lambda_0D_L,
\]
which arises in the final left-coset decomposition of $B_{t}$. Here, \(D_L\) is a basis matrix for the lattice \(L\), and the rows of \(C\) are considered modulo \(L\).\\

We first count the possible upper-right matrices $C$. For a fixed $M$, from \eqref{inclusion 2} and \eqref{inclusion 3}, there are
    no contributions unless $pM\subseteq L$. Assume
    this inclusion holds and put
    \begin{equation}\label{k dimension}
        k:=\dim_{\mathbb{F}_p}\bigl(M/(M\cap L)\bigr) = \dim_{\mathbb{F}_p}\bigl((M+L)/L\bigr).
    \end{equation}
    The possible intermediate lattices $J$ satisfy (again from \eqref{inclusion 2}, \eqref{inclusion 3})
    \begin{equation*}
        pM\subseteq J\subseteq M\cap L,
        \qquad [M:J]=p^t.
    \end{equation*}
    Hence their number is
    \begin{equation*}
        {n-k\brack n-t}_p
        =
        {n-k\brack t-k}_p.
    \end{equation*}

    Moreover, the map
    \begin{equation*}
        M/J\longrightarrow (M+L)/L
    \end{equation*}
    has kernel $(M\cap L)/J$, of cardinality $p^{t-k}$.
    Therefore, if the rows of $C$ belong to $M+L$,
    each final upper-right matrix $C$ has
    $p^{n(t-k)}$ lifts back to $M/J$. From \eqref{B three factor expression}, the coefficient that
    this $M$ contributes to the single coset is therefore 
    \begin{equation}\label{weight of temporary lattice}
        c_tc_j\,w_t(k),
        \qquad
        w_t(k):=
        p^{n(t-k)}{n-k\brack t-k}_p,
    \end{equation}
    provided the rows of $C$ belong to $M+L$, and is zero
    otherwise. We put $w_t(k)=0$ if $k>t$.\\

    We now claim that if $L \neq \Lambda_0$, then the corresponding coefficient of the coset is zero. First of all, we may compute 
    \begin{equation*}
        \textup{covol}(L) = \textup{covol}(\Lambda_0).
    \end{equation*}
    Since we assume $L \neq \Lambda_0$, $L$ is not contained in $\Lambda_0$. Hence, we may choose $x \in L\backslash \Lambda_0$ such that $px \in \Lambda_0$. Such an element exists since we may take any $x \in L \backslash \Lambda_0$ and multiply it by powers of $p$ until its next multiple lies in $\Lambda_0$. We keep this $x$ fixed. For every $M$ satisfying \eqref{inclusion 1}, define
    \begin{equation*}
        \widetilde{M} := M + \mathbb{Z}_px.
    \end{equation*}
    This lattice then also satisfies \eqref{inclusion 1}, and additionally, since $x \in L$, we have $M + L = \widetilde{M}+L$. We now consider the collection $\mc{C}(\widetilde{M})$ of all lattices $M$ satisfying \eqref{inclusion 1} that produce the same $\widetilde{M}$ via the procedure above (recall $x$ is fixed). Put now
\[
V=\widetilde M/\Lambda_0,\qquad \dim V=m.
\]
This is an \(\mathbb F_p\)-vector space. The element \(x\) determines a line
\[
\ell=(\Lambda_0+\mathbb Z_px)/\Lambda_0\subseteq V.
\]
and the lattices $M \in \mc{C}(\widetilde{M})$ are in one-to-one correspondence with $\mathbb{F}_p$-subspaces 
\[
U=M/\Lambda_0\subseteq V.
\]
The condition that $M$ gives $\widetilde{M}$ is $M+\mathbb Z_px=\widetilde M$, or equivalently, $U+\ell=V$. If $U$ contains $\ell$, then we must have $U=V$, and there is only one such choice of $M$. If $U$ does not contain $\ell$, then $U \cap \ell= \emptyset$, so this forces $\dim_{\mathbb F_p}(U) = m-1$. Since $M+ L = \widetilde{M}+L$ for all $M \in \mc{C}(\widetilde{M})$, all these lattices have the same \(k\), the same factor \(w_t(k)\), and the same condition on \(C\). Their total contribution is
\[
c_tw_t(k)
\bigl(c_m+p^{m-1}c_{m-1}\bigr)=0,
\]
since $c_m = (-1)^mp^{\langle m-1 \rangle}$ (see \eqref{initial defns}). It now suffices to show that $m \leq t$, so that the lattices $M$ occur in \eqref{B three factor expression}. But $L$ and $\Lambda_0$ have equal covolumes, so from \eqref{k dimension}
\begin{equation*}
    p^k = [\widetilde M:\widetilde M\cap L] =p^m[L:\widetilde M\cap L].
\end{equation*}
Therefore $m \leq k$. Also, the map $M/J \longrightarrow M/(M\cap L)$ is surjective, so $k \leq t$. This shows $m \leq t$.\\

It remains to consider $L=\Lambda_0$. In this case,
    $M\cap L=\Lambda_0$ and $k=j$. Thus the factor
    in \eqref{weight of temporary lattice} becomes
    \begin{equation}\label{new factor}
        c_tc_jp^{n(t-j)}
        {n-j\brack t-j}_p.
    \end{equation}
    From \eqref{inclusion 1}, summing over $M$ is equivalent to summing over $j$-dimensional subspaces $U_M := pM/p \Lambda_0 \subseteq \mathbb{F}_p^{n}$. Now, for a fixed \(M\), the rows of the upper-right matrix \(C\), taken modulo \(\Lambda_0\), belong to \(M/\Lambda_0\). Now $pC \in M_n(\mathbb{Z}_p)$. We define $Z  := pC \pmod p \in M_n(\mathbb{F}_p)$ and we choose a lift $\widetilde{Z} \in M_n(\mathbb{Z}_p)$ of $Z$. The condition on a row \(c\) of \(C\) becomes
\[
c+\Lambda_0\in M/\Lambda_0
\quad\Longleftrightarrow\quad
pc\bmod p\Lambda_0\in U_M.
\]
Consequently, every row of \(Z\) must belong to \(U_M\). Hence, up to the factor \eqref{new factor}, we obtain that the remaining sum is
\[
\sum_{\substack{U\leq\mathbb F_p^n\\\dim U=j}}
\ \sum_{\substack{Z\in M_n(\mathbb F_p)\\
\textup{rows of } Z\subseteq U}}
\Gamma \mbf{n}\left(\widetilde Z/p\right).
\]
Now fix one matrix \(Z\) of rank \(a\). The number of times its coset occurs in the above expression is equal to the number of $j$-dimensional subspaces $U$ which contain $\mc{R} = \{\textup{rows of }Z\}$, with $\dim \mc{R} = a$. These subspaces correspond bijectively to the \((j-a)\)-dimensional subspaces
\[
U/\mc{R}\subseteq \mathbb F_p^n/\mc{R}.
\]
Since the quotient has dimension \(n-a\), their number is
\[
{n-a\brack j-a}_p.
\]
If \(a>j\), there are none. Since the dimension of the column space is equal to the row space of $Z$, \eqref{B finite unipotent expansion} follows.\\

We will now consider the action of each $\mathcal{U}_j$ on $\epsilon_{\mathbb{A}}$. Recall,
\begin{equation*}
        \epsilon_{\mathbb{A}}((q, s); (\mc{L},\tau)) := \beta(q)\sum_{u \in \textup{GL}_{n}(K)} \lambda'(\xi^{-1}uq)\Phi(u, q, s, \tau).
\end{equation*}
The Hermitian matrix corresponding to $\mbf{n}(p^{-1}\widetilde{Z})$ has the form $\m{1_n & S \\ 0_n & 1_n}$,  where $S := (p^{-1}\widetilde{Z}, p^{-1}\widetilde{Z}^{t})$. Then $(q,s)\mbf{n}(p^{-1}\widetilde{Z})^{-1}= (q, s - q S q^{*})$. Therefore, acting with this operator via the rule \eqref{hecke action} produces an additional factor 
\begin{equation}\label{exponential factors}
    e_{p}(-p^{-1}\textup{tr}(T_u\widetilde{Z})),
\end{equation}
where for $u \in \textup{GL}(K)$, we set $T_u :=(q_p^{*}u_p^{*}\tau_pu_pq_p)_{\mathfrak p}$, with the subscript $_{\mf{p}}$ denoting its first componenent. We now show \(T_u\) is integral and has reduced rank \(r\). Recall $\mc{L} = \xi \mc{L}_0$ is a $\mathbb{Z}$-integral lattice. Let us write 
\[
A=(\xi^{*}\tau\xi)_{\mathfrak p}\in M_n(\mathbb Z_p),
\qquad
r=\operatorname{rank}_{\mathbb F_p}(\overline A),
\]
following the assumptions of the proposition. For a nonzero summand, we have $h_u:= \xi^{-1}uq \in E$. We then observe
\begin{equation*}
    T_{u} = h_{u,2}^{t} Ah_{u,1},
\end{equation*}
where $e_{u,i}$ are the two components of $e$ in $\textup{GL}_n(\mathbb{Z}_p)$. This shows $T_u \in M_n(\mathbb{Z}_p)$ with $\textup{rank}(\overline{T}_u)= r$. This rank is therefore the same for every nonzero summand of the epsilon series. Define now the nontrivial character of \(\mathbb F_p\)
\[
\overline\psi_p(a)=e_p(\widetilde a/p).
\]
For fixed \(U\), the sum of the factors \eqref{exponential factors} is
\begin{equation}\label{character sum}
\sum_{\operatorname{im}Z\subseteq U}
\overline\psi_p\!\left(-\operatorname{tr}(\overline T_uZ)\right).
\end{equation}
The columns of \(Z\) run independently through \(U\), so there are \(p^{nj}\) such matrices. Moreover, the character $\overline\psi_p(a)$ is trivial on all of them exactly when
\[
U\subseteq\ker\overline T_u.
\]
Indeed, writing the columns as \(z_1,\dots,z_n\), we have
\[
\operatorname{tr}(\overline T_uZ)
=
\sum_{i=1}^n e_i^{\,t}\overline T_u z_i.
\]
This vanishes for every choice of columns in \(U\) precisely when \(\overline T_u\) vanishes on \(U\). By additive-character orthogonality, \eqref{character sum} is therefore
\[
\begin{cases}
p^{nj},&U\subseteq\ker\overline T_u,\\
0,&\text{otherwise}.
\end{cases}
\]
Since \(\ker\overline T_u\) has dimension \(n-r\), it contains exactly
${n-r\brack j}_p$
subspaces of dimension \(j\). Thus every nonzero summand of \(\epsilon_{\mathbb A}\) is multiplied by the same scalar
\[
p^{nj}{n-r\brack j}_p.
\]
This shows 
\[
\epsilon_{\mathbb A}(\mathcal L;\tau)\mid\mathcal U_j
=
p^{nj}{n-r\brack j}_p\,
\epsilon_{\mathbb A}(\mathcal L;\tau).
\]
Using now \eqref{B finite unipotent expansion}, as well as the identity
\[
\sum_{j=0}^t
(-1)^jp^{\binom j2}
{n-r\brack j}_p
{n-j\brack t-j}_p
=
p^{(n-r)t}{r\brack t}_p,
\]
the Proposition follows.
\end{proof}

\begin{lemma}\label{action of U}
    We have
    \begin{equation*}
        \epsilon_{\mathbb{A}}(\mc{L}; \tau)\mid Q(v,w) = \sum_{r_1, r_2 \geq 0} \epsilon_{\mathbb{A}} \mid \mathbf{T}(-(r_1+r_2))t^{+}(\mf{p}^{r_1})t^{+}(\ol{\mf{p}}^{r_2})v^{r_1}w^{r_2} + F((q,s); v,w),
    \end{equation*}
    where $F((q,s); v,w)$ is a polynomial in $v,w$, whose coefficients are finite linear combinations of the functions $\epsilon_{\mathbb{A}}((q,s); (\mc{L}', \alpha\tau))$ with $\mc{L}'$ are $\mathbb{Z}$-integral $\mathcal{O}_K$-lattices, and which satisfy
    \begin{equation*}
        F((q,s); \chi(\mf{p}),\chi(\ol{\mf{p}})) = 0.
    \end{equation*}
\end{lemma}
\begin{proof}
    Given the previous results, the proof of this Lemma is quite similar to that of \cite[Lemma 5.1.11]{andrianov_book}. In the following, we drop the notation $(q,s)$ from $\epsilon_{\mathbb{A}}$, since it is always the same. We have
    \begin{equation*}
        \epsilon_{\mathbb{A}}(\mc{L}, \tau) \mid B(v,w) = B_{p}(v,w; (\mc{L},\tau))\epsilon_{\mathbb{A}}(\mc{L}, \tau),
    \end{equation*}
    which also gives
    \begin{equation*}
        \epsilon_{\mathbb{A}}(\mc{L}, \tau) \mid B(v,w)^{-1} = B_{p}(v,w; (\mc{L},\tau))^{-1}\epsilon_{\mathbb{A}}(\mc{L}, \tau).
    \end{equation*}
    Now, all coefficients of $Y^{+}(v,w)$ (see \eqref{Y polynomials}), except the constant term, are $p$-homogeneous elements of positive $p$-signature in the sense of \eqref{p homoegeneous}. Hence, from Lemma \ref{pos signature}, we obtain
    \begin{equation}\label{step}
        \epsilon_{\mathbb{A}}(\mc{L}, \tau)\mid Q(v,w) = B_{p}(v,w; (\mc{L},\tau))\epsilon_{\mathbb{A}} \mid \mathbf{T} Y_{+}(v,w)^{-1}B(v,w)^{-1}. 
    \end{equation}
    Writing now the decomposition into $p$-homogeneous components $$\mathbf{T} = \sum_{\sigma \in \mathbb{Z}} \mathbf{T}(\sigma),$$ and substituting $Y(v,w)$, we obtain
    \begin{equation*}
        \epsilon_{\mathbb{A}}(\mc{L}, \tau)\mid Q(v,w) = B_{p}(v,w;(\mc{L},\tau))\epsilon_{\mathbb A} \mid \sum_{\substack{\sigma \in \mathbb{Z}\\r_1, r_2 \geq 0}}\mathbf{T}(\sigma) t^{+}(\mf{p}^{r_1}) t^{+}(\ol{\mf{p}}^{r_2})v^{r_1}w^{r_2} \mid B(v,w; (\mc{L},\tau))^{-1}.
    \end{equation*}
    From Proposition \ref{action of hecke on epsilon}, for each $\sigma, r_1, r_2$, we may write
    \begin{equation}\label{expansion}
        \epsilon_{\mathbb A}(\mc{L}, \tau) \mid \mathbf{T}(\sigma)t^{+}(\mf{p}^{r_1}) t^{+}(\ol{\mf{p}}^{r_2}) = \sum_{\lambda} c(\sigma, r_1, r_2, \lambda)\epsilon_{\mathbb{A}}(\mc{L}_{\sigma, r_1, r_2, \lambda}, \alpha\tau),
    \end{equation}
    where the sum is finite and all the lattices $(\mc{L}_{\sigma, r_1, r_2, \lambda}, \alpha \tau)$ are $\mathbb{Z}$-integral $\mathcal{O}_K$-lattices. Here $\mu = \alpha r_{\infty}\omega$, where $\mu \in \mathbb{Q}_p^{\times}$ is the similitude of $T$ in the usual notation (see statement of Proposition \ref{action of hecke on epsilon}). Now, if $ \mathbf{T}(\sigma)t^{+}(\mf{p}^{r_1}) t^{+}(\ol{\mf{p}}^{r_2})$ is not equal to zero, then this is a $p$-homogeneous operator of $p$-signature $\sigma+r_1+r_2$. Hence, from Lemma \ref{pos signature}, we obtain
    \begin{equation*}
        \epsilon_{\mathbb A}(\mc{L}, \tau) \mid \mathbf{T}(\sigma)t^{+}(\mf{p}^{r_1}) t^{+}(\ol{\mf{p}}^{r_2}) = 0 \qquad \textup{if }\sigma+r_1+r_2 > 0. 
    \end{equation*}
    Substituting \eqref{expansion} into \eqref{step}, we obtain
    \begin{multline}\label{step2}
        \epsilon_{\mathbb{A}}(\mc{L}; \tau)\mid Q(v,w)\\
        = B_{p}(v,w; (\mc{L},\tau))\sum_{\substack{\sigma, r_1, r_2, \lambda \\ \sigma+r_1+r_2=0}}B_{p}(v,w; (\mc{L}_{\sigma, r_1, r_2, \lambda},\alpha\tau))^{-1}c(\sigma, r_1, r_2, \lambda)\epsilon_{\mathbb{A}}(\mc{L}_{\sigma, r_1, r_2, \lambda}, \alpha\tau)v^{r_1}w^{r_2}\\
        +  B_{p}(v,w; (\mc{L},\tau))\sum_{\substack{\sigma, r_1, r_2, \lambda \\ \sigma+r_1+r_2<0}}B_{p}(v,w, (\mc{L}_{\sigma, r_1, r_2, \lambda}, \alpha\tau))^{-1}c(\sigma, r_1, r_2, \lambda)\epsilon_{\mathbb{A}}(\mc{L}_{\sigma, r_1, r_2, \lambda}, \alpha \tau)v^{r_1}w^{r_2}.
    \end{multline}

     Let now $X := \mathbf{T}(\sigma)t^{+}(\mf{p}^{r_1}) t^{+}(\ol{\mf{p}}^{r_2})$. Assume $\sigma+r_1+r_2 <0$ and consider a lattice $(\mc{L}', \rho) := (\mc{L}_{\sigma, r_1, r_2, \lambda}, \alpha\tau)$. We claim that this lattice is degenerate. Since $\mc{L} = \xi \mc{L}_0$, let $A := (\xi^{*}\tau \xi)_{\mf{p}} \in M_n(\mathbb{Z}_p)$ denote the first component of its associated Gram matrix. We assume $\mc{L}$ is non-degenerate, hence we have $v_{p}(\det(A))=0$. Pick a double coset in $X$ and let $D$ its lower right block. From Proposition \ref{action of hecke on epsilon}, we obtain that there is some $\gamma \in \textup{U}(\tau)$ such that
    \begin{equation*}
        \gamma \mc{L'} = \xi \widehat{R}\mc{L}_0,
    \end{equation*}
    where $R \in EDE$. This means that the first component of the associated Gram matrix of $(\mc{L}',\rho)$ is given by 
    \begin{equation*}
        B := \alpha R_1^{-1}A(R_2^{-1})^{t}.
    \end{equation*}
    Consider then 
    \begin{equation*}
        v_{p}(\det (\alpha R_1^{-1}A(R_2^{-1})^{t})) = n v_{p}(\alpha)-v_{p}(\det(R_1)) - v_{p}(\det(R_2)) >0,
    \end{equation*}
    as $$0 > \sigma(X) = v_{p}(\det(R_1))-v_{p}(\det (\mu (R_2^{-1})^{t})) = v_{p}(\det(R_1)) + v_{p}(\det(R_2))-n v_{p}(\alpha),$$ since $v_{p}(\alpha) = v_{p}(\mu)$ (recall $\mu = r_{\infty}\alpha \omega$).\\

    Hence, according to Proposition \ref{action of B}, we obtain that $B_p(v,w, (\mc{L},\tau))B_{p}(v,w, (\mc{L}',\alpha\tau))^{-1}$ is a polynomial in $vw$, divisible by $(1- p^{n}vw)$. Recalling now that $\varphi_{\mathbf{h}}(p\mathcal{O}_K)=1$, hence $\psi(p\mathcal{O}_K)=1$, since $\varphi_{\mathbf{h}}$ restricted to $\mathbb{Q}$ coincides with the non-trivial quadratic character attached to $K/\mathbb{Q}$ and $p$ is a split prime, it follows that the second expression in \eqref{step2} is a poloynomial in $v,w$, whose cofficients are finite linear combinations of $\epsilon_{\mathbb{A}}(\mc{L}_{\sigma, r_1, r_2, \lambda}, \alpha\tau)$, and the polynomial is identically equal to zero if $v=\chi(\mf{p}), w = \chi(\ol{\mf{p}})$, since then $vw = \psi(p\mathcal{O}_K)p^{-n} = p^{-n}$.\\

    If now $\sigma+r_1+r_2 =0$, then $X \in \widetilde{\mathcal{H}}_{0,p}$ (see \eqref{even operators}) and hence from Proposition \ref{level divisor det}, we obtain that each of the lattices $\mc{L}_{\sigma, r_1, r_2, \lambda}$ has the same volume, level and divisor as that of $\mc{L}$. Using the same reasoning as above, $\mc{L}_{\sigma, r_1, r_2, \lambda}$ is non-degenerate (its associated Gram matrix has zero $p$-valuation). In that case $B_{p}(v,w, (\mc{L},\tau)) = B_{p}(v,w, (\mc{L}_{\sigma, r_1, r_2, \lambda}, \alpha\tau))$, showing that the first expression of \eqref{step2} is as claimed.
\end{proof}
We are now ready to state the main Theorem of this paper. We first recall the notion of isometry classes as in \eqref{isometric lattices}. Lemma \ref{invariance of epsilon under U} shows that $\epsilon_{\mathbb{A}}$ depends only on the isometry class $[(\mc{L}, \tau)]$ and similarly for $\theta_{\mb A}$. For our fixed positive definite $\mathbb{Z}$-integral $\mc{O}_K$-lattice $(\mc{L},\tau)$ in $V$, we define the set 
\begin{multline}\label{lambda set}
\Lambda(\ell, \delta,d):= \left\{[(\mathcal{L'}, \rho)] \mid (\mc{L}', \rho) \,\, \mathbb{Z}\textup{-integral } \mathcal{O}_K \textup{-lattice in }\, V, \right.\\\left. \rho \in \mathbb{Q}_{>0}\tau,\,\, \ell((\mathcal{L'}, \rho))= \ell,\,\, \delta((\mathcal{L'}, \rho)) = \delta,\,\, \textup{vol}((\mathcal{L'}, \rho)) = d\right\}.
\end{multline}
Let also $\textup{Gen}(\mc{L}, \tau) \subseteq \Lambda(\ell, \delta,d)$ consists of classes which are locally isometric to $(\mc{L}, \tau)$ for every finite prime $p$. Note that $\Lambda(\ell, \delta, d)$ is finite, as the following Lemma shows.
\begin{lemma}\label{lambda finite}
    The set $\Lambda(\ell, \delta, d)$ is finite.
\end{lemma}
\begin{proof}
To each pair $(\mathcal L',\rho)$, associate the integral quadratic lattice
\[
    q_\rho(x):=\rho[x], \qquad x\in\mathcal L'.
\]
Its underlying $\mathbb Z$-module has rank $2n$, and its associated bilinear form is
\[
 B(x,y) = \textup{tr}_{K/\mathbb Q}(\rho(x,y)).
\]
The determinant of its Gram matrix is
\[
\det(B(x,y)) = |d_K|^{n}\vol(\mc{L', \rho})^{2},
\]
where we identify the volume with the positive generator of its fractional ideal. Thus these quadratic lattices have fixed rank and determinant. By the classical finiteness theorem for positive-definite integral quadratic lattices, there are only finitely many possible $\mathbb Z$-isometry classes.
It now remains to check that adding the $\mathcal O_K$-module structure introduces only finitely many possibilities. Write $\mathcal O_K=\mathbb Z[\alpha]$. Fix one of the quadratic lattices above and a $\mathbb Z$-basis, and let $B$ be its bilinear Gram matrix. Multiplication by $\alpha$ is represented by an integral matrix $A$, because $\mc{L}'$ is an $\mathcal O_K$-lattice. Moreover, $q_\rho(\alpha x) = N_{K/\mathbb Q}(\alpha)q_\rho(x)$, so
\[
A^{t}BA=N_{K/\mathbb Q}(\alpha)B.
\]
In particular, the $i$-th column $a_i$ of $A$ satisfies
\[
a_i^{t}Ba_i = N_{K/\mathbb Q}(\alpha)B_{ii}.
\]
Since $B$ is positive definite, there are only finitely many integral vectors satisfying each such equation. Hence there are only finitely many possible matrices $A$. Finally, a $\mathbb Z$-isometry that intertwines these $\mathcal O_K$-actions is $K$-linear and preserves the Hermitian form. Consequently, there are only finitely many Hermitian isometry classes, proving the assertion.
\end{proof}

We now define the set of ``even Hecke operators'' as follows. Let
\begin{equation*}
    \mbf{E}^{+} := \left\{M \in \mbf{G}_{\mathbb{A}}  \mid \alpha \in N_{K/\mb Q}(K^{\times}),\, \,\, M_p \in \mathbf{D}_{p} \,\  \forall p \mid \mf{bc}\right\},
\end{equation*}
where we write $\mu = \alpha r_{\infty}\omega$ as usual. We define Hecke algebras of even Hecke operators as 
\begin{equation*}
    \mc{H}^{+} := \mc{H}(\mbf{D}, \mbf{E}^{+}), \qquad  \mc{H}_p^{+} := \mc{H}^{+} \cap \mc{H}_p.
\end{equation*}
Also, for any 
\begin{equation*}
    X = \sum_{i}a_i \mbf{D}_{0} \m{A_i & B_i \\ 0 & D_i }\mbf{D}_0 \in \mc{H}_0,
\end{equation*}
we define
\begin{equation}\label{eta map}
    \eta(X) := \sum_{\mf{a} \subseteq \mathcal{O}_K}\sum_{i : N(\mathfrak{a}) \det(D_iA_i^{-1})=\mathcal{O}_K} a_i \chi(\mf{a})  \mbf{D}_{0}\m{A_i & B_i \\ 0 & D_i }\mbf{D}_{0} t^+(\mathfrak{a}).
\end{equation}

\begin{theorem} \label{theorem split}
Let $p\mathcal{O}_K = \mf{p}\ol{\mf{p}}$ with $p \nmid \mf{bc}$. Let $(\mc{L}, \tau)$ be a positive definite $\mathbb{Z}$-integral $\mathcal{O}_K$-lattice inside $V$ with invariants $(\mc{\ell}, \delta, d)$. Then for every Hecke operator $T \in \mathcal{H}_p$ we have
\begin{equation}\label{eq 1 theorem}
(\theta_{\mathbb{A}}\mid T)(x;(\mathcal{L}, \tau)) = \sum_{[(\mathcal{L}', \rho)] \in \Lambda(\ell, \delta,d)} \epsilon((\mathcal{L},\tau), (\mathcal{L'}, \rho),\eta(\mathbf{T})) \theta_{\mathbb{A}}(x;(\mathcal{L}', \rho)), \qquad x \in G_{\mb A},
\end{equation}
where $\mathbf{T}:= \iota_p(T) \in \mc{H}_{0,p}$. Moreover, if $T \in \mc{H}^{+}_{p}$ is an even operator, we have
\begin{equation}\label{eq 2 theorem}
(\theta_{\mathbb{A}}\mid T)(x;(\mathcal{L},\tau)) = \sum_{[(\mathcal{L}',\rho)] \in \textup{Gen}(\mathcal{L}, \tau)} \epsilon((\mathcal{L},\tau),(\mathcal{L'},\rho),\eta(\mathbf{T})) \theta_{\mathbb{A}}(x;(\mathcal{L}', \rho)).
\end{equation}

\end{theorem}

\begin{proof}
    By linearity, we may assume $T$ is supported on a single double coset. We first observe
    \begin{equation*}
        \eta(\mbf{T}) \in \begin{cases}
            \widetilde{\mc{H}}_{0,p} \textup{ if } T \in \mc{H}_p,\\
            \widetilde{\mc{H}}_{0,p}^{+} \textup{ if } T \in \mc{H}^{+}_{p}
        \end{cases}.
    \end{equation*}
    Hence, by Proposition \ref{level divisor det}, the coefficient $\epsilon((\mc{L},\tau), (\mc{L'},\rho), \eta(\mbf{T}))$ equals zero unless $[(\mc{L'}, \rho)] \in \Lambda(\ell, \delta, d)$ or $[(\mc{L}', \rho)] \in \textup{Gen}(L, \tau)$ and $\mbf{T} \in \widetilde{H}_{0,p}^{+}$. In particular, \eqref{eq 2 theorem} is a consequence of \eqref{eq 1 theorem}, and we may prove the two formulae in the form
    \begin{equation}\label{uniform equation}
        (\theta_{\mathbb{A}}\mid T)(x;(\mathcal{L},\tau)) = \sum_{[(\mathcal{L}',\rho)]} \epsilon((\mathcal{L},\tau),(\mathcal{L'},\rho),\eta(\mathbf{T})) \theta_{\mathbb{A}}(x;(\mathcal{L}', \rho)),
    \end{equation}
    where the sum runs over classes of $\mathbb{Z}$-integral $\mathcal{O}_K$-lattices.\\
    
    We return to \eqref{start of epsilon action}. Using \eqref{next expression} and Lemma \ref{action of U}, we obtain
    \begin{multline*}
         \epsilon_{\mathbb{A}}\mid \sum_{r,s=0}^{\infty}t^{-}(\mathfrak{p}^r)t^{-}(\ol{\mf{p}}^{s})v^{r}w^{s}\mathbf{T} = \epsilon_{\mathbb{A}} \mid Q(v,w)\mid \sum_{r,s=0}^{\infty}t^{-}(\mathfrak{p}^r)t^{-}(\ol{\mf{p}}^{s})v^{r}w^{s} \\= \left(\sum_{r_1, r_2 \geq 0} \epsilon_{\mathbb{A}} \mid \mathbf{T}(-(r_1+r_2))t^{+}(\mf{p}^{r_1})t^{+}(\ol{\mf{p}}^{r_2})v^{r_1}w^{r_2}\right)\mid\sum_{r,s=0}^{\infty}t^{-}(\mathfrak{p}^r)t^{-}(\ol{\mf{p}}^{s})v^{r}w^{s} \\+ F(v,w)\mid \sum_{r,s=0}^{\infty}t^{-}(\mathfrak{p}^r)t^{-}(\ol{\mf{p}}^{s})v^{r}w^{s}.
    \end{multline*}
    The formal series
    \begin{equation*}
        F(v,w)\mid \sum_{r,s=0}^{\infty}t^{-}(\mathfrak{p}^r)t^{-}(\ol{\mf{p}}^{s})v^{r}w^{s},
    \end{equation*}
    is defined at $(v,w) = (\chi(\mf{p}), \chi(\ol{\mf{p}}))$ and vanishes there identically. Hence,
    \begin{align*}
         &\epsilon_{\mathbb{A}}\mid \sum_{r,s=0}^{\infty}t^{-}(\mathfrak{p}^r)t^{-}(\ol{\mf{p}}^{s})v^{r}w^{s}\mathbf{T} \mid_{(v,w)=(\chi(\mf{p}), \chi(\ol{\mf{p}})))} \\
         &=\sum_{r,s=0}^{\infty}\epsilon_{\mathbb{A}} \mid  \left(\sum_{r_1, r_2 \geq 0} \mathbf{T}(-(r_1+r_2))t^{+}(\mf{p}^{r_1})t^{+}(\ol{\mf{p}}^{r_2})\chi(\mf{p}^{r_1})\chi(\ol{\mf{p}}^{r_2})\right) \mid t^{-}(\mathfrak{p}^r)t^{-}(\ol{\mf{p}}^{s})\chi(\mf{p}^{r})\chi(\ol{\mf{p}}^{s})\\
         &= \!\begin{multlined}[t][10.5cm]\sum_{r,s=0}^{\infty}\left(\sum_{[(\mc{L'}, \rho)]}\epsilon\left((\mc{L},\tau), (\mc{L}',\rho) \sum_{r_1, r_2 \geq 0} \mathbf{T}(-(r_1+r_2))t^{+}(\mf{p}^{r_1})t^{+}(\ol{\mf{p}}^{r_2})\chi(\mf{p}^{r_1})\chi(\ol{\mf{p}}^{r_2})\right)\epsilon_{\mb A} (\mc{L'}, \rho)\right) \mid \\ t^{-}(\mathfrak{p}^r)t^{-}(\ol{\mf{p}}^{s})\chi(\mf{p}^{r})\chi(\ol{\mf{p}}^{s})\end{multlined},
    \end{align*}
    using Proposition \ref{action of hecke on epsilon}.\\
    
    According to the definition of the $p$-homogeneous component, $\mathbf{T}(-(r_1+r_2))$ is the sum of those summands $a_i \m{A_i & B_i \\ 0& D_i}$ of $\mathbf{T}$ which satisfy $\det (D_i)\det(A_i)^{-1} = p^{-(r_1+r_2)}$. Hence,
    \begin{equation*}
        \sum_{r_1, r_2 \geq 0} \mathbf{T}(-(r_1+r_2))t^{+}(\mf{p}^{r_1})t^{+}(\ol{\mf{p}}^{r_2})\chi(\mf{p}^{r_1})\chi(\ol{\mf{p}}^{r_2}) = \eta(\mathbf{T}) \in \widetilde{\mathcal{H}}_{0,p}.
    \end{equation*}
    Returning then to \eqref{action on theta},
    \begin{equation*}
        (\theta^{(n)}_{\mathbb{A}}\mid \mathbf{T})((q,s); (\mc{L}, \tau))  = \sum_{[(\mc{L}', \rho)]} \epsilon((\mc{L},\tau), (\mc{L'},\rho), \eta(\mathbf{T}))\theta_{\mathbb A}^{(n)}((q,s); (\mc{L}', \rho)),
    \end{equation*}
    which gives \eqref{uniform equation}. By Proposition \ref{level divisor det}, we may restrict the sum to classes in $\Lambda(\ell,\delta,d)$. We now define for $x \in \mbf{G}_{\mb{A}}$,  $$G(x) := (\theta_{\mathbb{A}}\mid T)(x; (\mc{L}, \tau))  - \sum_{[(\mc{L}', \rho)]\in \Lambda(\ell, \delta, d)} \epsilon((\mc{L},\tau), (\mc{L'},\rho), \eta(\mathbf{T}))\theta_{\mathbb A}(x; (\mc{L}', \rho)).$$
    From Lemma \ref{congruence subgroup theta} and Proposition \ref{level divisor det}, all theta series occurring in this expression have the same weight $n$, level and character $\varphi^{-n}$. The good-prime Hecke action preserves the space of such modular forms and the sum is finite. Hence $G$ is an automorphic form on $\gu(n,n)$ of weight $n$.\\
    
    Using now the fact that the leading Fourier part commutes with the parabolic Hecke action shows that $G^{(n)}(y)=0$ for all $y \in \mbf{P}_{\mb A}$. We claim that $G=0$. By Lemma \ref{unitary basis}, there is a unitary basis which consists entirely of elements in $\mbf{P}_{\mb{A}}$. By considering the components of $G$ corresponding to each element of the unitary basis in the sense of Definition \ref{automorphic form on GU}, we deduce that each has a leading Fourier part equal to zero and is hence singular. This means that $G\mid_{\un(n,n)_{\mb{A}}}$ is an automorphic form for $\un(n,n)$ of weight $n$, which is also singular. However, \cite[Lemma 20.8]{shimura2000arithmeticity} shows that any nonzero adelic form of weight $n$ has a nonzero Fourier coefficient of rank $n$ in some finite component. Hence, $G=0$.
\end{proof}

\section{Inert-prime Hecke action on theta series} \label{Section:Inert primes}
This Section is completely analogous to the work of Andrianov \cite{andrianov_book}. Actually, in some sense, it is even simpler. This is due to the fact that in the Hermitian case, a Hermitian vector space over a finite field $\mathbb{F}_{p^2}$ is completely determined by its dimension (even in the characteristic two case). We state the local results first and give their proofs in the Appendix.\\

Let $p\mathcal{O}_K=\mathfrak{p}$ be inert and assume $p\nmid\mathfrak{bc}$. We use $p$ as a uniformiser of $K_{\mathfrak{p}}$, write $\mathcal O=\mathcal O_{K_{\mathfrak p}}$, and identify $\mathcal O/p\mathcal O$ with $\mathbb F_{p^2}$. Write $\Gamma:=\mathbf{G}_p\cap\textup{GL}_{2n}(\mathcal O)$, $\Gamma_P:=\Gamma\cap\mathbf{P}_p$, and $\langle g\rangle_P=\Gamma_Pg\Gamma_P$. We will identify $\mathcal H_p$ with its image in $\mathcal H_{0,p}$. Put
\[
 D_i:=\textup{diag}(1_{n-i},p1_i),\qquad
 \Pi_i:=\langle\textup{diag}(pD_i^{-1},D_i)\rangle_P,
\]
\[
 \Pi_-:=\Pi_0,\quad \Pi_+:=\Pi_n,\quad [p]:=\langle p1_{2n}\rangle_P.
\]
For $H=H^*$ write $n(H)=\m{1_n&H\\0&1_n}$ and set
\[
 Z_a:=\langle n(p^{-1}\textup{diag}(1_a,0_{n-a}))\rangle_P,
 \qquad Z_0:=1,\qquad
 \varphi_j(z):=\prod_{\ell=1}^j(z^\ell-1),\quad \varphi_0(z):=1.
\]

We also note that, as in the split case and similar to \cite[Prop 3.3.28]{andrianov_book}, one can define, using Lemma \ref{triangular decomposition}, a map $\Omega=\Omega_p: \mathcal{H}_{0,p} \longrightarrow \mathcal{H}(\textup{GL}_n(\mathcal{O}_\mathfrak{p}),\textup{GL}_n(K_{\mathfrak{p}}))[x_0^{\pm}]$. We also write $\Phi_n$ for the composition of $\Omega_p$ with the Satake map for $\textup{GL}_n(K_{\mathfrak p})$, extended by fixing $x_0$. Its values lie in $\mathbb Q[x_0^{\pm1},x_1^{\pm1},\ldots,x_n^{\pm1}]$. Set $y_j:=p^{-1}x_j$. We use the standard polynomial $R_p(v)\in\mathcal H_p[v]$ normalised by
\[
 \Phi_n(R_p(v))=\prod_{j=1}^n(1-y_j^{-1}v)(1-y_jv).
\]

For a double coset $X=\langle\m{A&B\\0&D}\rangle_P$, we define its $p$-signature by
\begin{equation}\label{inert signature}
 \sigma_p(X):=v_{\mathfrak p}(\det D)-v_{\mathfrak p}(\det A).
\end{equation}
The signature is well-defined on double cosets and additive in products. As in \eqref{p homoegeneous}, $X(\sigma)$ denotes the $p$-homogeneous component of signature $\sigma$.

% Setting $\mathcal{C}_{\pm}:= \{X \in \mathcal{H}_{0,p} \,\,|\,\, \Pi_{\pm}X = X \Pi_{\pm}\}$ it can be shown (see \cite{heyer2}) that $\Omega$ is injective on $\mathcal{C}_{\pm}$.\newline
\begin{proposition}\label{standard polynomial inert}
We have $R_p(v)=X_-(v)B_p(v)X_+(v)$, where
\begin{align*}
 X_\pm(v)&=\sum_{i=0}^{n}(-1)^ix_{\pm i}v^i, &
 x_{-i}&=(p^{n^2}[p])^{-1}\Pi_-\Pi_{n-i}, &
 x_{+i}&=(p^{n^2}[p])^{-1}\Pi_i\Pi_+,\\
 B_p(v)&=\sum_{i=0}^{n}(-1)^ib_iv^i, &
 b_i&=p^{-i(2n-i)}\sum_{j=0}^{i}\alpha_{ij}(p)Z_{i-j}, &
 \alpha_{ij}(p)&=\frac{\varphi_{n-i+j}(p^2)}
 {\varphi_{n-i}(p^2)\varphi_j(-p)}.
\end{align*}
The coefficients $x_{\pm i}$ have $p$-signature $\pm2i$, and the coefficients of $B_p(v)$ have signature zero.
\end{proposition}
\begin{proof}
See \ref{appendix inert factorisation}.
\end{proof}

\begin{lemma}
    We have the following identities in the ring of formal power series.
    \begin{equation*}
        X_{\pm}(v)^{-1} = \sum_{r=0}^{\infty} p^{-(2n-1)r}t^{\pm}(\mathfrak{p}^{r})v^{r},
    \end{equation*}
    where $t^{\pm}(\mathfrak{p}^r)$ are the (local) Hecke operators defined by equations \eqref{t minus original}
and \eqref{t plus original}.\end{lemma}
\begin{proof}
    See \ref{proof inert inverse factors}.
\end{proof}
We now consider the action of the middle polynomial $B_p(v)$ occurring in Proposition \ref{standard polynomial inert} on the $\epsilon$ series. 
\begin{proposition}
    We have 
    \[
(\epsilon_{\mathbb A} \mid B_{p}(v))((q,s);\mathcal{L},\tau)
=
\prod_{j=0}^{r-1}\bigl(1+(-p)^{-j}v\bigr)
\,\epsilon_{\mathbb A}((q,s);\mathcal{L},\tau),
\]
where $r$ equals the $\mathbb{F}_{p^2}$ rank of the (local) matrix $\overline{(\xi^{*}\tau \xi)_{\mf{p}}}$ (recall $\mathcal{L} = \xi \mathcal{L}_0$).
\end{proposition}
\begin{proof}
See \ref{appendix inert middle}.
\end{proof}
The following Lemma can be proved in a similar way to Lemma \ref{pos signature}.
\begin{lemma}\label{pos signature inert}
    Assume $X$ is $p$-homogeneous. Then if $s_p(X) >0$ and $v_{p}(\vol(\mc{L})) = 0$, we have
    \begin{equation*}
        \epsilon_{\mathbb{A}}(\mathcal{L}, \tau) \mid X  = 0.
    \end{equation*}
\end{lemma}
Using the above Propositions/Lemmas, the following theorem is proven in complete analogy with \cite[Theorem 5.1.10]{andrianov_book}.
\begin{theorem} \label{theorem inert}

Let $p\mc{O}_K = \mf{p}$ is an inert prime with $p \nmid \mf{bc}$. Let $(\mathcal L,\tau)$ be a positive definite $\mathbb Z$-integral $\mathcal O_K$-lattice with invariants $(\ell,\delta,d)$. For every $T\in\mathcal H_p$ and $x\in\mathbf G_{\mathbb A}$, we have
\[
 (\theta_{\mathbb A}\mid T)(x;(\mathcal L,\tau))
 =\sum_{[(\mathcal L',\rho)]\in\Lambda(\ell,\delta,d)}
 \epsilon((\mathcal L,\tau),(\mathcal L',\rho),\eta(\mathbf T))
 \theta_{\mathbb A}(x;(\mathcal L',\rho)),
 \qquad \mathbf T=\iota_p(T).
\]
If $T\in\mathcal H_p^+$, the sum may be taken over $\textup{Gen}(\mathcal L,\tau)$.
\end{theorem}
\begin{proof}
This follows similarly to Theorem \ref{theorem split}/\cite[Theorem 5.1.10]{andrianov_book} using the results of this Section.    
\end{proof}
Finally, to obtain the global version of Theorems \ref{theorem split} and \ref{theorem inert}, we will also need the following Lemma (see Theorem \ref{Theorem m neq n} for the precise statement).
\begin{lemma}\label{eta lemma}
    Let $p_1, \ldots, p_{r}$ denote pairwise distinct prime numbers, not dividing $\mf{bc}$. For $\mathbf{T}_i \in \iota(\mc{H}_{p_i})$ we have
    \begin{equation*}
        \eta(\mathbf{T}_1 \cdots \mathbf{T}_{r}) = \eta(\mathbf{T}_1) \cdots \eta(\mathbf{T}_r).
    \end{equation*}
\end{lemma}
\begin{proof}
Recall 
\begin{equation*}
    \eta(X) = \sum_{\mf{a} \subseteq \mathcal{O}_K}\sum_{i : N(\mathfrak{a}) \det(D_iA_i^{-1})=\mathcal{O}_K} a_i \chi(\mf{a})  \mbf{D}_{0}\m{A_i & B_i \\ 0 & D_i }\mbf{D}_{0} t^+(\mathfrak{a}).
\end{equation*}
Write
\begin{equation*}
    \mathbf{T}_i = \sum_{\sigma} \mathbf{T}_i(\sigma)_{p_i}.
\end{equation*}
Any double coset $M$ occuring in $\mathbf{T}_1(\sigma_1) \ldots T_{r}(\sigma_r)$ has $\det(D_{M}A_M^{-1})\mathcal{O}_K = \prod_{i=1}^{r}p_i^{\sigma_i}\mathcal{O}_K$. Therefore, the condition $N(\mf{a})\det(D_MA_{M}^{-1}) = \mathcal{O}_K$ implies that $\mf{a}$ can only be supported on prime ideals above the $p_i$'s and $\sigma_i = -v_{p_i}(N(\mf{a}))$ for all $i$. We let $\mf{a} := \prod_{i=1}^{r}\mf{a_i}$ with the support of $\mf{a_i}$ being on prime ideals of $K$ above $p_i$. We then have
\begin{equation*}
\eta(\mathbf{T}_1\ldots\mathbf{T}_r) = \sum_{\mf{a}_i} \left(\prod_{i=1}^{r}\chi(\mf{a_i}) \mathbf{T}_i(-\alpha_i)\right) t^{+}(\mf{a_1}\cdots \mf{a}_r).
\end{equation*}
But $t^+(\mf{a_1}\cdots \mf{a}_r) = t^{+}(\mf{a}_1)\cdots t^{+}(\mf{a}_r)$, since the $\mf{a}_i$ are coprime. Also, for $i \neq j$, $$\mathbf{T}_i(-\alpha_i)t^{+}(\mf{a}_j) = t^{+}(\mf{a}_j)\mathbf{T}_i(-\alpha_i).$$ Hence
\begin{equation*}
    \eta(\mathbf{T}_1\ldots\mathbf{T}_r) = \sum_{\mf{a}_i}\prod_{i=1}^{r}\left(\chi(\mf{a}_i)\mathbf{T}_{i}(-\alpha_i)t^{+}(\mf{a}_i)\right) = \prod_{i=1}^{r} \sum_{\mf{a}_i}\left(\chi(\mf{a}_i)\mathbf{T}_{i}(-\alpha_i)t^{+}(\mf{a}_i)\right) = \eta(\mathbf{T}_1) \cdots \eta(\mathbf{T}_r).\qedhere
\end{equation*}
\end{proof}

\section{Hecke Operators and the Siegel Operator}\label{m neq n}
In this Section, we extend the above results to the case $m\neq n$. We first need to discuss the Siegel Operator in our setting and its relation to Hecke operators. To denote the degree of the objects considered, we will often add a superscript/subscript. We define Siegel's $\Phi$ operator 
%$\Phi^{\textup{Sie}}: \mathcal{M}_k(D^{(n)},k,\chi) \longrightarrow \mathcal{M}_k(D^{(n-1)},k,\chi)$ componentwise. By Lemma \ref{unitary basis}, we know that
%\[
%\mathcal{M}_k(D^{(n)},k,\chi)\cong \bigoplus_{b \in \mathcal{B}} M_k(\Gamma^{(n)}_b,k,\chi),
%\]
%where $\Gamma^{(n)}_b := p_b^{-1}D^{(n)}p_b \cap G^{(n)}_\mathbb{Q}$. We note that the size of the set $\mathcal{B}$ is equal to the class number of $K$ and, in particular, is independent of $n$. Moreover, we may select the $b \in \textup{GL}_n(K)_{\mathbf{h}}$ to be of the form $b = \textup{diag}(\tilde{b},1,\ldots,1)$ for an adelic $\tilde{b} \in K_{\mathbf{h}}$. \\
%We then define the map $\Phi^{\textup{Sie}}$ component-wise. That is for an $\mbf{f} \in \mathcal{M}(D^{(n)},k,\chi) $ we write $f_b$ for the $b$-component of $\mbf{f}$. Then we define $\Phi_b^{\textup{Sie}} : M_k(\Gamma^{(n)}_b,k,\chi) \longrightarrow M_k(\Gamma^{(n-1)}_b,k,\chi)$ as in \cite[Equation (27.4)]{shimura2000arithmeticity}. Indeed, by \cite[Lemma 27.2 (4)]{shimura2000arithmeticity} this is well-defined (\textcolor{red}{to check!}). We now note that this definition is equivalent to
\begin{align*}
\Phi^{\textup{Sie}} : \mc{M}(D^{(n)},k,\chi) &\longrightarrow \mc{M}(D^{(n-1)},k,\chi)\\
\Phi^{\textup{Sie}}(\mbf{f})((q,s)) &:= \lim_{\rho \rightarrow \infty} \rho^{-k/2}\mbf{f}\left((q', s')\right),
\end{align*}
where 
\begin{equation*}
    (q_{\infty}',q'_{\mathbf{h}}) = \left(\m{q_{\infty} & 0\\ 0 & \rho^{1/2}}, \, \m{q_{\mathbf{h}}& 0 \\ 0 & 1 }\right),\qquad (s_{\infty}', s_{\mathbf{h}}') = \left(\m{s_{\infty} & 0\\ 0& 0},\, \m{s_{\mathbf{h}}& 0\\ 0& 0}\right).
\end{equation*}
Note that if $F \in  \mc{M}(D^{(n)},k,\chi)$, this definition uniquely determines $\Phi^{\textup{Sie}}(F)$ on the whole $G_{\mb{A}}^{(n-1)}$ because of Lemma \ref{unitary basis}. In particular, it is straightforward to show
\begin{equation}\label{phi on theta}
    \Phi^{\textup{Sie}}( \theta_{\mb{A}, r}) = \theta_{\mb{A}, r-1}, \qquad r \geq 2,
\end{equation}
where by the subscript we denote the degree of the theta series.\\

% In particular if we set $\mbf{g}:= \Phi \mbf{f}$, then we have
% \[
% \mbf{g}(q,s) =  
%  \sum_{\xi \in S^{(n-1)}_{\geq0}} c(\textup{diag}[\xi, 0],\textup{diag}[q,1] ;\mathbf{f}) e(\imath \, tr(q^*_{\infty}\xi q_{\infty} ))e_{\mathbb{A}}(tr(\xi s)),
% \]
% where
% \[
% \mathbf{f}\left((q,s)\right) 
% =  
% \sum_{\xi \in S^{(n)}_{\geq0}} c(\xi, q;\mathbf{f}) e(\imath \, tr(q^*_{\infty}\xi q_{\infty} ))e_{\mathbb{A}}(tr(\xi s)).
% \]

We now define the analogue of the $\Psi$ map of Andrianov \cite[(4.2.64)]{andrianov_book}. Assume first that $p$ denotes a split prime and keep the notation of Section \ref{Hecke rings of the general linear group}. We construct a ring homomorphism 
\[
\Psi_p^{(n)}: \mathcal{H}(\mbf{D}^{(n)}_{0,p},  \mbf{G}^{(n)}_{0,p})\longrightarrow \mathcal{H}(\mbf{D}^{(n-1)}_{0,p},  \mbf{G}^{(n-1)}_{0,p})[u^{\pm},w^{\pm}],
\]
or equivalently
\[
\Psi_p^{(n)} : \mathcal{H}(\Gamma_{2n,n}, G_{2n,n})[x^{\pm}]\longrightarrow \mathcal{H}(\Gamma_{2n-2,n-1}, G_{2n-2, n-1})[x^{\pm}, u^{\pm},w^{\pm}],
\]
as follows: Given an element $$X = \sum_i \alpha_i \Gamma_{2n,n}\m{A_i & B_i \\ 0 & D_i}x^{\delta_i} \in \mathcal{H}(\textup{GL}_{2n,n}(\mathbb{Z}_p), \textup{GL}_{2n,n}(\mathbb{Q}_p))[x^{\pm}],$$ we may select $A_i = \m{A_i' & 0 \\ * & p^{a_i}}$ and $D_i = \m{D_i' & * \\ 0 & p^{d_i}}$. Then we define
\[
\Psi^{(n)}_p(X,u,w) := \sum_i \alpha_i (up^{-n})^{\delta_i-a_i} (wp^{-n})^{d_i} \Gamma_{2(n-1),n-1} \m{A_i' & B_i' \\ 0 & D_i'}x^{\delta_i},
\]
where $B_i'\in M_{n-1}(\mathbb{Q}_p)$ such that $B_i = \m{B_i' & * \\ * & *}$. If $n=1$ then we set $$\Psi_p^{(n)}(X,u,w) := \sum_i \alpha_i (up^{-1})^{a_i}(wp^{-1})^{d_i}.$$

\begin{proposition} With notation as above, we have
\begin{enumerate}
\item For any $X \in \mathcal{H}(\Gamma_{2n,n}, G_{2n,n})[x^{\pm}]$ we have that
\[
\Psi^{(n)}_p(X,u,w) \in \mathcal{H}(\Gamma_{2n-2,n-1}, G_{2n-2,n-1})[x^{\pm}, u^{\pm},w^{\pm}].
\]
\item If $\Phi_{2n,n} :\mathcal{H}(\Gamma_{2n,n}, G_{2n,n})[x^{\pm}] \longrightarrow \mathbb{Q}[x^{\pm}, x_1^{\pm},\ldots,x^{\pm}_{n},y_1^{\pm},\ldots, y^{\pm}_{n}]$ is the Satake map of Definition \ref{spherical map}, we have
\begin{equation}\label{commutativity}
\overline{\Phi}_{2n-2,n-1}\circ \Psi_p^{(n)} = \Xi_n \circ\Phi_{2n,n} ,  
\end{equation}
where $\Xi_n: \mathbb{Q}[x^{\pm},x_1^{\pm},\ldots,x_n^{\pm},y_1^{\pm},\ldots, y^{\pm}_{2n}] \longrightarrow \mathbb{Q}[x^{\pm},x_1^{\pm},\ldots,x_{n-1}^{\pm},u^{\pm}, y_{1}^{\pm},\ldots  y^{\pm}_{2n-1},w^{\pm}]$ is defined as
\[
\Xi_n(x) = up^{-n}x, \qquad \Xi_n(x_{n})= p^{n+1}u^{-1}, \qquad \Xi_n(y_{n}) = w,
\]
\[
\Xi_n(x_i) = px_i, \qquad \Xi_n(y_i) = y_i, \qquad 1 \leq i < n.
\]
and $\overline{\Phi}_{2n-2,n-1}$ is the $\mathbb{Q}[u^{\pm},w^{\pm}]$-linear extension of $\Phi_{2n-2, n-1}$.

\end{enumerate}
    
\end{proposition}
\begin{proof}
This can be proved similarly to \cite[Proposition 4.2.17]{andrianov_book} and follows directly from the definitions of $\Phi_{2n,n}$, $\Xi_n$ and $\Psi^{(n)}_p$.
\end{proof}

We now have the following Proposition.

\begin{proposition}\label{surjectivity psi local split}
    Let $p\nmid \mf{bc}$ be a split prime. For fixed $u_0, w_0 \in \mathbb{C}^{\times}$, the map
    \begin{equation*}
        \Psi_p^{(n)}(\cdot,  u_0, w_0) : \iota_p(\mc{H}_p^{(n)}) \otimes_{\mathbb{Q}}\mathbb{C}\longrightarrow \iota_p(\mc{H}_p^{(n-1)}) \otimes_{\mathbb{Q}} \mathbb{C},
    \end{equation*}
    is a surjective ring homomorphism.
\end{proposition}
\begin{proof}
    First of all, in a similar way to \cite[Proposition 4.2.17 (2)]{andrianov_book}, we have $\Psi_p^{(n)}(\iota_p(\mc{H}_p^{(n)})) \subseteq \Psi_p^{(n-1)}(\iota_p(\mc{H}_p^{(n-1)}))$. From now on, we identify the global Hecke algebras with their images in the parabolic. Since the restriction of $\Phi_{2n,n}$ to the global Hecke algebra is injective, using the commutativity of \eqref{commutativity}, it suffices to show
    \begin{equation}\label{suffices to prove}
        {\Phi}_{2n-2,n-1}(\mc{H}_p^{(n-1)}) = \Xi_n(\Phi_{2n,n}(\mc{H}_p^{(n)})).
    \end{equation}
    Recall the definition of the standard Hecke polynomial \ref{standard polynomials gl} 
    \begin{equation*}
        R_{1,p}^{(2n)}(v) = \sum_{i=0}^{2n} (-1)^{i}p^{\langle i-1\rangle } \pi_{i,2n} v^{i}.
    \end{equation*}
    Applying $\Phi_{2n,n}$, we obtain (cf. \cite[Example 2]{parabolic_gl_II})
    \begin{equation*}
        \Phi_{2n,n}(R_{1,p}^{(2n)}(v)) = \prod_{i=1}^{n}(1-p^{-1}x_i v)(1-p^{n-1}y_iv).
    \end{equation*}
    Applying $\Xi_n$, we obtain
    \begin{multline}\label{relation of standard different genus}
        \Xi_n(\Phi_{2n,n}(R_{1,p}^{(2n)}(v))) = (1-p^nu_0^{-1}v)(1-p^{n-1}w_0v)\prod_{i=1}^{n-1} (1-x_i v)(1-p^{n-1}y_i v) \\= (1-p^nu_0^{-1}v)(1-p^{n-1}w_0v)\Phi_{2n-2,n-1}(R_{1,p}^{(2n-2)}(pv)).
    \end{multline}

That is,
\[
\Phi_{2n-2,n-1}(R_{1,p}^{(2n-2)}(pv))=\left((1-p^nu_0^{-1}v)(1-p^{n-1}w_0v)\right)^{-1}\Xi_n(\Phi_{2n,n}(R_{1,p}^{(2n)}(v)))
\]
But we have that $\left((1-p^nu_0^{-1}v)(1-p^{n-1}w_0v)\right)^{-1} \in \mathbb{C}[[v]]$, the formal power series in $v$, and hence by equating powers of $v$ we deduce that each coefficient of $\Phi_{2n-2,n-1}(R_{1,p}^{(2n-2)}(pv))$ belongs to $\Xi_n(\Phi_{2n,n}(\mc{H}_p^{(n)})) \otimes_{\mathbb{Q}} \mathbb{C}$. This means that each $\Phi_{2n-2,n-1}(\pi_{i,2n-2})$ belongs in this image for all $0 \leq i \leq 2n-2$. Additionally, 
    \begin{equation*}
        \Xi_{n}(\Phi_{2n,n}(\pi_{2n,2n})) = u_0^{-1}w_0\Phi_{2n-2,n-1}(\pi_{2n-2,2n-2}),
    \end{equation*}
    and since $\pi_{2n,2n}$ is invertible, we also obtain $\Phi_{2n-2,n-1}(\pi_{2n-2,2n-2}^{-1}) \in \Xi_n(\Phi_{2n,n}(\mc{H}_p^{(n)})) \otimes_{\mathbb{Q}}\mathbb{C}$. Finally, since $\Xi_n(x) = u_0p^{-n}x$, the same is true for $x^{\pm}$. Hence, all generators of $\mc{H}_p^{(n-1)}$ belong in $\Xi_n(\Phi_{2n,n}(\mc{H}_p^{(n)})) \otimes_{\mathbb{Q}} \mathbb{C}$, which proves \eqref{suffices to prove}, and hence the result.
\end{proof}

We now give the analogous results when $p$ is inert. These are very close to the case treated by Andrianov in \cite{andrianov_book}, so we give only the statements without proofs, which are mostly identical to those of the symplectic case. We start with the definition of $\Psi^{(n)}_p$. 

Given an element $X = \sum_i \alpha_i \mathbf{D}^{(n)}_{0,p} \m{p^{\delta_i}\widehat{D_i}   & B_i \\ 0 & D_i} \in \mathcal{H}_{0,p}$, we may select the elements $D_i$ to be of the form 
\[
D_i = \m{D_i' & * \\ 0 & p^{d_i}}.
\]
We then define
\[
\Psi_p^{(n)}(X,u) :=\sum_i \alpha_i u^{-\delta_i} (up^{-2n})^{d_i} \mathbf{D}^{(n-1)}_{0,p} \m{p^{\delta_i} \widehat{D_i'} & B_i' \\ 0 & D_i'},
\]
where $B_i'$ is such that $B_i = \m{B_i' & * \\ * & *}$. We write $\Phi_n :\mathcal{H}^{(n)}_{0,p} \rightarrow \mathbb{Q}[x_{0}^{\pm},x^{\pm}_{1},\ldots,x_n^{\pm}]$ for the map obtained as the composition of the map $\Omega_p$ of Section \ref{Section:Inert primes} and the Satake map on the Hecke algebra $\mathcal{H}(\textup{GL}_n(\mathcal{O}_\mathfrak{p}),\textup{GL}_n(K_{\mathfrak{p}}))[x_0^{\pm}]$. We then have the following Proposition, which is the analogue of \cite[Proposition 4.2.17]{andrianov_book}. We omit the proof since it follows exactly as in \cite{andrianov_book}. 

\begin{proposition}\label{Psi map inert}
Let $p\nmid \mf{bc}$ be an inert prime and write $\mathcal{H}^{(n)}_{0,p}=\mathcal{H}(\mathbf{D}^{(n)}_{0,p},\mathbf{P}^{(n)}_p)$. Then, with notation as above, we have
\begin{enumerate}
\item For any $X \in \mathcal{H}^{(n)}_{0,p} $ we have that
\[
\Psi^{(n)}_p(X,u) \in \mathcal{H}^{(n-1)}_{0,p}[u^{\pm}].
\]
\item We have
\begin{equation}\label{commutativity}
\overline{\Phi}_{n-1}\circ \Psi_p^{(n)} = \Xi_n \circ\Phi_{n} ,  
\end{equation}
where $\Xi_n: \mathbb{Q}[x^{\pm},x_1^{\pm},\ldots,x_n^{\pm}] \longrightarrow \mathbb{Q}[x^{\pm},x_1^{\pm},\ldots,x_{n-1}^{\pm},u^{\pm}]$ is defined as
\[
\Xi_n(x) = u^{-1}x, \qquad \Xi_n(x_{n})= u, \qquad \Xi_n(x_i) = x_i, \qquad 1\leq i < n.
\]
and $\overline{\Phi}_{n-1}$ is the $\mathbb{Q}[u^{\pm}]$-linear extension of $\Phi_{n-1}$. 

\item
 For fixed $u_0 \in \mathbb{C}^{\times}$, the map
    \begin{equation*}
        \Psi_p^{(n)}(\cdot,  u_0) : \iota_p(\mc{H}_p^{(n)}) \otimes_{\mathbb{Q}}\mathbb{C}\longrightarrow \iota_p(\mc{H}_p^{(n-1)}) \otimes_{\mathbb{Q}} \mathbb{C},
    \end{equation*}
    is a surjective ring homomorphism, except in the case where $u_0 =-p$. In this case, the image is the even subring $\iota_p(\mc{H}_p^{(n-1)+}) \otimes_{\mathbb{Q}} \mathbb{C}$ of $\iota_p(\mc{H}_p^{(n-1)}) \otimes_{\mathbb{Q}} \mathbb{C}$.\\
    
    \noindent Additionally, restricted to the even rings $\iota_{p}(\mc{H}_{p}^{(n)+}\otimes \mb{C}) \longrightarrow \iota_{p}(\mc{H}_{p}^{(n-1)+} \otimes \mb{C})$, the map is always surjective.
\end{enumerate}
\end{proposition}
We now have the following Theorem.
\begin{theorem} \label{Psi and Phi}%Let $p\mc{O}_K = \mf{p}\ol{\mf{p}}$ with $\mf{p} \neq \mf{\ol{p}}$. 
For all $X \in \mc{H}_{0,p}$, we have
\begin{equation*}
        \Phi^{\textup{Sie}}(\theta_{\mb A,n} \mid X) = \begin{cases} \Phi^{\textup{Sie}}(\theta_{\mb A,n}) \mid \Psi_p^{(n)}(X; p^{n-m/2}\vphi_{\mf{p}}(p\mc{O}_K)^{m}, p^{n-m/2}\vphi_{\ol{\mf{p}}}(p\mc{O}_K)^{m}), \,\,\text{if $p$ splits in $K$} \\ 
        \Phi^{\textup{Sie}}(\theta_{\mb A, n}) \mid \left(\sum_{\delta}p^{(2n-m)\delta} \Psi_p^{(n)}(X_{\delta};p^{2n-m}{\vphi_{\mathfrak{p}}(p\mc{O}_K)^{m}} )\right), \,\, \text{if $p$ is inert in $K$}, \end{cases}
    \end{equation*}
    where $X= \sum_{\delta}X_{\delta}$ is the decomposition of $X$ into parts with similitude valuation $\delta$.
    \end{theorem}
\begin{proof}
We give the proof in the split case, the inert case can be done similarly.
    We take a single left coset in the decomposition of $X$, represented by $M' = \m{\mu \wh{L'} & N' \\ 0 & L'}$. We then have, in analogy to the proof of Proposition \ref{action of hecke on epsilon},
    \begin{multline*}
     (\theta_{\mb A,n} \mid M)((q', s')) = \vphi_{\mf{c}}^{-mn}(\omega)\beta_{n}(q'L'^{*})\\\times \sum_{u' \in M_{m,n}(K)} \lambda_n(\xi^{-1}u'q'L'^{*})\Phi(u',q',s', \alpha\tau)e_{\mathbf{h}}\{\alpha\tau[u'_{\mathbf{h}}](-\mu^{-1}q' L'^* N'q'^* )\}.
    \end{multline*}

We now look at the real part of the exponential, which is equal to $e_{\infty}\{\tau[uq'_{\infty}] + \tau[u]s\}$. 
By writing $u' = (u, u_n) \in V^{n-1} \times V$, we observe that this becomes
\begin{equation*}
    \textup{exp}(-2\pi\alpha  \textup{tr}(\tau[u_{\infty}q_{\infty}]))\exp(-2\pi\alpha \rho \tau[u_n]),
\end{equation*}
which tends to $0$ as $\rho \to \infty$, unless $u_n = 0$ (since $\alpha>0$). Assume now that
\begin{equation*}
    M' \longleftrightarrow \m{\m{A' & B' \\0 & D'}, p^{\delta}},
\end{equation*}
under the usual identification. Split $A' = \m{A & 0 \\ * & p^{a}}$ and $D' = \m{D' & * \\ 0 & p^{d}}$. Since $L' = (D', p^{\delta}(A'^{-1})^{t})$, we may write
\begin{equation*}
    L' = \m{L & x \\ 0 & y}, \qquad N' = \m{N & * \\ * & *},
\end{equation*}
with $y = (p^{d}, p^{\delta-a})$.
We observe that 
\begin{equation*}
    \beta_{n}(q'(L')^{*}) = \rho^{m/2}\beta_{n-1}(qL^{*})|\overline{y}|_p^{m/2}\vphi_p(\ol{y})^{m} := \rho^{m/2}\beta_{n-1}(qL^{*})C(M).
\end{equation*}
Now, since $\mu = p^{\delta}$, we may choose $\alpha = p^{\delta}$, $\omega_p=1$ and $\omega_q = p^{-\delta}$ for $q \neq p$. Note that $\mf{c}$ contains all the ramified primes and $p$ splits, hence $\vphi_p(p)=1$. Hence, in total
\begin{align*}
     \Phi^{\textup{Sie}}(\theta_{\mb A, n} \mid M)((q', s')) &= \!\begin{multlined}[t][10.5cm]\vphi_{\mf{c}}^{-m}(\omega)C(M)\vphi_{\mf{c}}(\omega)^{-m(n-1)}\beta_{n-1}(qL^{*})\\\times \sum_{u \in M_{m,n-1}(K)} \lambda_{n-1}(\xi^{-1}uqL^{*})\Phi(u,q,s, \alpha\tau)e_{\mathbf{h}}\{\alpha\tau[u_{\mathbf{h}}](-\mu^{-1}q L^* Nq^* )\}\end{multlined}\\
     &= \vphi_{\mf{c}}^{-m}(\omega)C(M) (\theta_{\mb{A}, n-1} \mid M)((q, s)),
    \end{align*}
    where $M := \m{\mu \wh{L} & N \\ 0 & L}$. Now, in the split case, $\vphi_{\mf{c}}(\omega)=1$, so we have
    \begin{equation*}
        C(M) = |\ol{y}|_p^{m/2}\vphi_p(\ol{y})^{m} = p^{-m(\delta+d-a)/2} \vphi_{\mf{p}}(p)^{m(\delta-a)}\vphi_{\ol{\mf{p}}}(p)^{md}.
    \end{equation*}
    Hence, using \eqref{phi on theta}, we obtain
    \begin{equation*}
        \Phi^{\textup{Sie}}(\theta_{\mb A,n} \mid M) = \Phi^{\textup{Sie}}(\theta_{\mb A, n}) \mid \Psi(M; p^{n-m/2}\vphi_{\mf{p}}(p)^{m}, p^{n-m/2}\vphi_{\ol{\mf{p}}}(p)^{m}).\qedhere
    \end{equation*}
\end{proof}
We now have the following Definition.
\begin{defn}\label{psi definition global}
    If $p \nmid \mf{bc}$ is inert in $K$, write $X = \sum_{\delta} X_{\delta}$, where $X_{\delta}$ is the part of $X$ with similitude valuation $\delta$. Set 
    \begin{equation*}
        \Psi_{p}^{(n)}(X) = \Psi^{(n)}_{m, \vphi, p}(X) := \sum_{\delta}p^{(2n-m)\delta} \Psi_p^{(n)}(X_{\delta};p^{2n-m}{\vphi_{\mathfrak{p}}(p\mc{O}_K)^{m}}).
    \end{equation*}
    If $p \nmid \mf{bc}$ is split in K, let
    \begin{equation*}
        \Psi_{p}^{(n)}(X) = \Psi^{(n)}_{m, \vphi, p}(X) := \Psi_p^{(n)}(X; p^{n-m/2}\vphi_{\mf{p}}(p\mc{O}_K)^{m}, p^{n-m/2}\vphi_{\ol{\mf{p}}}(p\mc{O}_K)^{m}).
    \end{equation*}
    We then define
    \begin{equation}\label{Psi global}
        \Psi^{(n)}(X) = \Psi^{(n)}_{m, \vphi}(X) := \bigotimes_{p \nmid \mf{bc}} \Psi^{(n)}_{p}(X).
    \end{equation}
\end{defn}
Finally, before we state the generalisation of Theorems \ref{theorem split} and \ref{theorem inert} to the case where $m \geq n$, we need the following Lemma.
\begin{lemma} \label{Surjectivity Psi even}
The map $\Psi^{(n)}: \mc{H}^{(n)}\otimes\mb{C} \longrightarrow \mc{H}^{(n-1)} \otimes \mb{C}$ is a surjective ring homomorphism if $m \neq 2n-1$. In the case $m=2n-1$, its image is generated by the local rings $\iota(\mc{H}_{p}^{(n-1)}\otimes \mb{C})$ for all split primes $p\nmid \mf{bc}$ and $\iota(\mc{H}_{p}^{(n-1)+}\otimes \mb{C})$ for inert primes $p \nmid \mf{bc}$. \\

Moreover, the map $\Psi^{(n)}: \iota(\mathcal{H}^{(n)+}\otimes \mb{C}) \longrightarrow \iota(\mathcal{H}^{(n-1)+} \otimes \mb{C})$, i.e. restricted to the even Hecke rings, is always a surjective ring homomorphism.
\end{lemma}

\begin{proof} 
Since $\Psi^{(n)}$ is the restricted tensor product of the local maps $\Psi^{(n)}_{p}$ for $p \nmid \mf{bc}$, its image is generated by the images of these local maps. The first assertion then follows from Proposition \ref{surjectivity psi local split} and Proposition \ref{Psi map inert} (3). \\

For the next assertion, by linearity, we may consider $T \in \mathcal{H}^{(n-1)+}$ consisting of a single double coset. Let us write $T_p \in \mathcal{H}_p^{(n-1)}$ for its local components. We write $\mu$ for the similitude of $T$, and use the usual notation $\mu = \alpha r_{\infty} \omega$, with $\alpha \in \mathbb{Q}_{>0}$. % By Hasse's principle for quadratic spaces we have that $\alpha \in N_{K/\mathbb{Q}}(K^{\times})$ if and only if $\alpha_p \in N_{K_p/\mathbb{Q}_p}(K^{\times}_p)$ for all finite primes $p$ and $\alpha>0$. Moreover since $\alpha_p = \mu_p \mathbb{Z}_p^{\times}$ and $\mu_p =1$ for $p$ ramified in $K$, we have that $\alpha_p$ is locally a norm if and only if $\mu_p$ is for all $p$ inert or split in $K$.
By definition, we may write 
\begin{equation*}
    \alpha = \prod_{p} p^{\delta_p} \in N_{K/\mathbb{Q}}(K^{\times}),
\end{equation*}
where $\delta_p \in \mathbb{Z}$. Now, $\Psi_{p}^{(n)}$ preserves the similitude valuation for all primes $p \nmid \mf{bc}$, and hence $\Psi^{(n)}$ preserves $\alpha$ and hence maps $\mc{H}^{(n)+}\otimes \mb{C}$ into $\mc{H}^{(n-1)+} \otimes \mb{C}$.\\

If now $p$ is inert, then $\delta_p$ is even, since $\alpha$ is a norm (hence is locally a norm by Hasse's principle). Since $p^{\delta_p} = N_{K/\mathbb{Q}}(p^{\delta_p/2})$, we have $T_p \in \mc{H}_p^{(n-1)+}$ and using local even surjectivity, there is $\widetilde{T}_p \in \mc{H}_p^{(n)+} \otimes \mb{C}$ with $\Psi_p^{(n)}(\widetilde{T}_p) = T_p$. If $p$ is split, then using Proposition \ref{surjectivity psi local split}, we obtain that there is $\widetilde{T}_p \in \mc{H}_p^{(n)} \otimes \mb{C}$ such that $\Psi_p^{(n)}(\widetilde{T}_p) = T_p$.\\

Since $\Psi_p^{(n)}$ preserves the similitude valuation, the component of $\widetilde{T}_p$ of valuation $\delta_p$ still maps to $T_p$ for all inert and split primes $p$. Hence, we may assume that every double coset occurring in \(\widetilde T_p\) has similitude valuation $\delta_p$. For the remaining primes where $T_p$ is the identity, we take the identity lift. Then, every double coset occuring in $\widetilde{T} := \otimes_p \widetilde{T}_p$ has $\widetilde{\alpha} = \alpha$. Hence, $\widetilde{T} \in \mc{H}^{(n)+}\otimes \mb{C}$ with $\Psi^{(n)}(\widetilde{T}) = T$.
% Without loss of generality we may assume that $\alpha \in \mathbb{Z}$ if we write $\alpha = \prod_{p} \alpha_p$ for the factorisation of $\alpha$, then $\alpha_p$ is the global similitude corresponding to $T_p$. For $p$ inert we know by ... that the map $\Psi^{(n)}_p$ is surjective on the even part. That is, there exists an even operator $\widetilde{T}_p \in \mathcal{H}_p^{(n)+}$ such that $\Psi_p(\widetilde{T}_p) = T_p$. Moreover we know that for the operator $T_{\textup{split}} := \prod_p T_p$, where the product is over the split pirmes, the  
\end{proof}

We can now state the generalisation of Theorems \ref{theorem inert} and \ref{theorem split} for any $m,n\geq 1$. Below, we write $\theta_{\mathbb{A},n}(\cdot;(\mathcal{L},\tau))$ for the theta series of degree $n$ attached to the Hermitian lattice $(\mathcal{L},\tau)$ of rank $m$. When the lattice is clear, we simply write $\theta_{\mathbb{A},n}$ to emphasise the degree of the theta series.
\begin{theorem}\label{Theorem m neq n}
Let $(\mc{L}, \tau)$ be a positive definite $\mathbb{Z}$-integral $\mathcal{O}_K$-lattice of rank $m \geq 1$ with invariants $(\mc{\ell}, \delta, d)$. For any $n\geq 1$, let $T \in \mc{H}^{(n)}$ and write $\mbf{T} := \iota^{(n)}(T)$. Define $\widetilde{\mbf{T}} \in \iota^{(m)}(\mc{H}^{(m)})\otimes \mb{C} $ as follows.
\begin{enumerate}
    \item If $m < n$, let $\widetilde{\mbf{T}} := \left(\Psi^{(m+1)}\circ \cdots \circ \Psi^{(n)}\right)(\mbf{T})$.
    \item If $m=n$, let $\widetilde{\mbf{T}} := \mbf{T}$.
    \item If $m > n$, choose any $\widetilde{T} \in \mc{H}^{(m)} \otimes \mb{C}$ (if it exists) whose parabolic image $\widetilde{\mbf{T}} = \iota^{(m)}(\widetilde{T})$ satisfies
        \begin{equation}\label{composite Psi}
            \left(\Psi^{(n+1)} \circ \cdots \circ \Psi^{(m)}\right)(\widetilde{\mbf{T}}) = \mbf{T}.
        \end{equation}
\end{enumerate}
Then, we have
\begin{equation}\label{global m geq n}
(\theta_{\mathbb{A},n}\mid T)(x;(\mathcal{L}, \tau)) = \sum_{[(\mathcal{L}', \rho)] \in \Lambda(\ell, \delta,d)} \epsilon((\mathcal{L},\tau), (\mathcal{L'}, \rho),\eta^{(m)}(\widetilde{\mathbf{T}})) \theta_{\mathbb{A},n}(x;(\mathcal{L}', \rho)), \qquad x \in G_{\mb A},
\end{equation}
where $\eta^{(m)}$ denotes the map $\eta$ of \eqref{eta map} in degree $m$. \\

If $T \in \mc{H}^{(n)+}$ is an even operator, then $\widetilde{\mbf{T}}$ is even for $m \leq n$. For $m > n$, an even lift $\widetilde{T} \in \mc{H}^{(m)+} \otimes \mb{C}$ always exists from Lemma \ref{Surjectivity Psi even}. In the cases when $\widetilde{\mbf{T}}$ is even, the sum \eqref{global m geq n} may be taken over $\textup{Gen}(\mc{L},\tau)$.
\end{theorem}
\begin{proof} 
For the operator $\widetilde{\mathbf{T}}$, we may apply Theorems \ref{theorem split} and \ref{theorem inert}, as well as Lemma \ref{eta lemma} and Lemma \ref{multiplicative epsilon} (3), to obtain the formula \eqref{global m geq n} for $m=n$, exactly as in \cite[p. 306]{andrianov_book}. Assume now $m > n$. Applying Theorem \ref{Psi and Phi} $m-n$ times, we obtain
\begin{equation*}
    (\Phi^{\textup{Sie}})^{m-n} \left(\theta_{\mb{A},m} \mid \widetilde{T}\right) = \theta_{\mb{A},n} \mid \left(\Psi^{(n+1)} \circ \cdots \circ \Psi^{(m)}\right)(\widetilde{\mbf{T}}) = \theta_{\mb{A},n} \mid T,
\end{equation*}
where we have used also the fact that $(\Phi^{\textup{Sie}})^{m-n} (\theta_{\mathbb{A},m}) = \theta_{\mathbb{A},n}$, which follows from \eqref{phi on theta}, we obtain formula \eqref{global m geq n} in degree $n$. For the even case, we apply Lemma \ref{Surjectivity Psi even} successively to find the Hecke operator $\tilde{\mathbf{T}}$.\\

Assume now $m < n$. As in the proof \cite[Theorem 5.1.15]{andrianov_book}, it suffices to show that if 
\begin{equation*}
    F := (\theta_{\mb{A},n}\mid T)(x; (\mc{L},\tau)) - \sum_{[(\mc{L}', \rho)]}\epsilon((\mc{L}, \tau), (\mc{L}', \rho), \eta^{(m)}(\widetilde{\mbf{T}}))\theta_{\mb{A},n}(x; (\mc{L}', \rho)),
\end{equation*}
satisfies $(\Phi^{\textup{Sie}})^{n-m} (F) = 0$, then $F =0$. Here the sum is finite. We consider the Fourier expansion of $F$ as in \cite[Proposition 20.2]{shimura2000arithmeticity}. In particular, we may write
\begin{equation*}
    F((q,s)) = \ol{\det(q_{\infty})}^m\sum_{A \in S_{+}^{(n)}} c(A, q_{\mbf{h}}) e(\textup{tr}(iq_{\infty}^{*}Aq_{\infty}))e_{\mb A}(\textup{tr}(As)),
\end{equation*}
where $S^{(n)}_{+} := \{A^{*}=A \in M_n(K) \mid A \geq 0\}$. The condition $(\Phi^{\textup{Sie}})^{n-m} F =0$, then implies 
\begin{equation}\label{vanishing from phi}
    c\left(\m{A_0 & \\ & 0_{n-m}}, \m{X & \\ & 1_{n-m}}\right) = 0 \qquad \forall A_0 \in S^{(m)}_{+}, \,\,\,\  X\in \gl_m(K)_{\mbf{h}}.
\end{equation}
Let $c(A, q)$ be any Fourier coefficient of $F$. Since $\textup{rk}(A) \leq m$, there is $U \in \textup{GL}_n(K)$ such that $A[U] = \textup{diag}(A_0,0_{n-m})$ with $A_0 \in S^{(m)}_{+}$. We can then assume $A = \textup{diag}(A_0,0_{n-m})$, since $c(A, q) = \ol{\det(U)}^{-m}c\left(A[U], U^{-1}q\right)$ for all $U \in \textup{GL}_n(K)$ from \cite[(20.9c)]{shimura2000arithmeticity}.\\

Using now Proposition \ref{adelic expansion of theta}, if $\mc{L}' = \xi' \mc{L}_0$, with $\xi' \in \gl_m(K)_{\mbf{h}}$ we observe that the Fourier coefficients of each $\theta_{\mb{A}}(\cdot; (\mc{L}',\rho))$ have the form
\begin{equation}\label{theta exp 1}
    c_{\theta}(A, q) := \beta(q)\sum_{\substack{u \in M_{m,n}(K)\\ \rho[u] = A}} \lambda(\xi'^{-1}uq).
\end{equation}
Also, if $M = (\mu, L, N)$ (see Section \eqref{The Hecke action on epsilon series}) is any coset representive in the decomposition of $T$, the Fourier coefficients of $\theta_{\mb{A},n}(\cdot; (\mc{L},\tau))\mid M$ have the form
\begin{equation}\label{theta exp 2}
    c_{\theta \mid M}(A, q) := K_M \beta(q)e_{\mbf{h}}(-\textup{tr}(A\mu^{-1}q L^{*}Nq^{*})) \sum_{\substack{u \in M_{m,n}(K)\\ \alpha\tau[ u] = A}} \lambda(\xi^{-1}uq L^{*}),
\end{equation}
where $K_M$ is independent of $q$. We observe that since $A = \textup{diag}(A_0, 0_{n-m})$, then $\rho[u] = A$ implies $u = (v, 0_{n-m})$. Additionally, it is possible to find $k \in \textup{GL}_n(\widehat{\mc{O}}_K)$ such that 
\begin{equation*}
    qk = \m{x & 0 \\ y &z}.
\end{equation*}
From \cite[(20.9d)]{shimura2000arithmeticity}, we can then assume $q$ has this form. Then $uq = (vx, 0) = u q_0$ and $A[q] = A[q_0]$, where $q_0:= \textup{diag}(x, 1_{n-m})$. These relations give the same lattice conditions in \eqref{theta exp 1} and \eqref{theta exp 2}, as well as the same exponential in \eqref{theta exp 2}. As $\beta(q)$ is common to all the terms defining $F$, we have $c(A,q) = \beta(q)\beta(q_0)^{-1}c(A,q_0)$. Using \eqref{vanishing from phi}, we obtain that $F$ vanishes on the Siegel parabolic. Using then Lemma \ref{unitary basis}, we obtain $F =0$.
\end{proof}
\begin{remark}\label{admissible operators}
    We remark that in the case $m > n$, Lemma \ref{Surjectivity Psi even} implies that an operator $\widetilde{T}$ satisfying \eqref{composite Psi} exists if $m$ is even or $n \leq m < 2n+1$. If $m$ is odd and $m \geq 2n+1$, such $\widetilde{T}$ exists if and only if $T$ is supported on double cosets with even similitude valuation at every inert prime $p \nmid \mf{bc}$. For $m \leq n$, there are no restrictions.
\end{remark}
In particular, if $\ell, \delta \in \mb{Z}_{>0}$ and $d \in \mb{Q}_{>0}$, introducing the spaces
\[
\Theta^n(\ell,\delta,d) := \textup{span}_{\mb{C}}\{ \theta_{\mathbb{A},n}(\cdot;(\mathcal{L}',\rho))\mid[(\mathcal{L}',\rho)] \in \Lambda(\ell, \delta,d) \},
\]
and for a fixed lattice $(\mathcal{L},\tau) \in \Lambda(\ell,\delta,d)$ the subspace
\[
\Theta^n(\textup{Gen}(\mathcal{L},\tau)) := \textup{span}_{\mb{C}}\{\theta_{\mathbb{A},n}(\cdot;(\mathcal{L}',\rho))\mid [(\mathcal{L}',\rho)] \in \textup{Gen}(\mathcal{L},\tau) \},
\]
we have the following Corollary.
\begin{corollary}\label{theorem}
    For any positive definite $\mb{Z}$-integral Hermitian $\mc{O}_K$-lattice $(\mc{L},\tau)$, the space $\Theta^n(\textup{Gen}(\mathcal{L},\tau))$ is $\mathcal{H}^{(n)+}$-stable. Moreover, if $m$ is even or $m < 2n +1$, then the space $\Theta^n(\ell,\delta,d)$ is $\mathcal{H}^{(n)}$-stable. If $m$ is odd and $m \geq 2n+1$, $\Theta^n(\ell,\delta,d)$ is stable under the subalgebra of $\mc{H}^{(n)}$ generated by $\mc{H}_p^{(n)}$ for split primes $p$ in $K$ and $\mc{H}_p^{(n)+}$ for inert primes $p$ in K, all away from the level.
\end{corollary}

% \begin{remark} \label{the case m<n} In Theorem \ref{Theorem m neq n} we are assuming that $m \geq n$. We can also prove a similar statement when $m < n$, but only when the class number of $K$ is one. Indeed, as in the proof of \cite[Theorem 5.1.15]{andrianov_book} in this case we need that the space of singular forms is spanned by theta series. In the Hermitian case, this is known only in the case where the class number of $K$ is equal to one \cite[Theorem 2]{kojima}.  
% \end{remark}
\section{Eichler Matrices and Representation Numbers}\label{Eichler Matrices and Representation Numbers}
The aim of this final Section is to define the analogous Eichler and Theta matrices of \cite[Section 5.2]{andrianov_book} in this setting, and to compute the theta matrix of certain Hecke operators. For the sake of space-saving, we won't specialise the results to the case when one can only consider classes in $\textup{Gen}(\mc{L},\tau)$. This is a straightforward generalisation when we restrict to even Hecke operators. \\

We still work with our fixed rank-$m$ $\mb{Z}$-integral $\mc{O}_K$-Hermitian lattice $(\mc{L}, \tau)$ with characteristics $(\ell, \delta, d)$. We recall that $\Lambda(\ell, \delta, d)$ of \eqref{lambda set} is finite from Lemma \ref{lambda finite}. Let us now consider a set of representatives $\{\mathscr{L}_i := (\mc{L}_i, \tau_i)\},\,\, i=1,\ldots, H$ for this set of classes with $\mathscr{L}_1 = \msc{L} = (\mc{L}, \tau)$. For a fixed $n \geq 1$, we consider the theta series $\mbf{g}_{i}(x) := \theta_{\mathbb{A}}(x, \msc{L}_i)$ attached to $\msc{L}_i$. We can obtain information about the representation numbers of these lattices by considering the various components of $\mbf{g}_{i}$ in the sense of Definition \ref{automorphic form on U}. In particular, let us write $\mc{L}_i = \xi_i \mc{L}_0$, where $\xi_i \in \textup{GL}_m(K)_{\mbf{h}}$ and $\xi_1 = \xi$. For $b \in \textup{GL}_n(K)_{\mbf{h}}$, let $x_{b} := \textup{diag}(b, \widehat{b})$ and $\mc{M}_b := b \mc{O}_K^n$. From \cite[(A5.10)]{shimura2000arithmeticity} we obtain that the $x_b$-component of $\mbf{g}_i$ is
\begin{equation*}
    g_{i}^{(b)}(z) := \beta(b) \sum_{u \in M_{m,n}(K)} \lambda(\xi_{i}^{-1}ub)e(\textup{tr}(\tau_i[u] z)) = \beta(b) \sum_{A \in S} r(\msc{L}_i,\mc{M}_b, A)e(\textup{tr}(Az)),\qquad z \in \mathbb{H}_n,
\end{equation*}
where $S = \{A \in M_n(K) \mid A^{*}=A\}$ and 
\begin{equation}\label{rep numbers}
    r(\msc{L}_i,\mc{M}_b, A) := \#\{u \in \xi_{i}Wb^{-1} \cap M_{m,n}(K) \mid \tau_{i}[u] = A\},
\end{equation}
with $W = \prod_{p \nmid \infty}M_{m,n}(\mathcal{O}_{K,p})$ the set of \eqref{R and E}. Using the fact $\mc{L}_k = \xi_k \mc{O}_K^m$ for all $k$, we obtain equivalently, that 
\begin{equation}\label{rep numbers 2}
    r(\msc{L}_i,\mc{M}_b, A) = \#\{u \in M_{m,n}(K) \mid u\mc{M}_b \subseteq \mc{L}_i \textup{ and } \tau_{i}[u] = A\}.
\end{equation}
We set 
\begin{equation}\label{automorphisms}
    e_i := |\textup{Aut}(\msc{L}_i)| = \{u \in \textup{GL}_m(K) \mid u\mc{L}_i = \mc{L}_i \textup{ and } \tau_i[u] = \tau_i\}, \qquad i=1,\ldots, H.
\end{equation}
\begin{lemma}\label{linear independence}
   If $n \geq m$, the adelic theta series $\mbf{g}_1, \ldots, \mbf{g}_{H}$ are linearly independent.
\end{lemma}
\begin{proof}
    We first consider the case $n=m$. By taking $b = \xi_j$, we have $\mc{M}_b = \mc{L}_j$. In that case, we have
    \begin{equation}\label{rep numbers lattices}
    r(\msc{L}_i, \msc{L}_j; \tau_j) = \begin{cases}
        e_i, & i=j,\\
        0& i\neq j.
    \end{cases}
\end{equation}
This is because $r(\msc{L}_i, \msc{L}_j; \tau_j)$ describes isometric embeddings of $\msc{L}_j$ into $\msc{L}_i$, and the inclusion becomes equality since both have the same volume.\\
    
Hence, if $\sum_{i}c_i \mbf{g}_i=0$, by considering the $\tau_{j}$-Fourier coefficient of the $x_{\xi_j}$-component, we obtain
\begin{equation*}
    0 = \beta(\xi_{j}) \sum_{i=1}^{H} c_{i}r(\msc{L}_i, \msc{L}_j, \tau_j) = \beta(\xi_j)e_{j}c_{j},
\end{equation*}
which implies $c_j=0$ for all $j$. If now $n > m$, applying $(\Phi^{\textup{Sie}})^{n-m}$ to $\sum_{i}c_i \mbf{g}_i=0$ reduces the Lemma to the case $m=n$.
\end{proof}
Using the notation $\mbf{g}_i(x) = \theta_{\mb A,n}(x, \msc{L}_i)$ as above, we now define the theta column of $\msc{L}$ as 
\begin{equation*}
    \mbf{t}(x) = \mbf{t}(x, \msc{L}) := \m{e_1^{-1} \mbf{g}_i(x)\\ \vdots \\e_{H}^{-1}\mbf{g}_H(x)}, \qquad x \in G_{\mathbb{A}}.
\end{equation*}
Its $x_{b}$-component is 
\begin{equation*}
    \mbf{t}_{b}(z) = \m{e_1^{-1}g_{1}^{(b)}(z) \\ \vdots \\ e_{H}^{-1} g_{H}^{(b)}(z)}, \qquad z \in \mathbb{H}_n.
\end{equation*}
If now $T \in \mc{H}^{(n)}$ is as in Remark \ref{admissible operators}, let $\mbf{T} := \iota^{(n)} (T)$ and put $\widetilde{\mbf{T}}$ the associated degree-$m$ operator of Theorem \ref{Theorem m neq n}. We then obtain
\begin{equation*}
    e_{\alpha}^{-1}(\mbf{g}_\alpha \mid T)(x) = \sum_{\beta=1}^{H}\theta_{\alpha\beta}(T)e_{\beta}^{-1} \mbf{g}_{\beta}(x),
\end{equation*}
where $\theta_{\alpha\beta}(T):= e_{\alpha}^{-1}e_{\beta} \times\epsilon(\msc{L}_{\alpha}, \msc{L}_{\beta}, \eta^{(m)}(\widetilde{\mbf{T}}))$. We can then rewrite the above relation as
\begin{equation}\label{theta matrix}
    \mbf{t} \mid T = \Theta(T) \mbf{t},
\end{equation}
where the Hecke operator acts componentwise and $\Theta(T):= (\theta_{\alpha\beta}(T)) \in M_{H}(\mb{C})$. We call $\Theta(T)$ the theta matrix of $T$. If $n \geq m$, it is uniquely determined by the linear independence of Lemma \ref{linear independence}. If $n<m$, it need not be determined uniquely.\\

We now define columns of representation numbers. For brevity, let $r_{ib}(A):= r(\msc{L}_i, \mc{M}_b, A)$ for $A \in S$. For each $b$, we define
\begin{equation*}
    \mbf{r}_{b}(A) := \m{e_1^{-1} r_{1b}(A) \\ \vdots \\ e_{n}^{-1}r_{Hb}(A)} .
\end{equation*}
We also define $(\mbf{r} \mid T)_{b}(A)$ for $T \in \mc{H}^{(n)}$ as follows. For its $i$-th entry, we consider the $x_b$-component of $\mbf{g}_i \mid T$ and take the $A$-Fourier coefficient multplied by $e_{i}^{-1}\beta(b)^{-1}$. Comparing now coefficients in \eqref{theta matrix}, we obtain
\begin{equation}\label{relation r theta matrix 1}
    (\mbf{r}\mid T)_{b}(A) = \Theta(T)\mbf{r}_b(A)
\end{equation}
for all $b$. We then have the following Lemma.
\begin{lemma}\label{relation theta and rep numbers}
    Assume $n \geq m$ and let 
    \begin{equation*}
        b_j := \textup{diag}(\xi_j, 1_{n-m}), \quad A_j := \textup{diag}(\tau_j, 0_{n-m}), \quad j=1,\ldots, H.
    \end{equation*}
    
    For all $T \in \mc{H}^{(n)}$, we have 
    \begin{equation*}
        \Theta(T) = \left((\mbf{r} \mid T)_{b_1}(A_1),  \ldots, (\mbf{r} \mid T)_{b_H}(A_H)\right).
    \end{equation*}
\end{lemma}
\begin{proof}
    We first observe that $r(\msc{L}_i, M_{b_j}; A_j) = e_{i}\delta_{ij}$. Indeed, writing a representation matrix $u = (u_0, u_1)$, with $u_0 \in M_{m}(K)$, we have $\tau_i[u] = A_j$ implies $u_1=0$. Hence $u\mc{M}_{b_j} \subseteq \mc{L}_i$ implies $u_0 \mc{L}_j \subseteq \mc{L}_i$. Hence, the claim follows as in \eqref{rep numbers lattices}. Hence, we obtain
    \begin{equation*}
        (\mbf{r}\mid T)_{b_j}(A_j) = \Theta(T)\mbf{r}_{b_j}(A_j) = \Theta(T)f_j,
    \end{equation*}
    where $f_{j}$ is the standard basis vector in $\mathbb{C}^{H}$. Hence, the $j$-th column of $\Theta$ equals $(\mbf{r}\mid T)_{b_j}(A_j)$, as required.
\end{proof}
Finally, we give the Theta matrices of certain Hecke operators of similitude $p$. We start with the split-prime case.\\

Let $T_{r}(p),\,\, 0 \leq r \leq 2n$ denote the Hecke operator corresponding to $(\pi_{r,2n}, p)$ under the identification $\mc{H}(\mbf{D}_p^{(n)}, \mbf{G}_p^{(n)}) \cong \mc{H}(\textup{GL}_{2n}(\mathbb{Z}_p), \gl_{2n}(\mb Q_p))[x^{\pm}]$ of Section \ref{Split-prime local unitary Hecke rings}. Here, $\pi_{r,2n}$ are the generators of the Hecke algebra of $\gl_{2n}$ as in \eqref{generators of gl}. We can compute their action explicitly, analogously to \cite[Theorem 5.2.5]{andrianov_book}. We first define for $0 \leq k \leq m$
\begin{equation*}
    R_{ij}^{(k)}(p) := \#\left\{u \in M_m(K)\mid  \tau_i[u] = p\tau_j, \,\, u\mc{L}_j \subseteq \mc{L}_i, \,\, \mc{L}_i/(u\mc{L}_j) \cong \left(\mc{O}_K/\mf{p}\right)^{k} \oplus \left(\mc{O}_K/\mf{\ol{p}}\right)^{m-k}  \right\}.
\end{equation*}
We also define for $0 \leq k \leq m$, the \textit{Eichler matrices}
\begin{equation*}
    M_k(p) := (e_{i}^{-1} R^{(k)}_{ij}(p))_{i,j=1}^{H}.
\end{equation*}
\begin{theorem}\label{eichler matrices split}
    For $0 \leq r \leq 2n$ we may take
    \begin{equation*}
        \Theta(T_r(p)) = \chi(\mathfrak p)^{\,n-r} \sum_{k=0}^{m} p^{-k(r-k+m-n)}{2n-m\brack r-k}_p\,M_k(p).
    \end{equation*}
    Here, for \(a\in\mathbb Z\) and \(b\geq0\), we set
\begin{equation}\label{extended binomial}
{a\brack b}_p
:=
\prod_{j=0}^{b-1}
\frac{p^{a-j}-1}{p^{b-j}-1},
\qquad
{a\brack0}_p=1,
\end{equation}
and set \({a\brack b}_p=0\) for \(b<0\).
\end{theorem}
\begin{proof}
    We first deal with the case $m=n$. Let $A_{a},\, 0 \leq a \leq n$ denote the matrices of \eqref{A matrices} and $D_{a,b} := (A_b, p(A_a^{-1})^{t})$. We then have from the proof of Lemma \ref{Embeddings in the parabolic Hecke algebra} that 
    \begin{equation*}
        \iota_p(T_r(p)) = \sum_{\substack{a+b=r\\0 \leq a,b \leq n}} \mbf{D}_{0,p} \m{p \wh{D}_{a,b} & \\ & D_{a,b}}\mbf{D}_{0,p}.
    \end{equation*}
    Let $M_{a,b}$ denote the corresponding matrix. We can then write from Lemma \ref{triangular decomposition}
    \begin{equation*}
        \mbf{D}_{0,p} \m{p \wh{D}_{a,b} & \\ & D_{a,b}}\mbf{D}_{0,p} = \sum_{\substack{L \in E_p \backslash E_p D_{a,b}E_p\\ B \in B(L, M_{a,b}) \pmod L}} \mbf{D}_{0,p} \m{p\widehat{L} & B\\ & L}.
    \end{equation*}
    Analogously to the proof of Proposition \ref{action of hecke on epsilon}, we obtain that the $x_j$-component of $\mbf{g}_i \mid T_{M_{a,b}}$ is given by (recall $x_j = \textup{diag}(\xi_j, \widehat{\xi}_j)$ and $\mu = p$, so $\alpha=p$ and $\omega$ is trivial)
    \begin{equation*}
        (\mbf{g}_{i} \mid T_{M_{a,b}})_{x_j}(z) = \beta(\xi_{j})\sum_{L}\beta(L^{*})\sum_{u \in M_n(K)}\lambda(\xi_i^{-1}u\xi_jL^{*})s((-\tau_i)[u\xi_jL^{*}], L, M_{a,b})e(\textup{tr}(p\tau_i[v]z)),
    \end{equation*}
    for $z \in \mb{H}_n$. Now, from Lemma \ref{relation theta and rep numbers}, to evaluate $\Theta(T_{r}(p))_{ij}$, we need to evaluate the $\tau_j$-Fourier coefficient. So, in the above sum, for fixed $L$, we look at the set
    \begin{equation}\label{counting reps set}
        \mc{R}_{ij}(L) := \{u \in M_n(K) \mid p\tau_i[u] = \tau_j \textup{ and } \xi_i^{-1}u\xi_jL^{*} \in W\}.
    \end{equation}
    Under these conditions, we observe that the $\tau_j$-Fourier coefficient of the $x_{j}$-component of $\mbf{g}_i \mid T_{M_{a,b}}$ equals
    \begin{equation*}
        \beta(\xi_j) \sum_{L} \beta(L^{*})s(-p^{-1}\tau_j[\xi_j L^{*}], L,M_{a,b}) \#\mc{R}_{ij}(L).
    \end{equation*}
    Recall now 
    \begin{equation*}
        s(-p^{-1}\tau_j[\xi_j L^{*}], L,M_{a,b}) = \sum_{N \in B(L, M_{a,b}) \pmod L} e_{p}\{-p^{-1}\tau_j[\xi_j]L^{*}N\}.
    \end{equation*}
    But according to Lemma \ref{gritsenko decomposition}, we may choose the representatives $N$ in the form $N = p\widehat{L}Y$, for some $Y \in \mc{S}_p$. Then
    \begin{equation*}
        e_{p}\{-p^{-1}\tau_j[\xi_j]L^{*}N\} = e_{p}\{\tau_j[\xi_j]Y\} =1,
    \end{equation*}
    since $Y \in \mc{S}_p$ and the lattice $(\mc{L}_j, \tau_j)$ is $\mathbb{Z}$-integral, hence $\tau_{j}[\xi_j] \in \widehat{S}_p$. Hence,
    \begin{equation*}
        s(-p^{-1}\tau_j[\xi_j L^{*}], L,M_{a,b}) = \#\{N \in B(L,M_{a,b}) \pmod{L}\} = p^{(n-a)b},
    \end{equation*}
    using the same reasoning as in the proof of Lemma \ref{Embeddings in the parabolic Hecke algebra}. Moreover, 
    \begin{equation*}
        \beta(L^{*}) = \beta(D_{a,b}^{*}) = |\det(D_{a,b})^{*}|_{K,p}^{n/2}\vphi_p(\det(D_{a,b}^{*}))^n,
    \end{equation*}
    from \eqref{betadef}. But $\det(D_{a,b}^{*}) = \mf{p}^{n-a}\ol{\mf{p}}^b$ and hence $\beta(L^{*}) = \chi(\mf{p})^{n-a}\chi(\ol{\mf{p}})^b$, using \eqref{chi defn}. Hence, in total,
    \begin{equation}\label{theta relation intermediate}
        \Theta(T_{r}(p))_{ij} = e_{i}^{-1}\sum_{\substack{a+b=r\\0\leq a,b \leq n}} p^{(n-a)b}\chi(\mf{p})\chi(\ol{\mf{p}})^b\sum_{L \in E_p\backslash E_pD_{a,b}E_p} \#\mc{R}_{ij}(L).
    \end{equation}
    We now need to count $\#\mc{R}_{ij}(L)$. Let us fix $u \in M_n(K)$ such that $p\tau_i[u]=\tau_{j}$, and count how many $L \in E_p\backslash E_p D_{a,b}E_p$ satisfy $\xi_i^{-1}u\xi_j L^{*} \in W$. Under the split-prime identification, we may write $L = (B,p(A^{-1})^{t})$ with $A \in \Gamma A_a \Gamma$ and $B \in \Gamma A_b \Gamma$, with $\Gamma = \textup{GL}_{2n}(\mb{Z}_p)$. Then, after writing $p\xi_i^{-1}u\xi_j = (C_1, C_2)$, the condition is equivalent to
    \begin{equation}\label{condition rep}
        C_1 A^{-1} \in M_n(\mb{Z}_p), \qquad C_2B^{t} \in pM_n(\mb{Z}_p).
    \end{equation}
    Let now $G_i = (\tau_i[\xi_i])_{\mf{p}}$ and similarly define $G_j$. Because of our assumptions, $G_i, G_j \in \textup{GL}_n(\mb{Z}_p)$. Now,  $p\tau_i[u] = \tau_j$, gives $C_2^{t}G_i C_1 = p G_j$. If \eqref{condition rep} is satisifed, then $C_i \in M_n(\mb{Z}_p)$ and from the previous relation we also get $pC_i^{-1} \in M_n(\mb{Z}_p)$, showing that the elementary divisors of $C_i$ are $1$ and $p$. \\
    
    Now, if $\dim _{\mb{F}_p} \ker(\overline{C_1}) := k$, with $0 \leq k \leq n$, by taking determninants in $C_2^{t}G_i C_1 = p G_j$ and using the fact about the elementary divisors of $C_i$, we obtain $\dim _{\mb{F}_p} \ker(\overline{C_2}) = n-k$. Then, it is well-known (cf. proof \cite[Theorem 5.2.5]{andrianov_book}) that the number of $B$ satisfying the second condition in \eqref{condition rep} is equal to 
    \begin{equation*}
        {n-k\brack n-b}_p = {n-k\brack b-k}_p.
    \end{equation*}
    Similarly, by setting $A' := p (A^{-1})^{t}$, we may show that the number of $A'$ satisfying the first condition of \eqref{condition rep} is equal to ${k \brack a}_p$.
    Hence, we obtain 
    \begin{equation*}
        \sum_{L \in E_p\backslash E_pD_{a,b}E_p}\mathcal \#\mc{R}_{ij}(L) = \sum_{\substack{u \in M_n(K)\\p\tau_i[u] = \tau_j\\pu\mc{L}_j \subseteq \mc{L}_i}}{k\brack a}_p {n-k\brack b-k}_p =  \sum_{k=0}^{n} {k\brack a}_p {n-k\brack b-k}_p R_{ij}^{(k)}(p).
    \end{equation*}
    We may write $R^{(k)}_{ij}(p)$ since the map $u \longmapsto pu$ is a bijection. Substituting into \eqref{theta relation intermediate}, and using that 
    \[
    \sum_{a+b=r} p^{-ab}{k\brack a}_p{n-k\brack b-k}_p = p^{-k(r-k)}{n\brack r-k}_p,
    \]
    we obtain the result.\\

    Assume now $m \neq n$. We keep $m$ fixed. For each $s \geq 1$, put
\begin{equation*}
    \Theta_r^{(s)}
    :=
    \chi(\mathfrak p)^{s-r}
    \sum_{k=0}^{m}
    p^{-k(r-k+m-s)}
    {2s-m\brack r-k}_p M_k(p),
    \qquad 0 \leq r \leq 2s.
\end{equation*}
We will show that these matrices satisfy
\begin{equation*}
    \mbf{t}_s \mid T_r^{(s)}(p)
    =
    \Theta_r^{(s)}\mbf{t}_s.
\end{equation*}
This has already been proved when $s=m$. For $s \geq 2$, comparing coefficients of the standard polynomials \eqref{relation of standard different genus} gives
\begin{equation}\label{psi on split generators}
\begin{aligned}
    \Psi_p^{(s)}(T_r^{(s)}(p))
    =
    \chi(\mf{p})p^r T_r^{(s-1)}(p)
    +(1+p^{2s-m-1})T_{r-1}^{(s-1)}(p)+\chi(\ol{\mf{p}})p^{2s-r}T_{r-2}^{(s-1)}(p),
\end{aligned}
\end{equation}
where $\Psi_p^{(s)}$ denotes the map of Definition \ref{psi definition global}. Here, operators with indices outside $0,\ldots,2s-2$ are interpreted as zero. We can now check that the matrices satisfy
\begin{equation}\label{theta matrix recurrence different degrees}
    \Theta_r^{(s)} = \chi(\mf{p})p^r\Theta_r^{(s-1)} +(1+p^{2s-m-1})\Theta_{r-1}^{(s-1)} +\chi(\ol{\mf{p}})p^{2s-r}\Theta_{r-2}^{(s-1)}.
\end{equation}
% For descent, we use this relation for $0 \leq r \leq 2s-2$.
% If $s>m$, it also holds for the two upper indices, with the
% out-of-range matrices in degree $s-1$ interpreted as zero.\\
Assume first that $m>n$. We proceed downwards from degree $m$.
Suppose that the formula holds in degree $s$. Applying Theorem
\ref{Psi and Phi}, we obtain
\begin{equation*}
    \mbf{t}_{s-1}\mid\Psi_p^{(s)}(T_r^{(s)}(p))
    =
    \Theta_r^{(s)}\mbf{t}_{s-1}.
\end{equation*}
We now prove the formula in degree $s-1$ by induction on
$r=0,\ldots,2s-2$. Substituting \eqref{psi on split generators},
using the already established formulas for the smaller indices,
and then applying \eqref{theta matrix recurrence different degrees},
gives
\begin{equation*}
    \chi(\mf{p})p^r
    \left(
        \mbf{t}_{s-1}\mid T_r^{(s-1)}(p)
        -
        \Theta_r^{(s-1)}\mbf{t}_{s-1}
    \right)=0.
\end{equation*}
Since $\chi(\mf{p})p^r \neq 0$, we obtain
\begin{equation*}
    \mbf{t}_{s-1}\mid T_r^{(s-1)}(p)
    =
    \Theta_r^{(s-1)}\mbf{t}_{s-1}.
\end{equation*}
Repeating this argument proves the result in degree $n$.\\

Assume now that $m<n$. We proceed upwards from degree $m$.
Suppose that the formula holds in degree $s-1$, where $s>m$.
By Theorem \ref{Psi and Phi}, identities
\eqref{psi on split generators} and
\eqref{theta matrix recurrence different degrees}, and the
induction hypothesis, we obtain
\begin{equation*}
    \Phi^{\textup{Sie}}
    \left(
        \mbf{t}_s\mid T_r^{(s)}(p)
        -
        \Theta_r^{(s)}\mbf{t}_s
    \right)=0.
\end{equation*}
Each entry of the expression in parentheses has Fourier indices of rank at most $m$, since this holds for the theta series and is preserved by the Hecke action. Thus its image under $(\Phi^{\textup{Sie}})^{s-m}$ also vanishes. The injectivity argument used in the proof of Theorem \ref{Theorem m neq n} then implies
\begin{equation*}
    \mbf{t}_s\mid T_r^{(s)}(p)
    =
    \Theta_r^{(s)}\mbf{t}_s.
\end{equation*}
Induction proves the formula in degree $n$, completing the proof.
\end{proof}

Next, let $p \nmid \mf{bc}$ be inert, and let $T_p := \mbf{D}_p \textup{diag}(p1_n, 1_n)\mbf{D}_p$. In this case, we define
\begin{equation*}
R_{ij}(p):= \#\left\{ u\in M_m(K)\mid  \tau_i[u]=p\tau_j,\quad u\mathcal L_j\subseteq\mathcal L_i \right\}, \qquad M(p):=\left(e_i^{-1}R_{ij}(p)\right)_{i,j=1}^{H}.
\end{equation*}
\begin{theorem}
    We have
    \begin{equation*}
        \Theta(T(p)) =
    \begin{cases}
        \displaystyle
        \frac{\prod_{j=0}^{n-1}(1+p^{2j+1-m})}
             {\prod_{j=1}^{m/2}(1+p^{2j-1})}M(p),
        &\qquad \textup{if } m \textup{ is even},\\[12pt]
        0,
        &\qquad \textup{if } m \textup{ is odd and } m<2n+1.
    \end{cases}
\end{equation*}
If $m \geq 2n+1$ is odd, then $T(p)$ does not belong in the Hecke algebra described by Corollary \ref{theorem}.
\end{theorem}
\begin{proof}
    This is similar to the proof of Theorem \ref{eichler matrices split} and \cite[Theorem 5.2.5]{andrianov_book}. The fact that the action of $T(p)$ is trivial in the case $m=n$ odd follows from the fact that counting representations of the form $p\tau_i[u] = \tau_j$ as in \eqref{counting reps set} implies $m + 2v_{\mf{p}}(\det u) =0 $ after taking determinants. This then gives the zero action for all $m < 2n+1$ odd.
\end{proof}
\appendix
\renewcommand{\thesubsection}{\thesection\arabic{subsection}}
\section{Proofs for inert primes}\label{appendix inert}
We prove the local identities used in Section \ref{Section:Inert primes}, keeping the notation of that Section. The argument for the factorisation follows \cite[Section 3.5]{andrianov_book}.

\subsection{The parabolic Satake map}\label{appendix inert Frobenius}

Let $\mathcal C_\pm$ be the centralisers of $\Pi_\pm$ in $\mathcal H_{0,p}$ and put
\[
 \mathcal V_-:=\mathcal C_-\mathcal H_p,\qquad
 \mathcal V_+:=\mathcal H_p\mathcal C_+,
\]
where products denote linear spans. We use the injectivity result of \cite[Theorem 5.17 and Proposition 5.21]{heyer_decomposition}. The standard self-dual lattice gives the hyperspecial subgroup $\Gamma$, and the relative roots in the Siegel unipotent radical are $e_i+e_j$ and $2e_i$ in the root system of type $C_n$. Their pairwise inner products are nonnegative, so the Siegel parabolic is non-obtuse in Heyer's sense.\\

For his strictly positive element we take $a_P=\textup{diag}(p1_n,1_n)$, which is central in the Levi and satisfies $a_P n(H)a_P^{-1}=n(pH)$. Its Hecke element is $\Pi_-$. Also, $\Pi_+=[p]\langle a_P^{-1}\rangle_P$, with $[p]$ central and invertible. Thus our $\mathcal C_-,\mathcal C_+$ are his $\mathcal C_P^+,\mathcal C_P^-$, respectively, and our $\mathcal V_-,\mathcal V_+$ are his $\mathcal O_P^+,\mathcal O_P^-$. His projection $\Theta_M^P$ is our $\Omega$ by the formula on left cosets. It follows that $\Phi_n$ is injective on $\mathcal V_\pm$ and hence on $\mathcal C_\pm$.\\

We use the $*$-anti-automorphism induced by $g\longmapsto\mu(g)g^{-1}$. It interchanges $\Pi_-$ and $\Pi_+$, sends $\Pi_i$ to $\Pi_{n-i}$, and fixes the spherical Cartan double cosets by interchanging the paired exponents. For an integral triangular matrix of similitude $p$ we have $pBD^{-1}=BA^*$. Left multiplication by $\Pi_-$ therefore clears $BD^{-1}$, giving diagonal left-coset representatives. By \cite[Lemma 5.9]{heyer_decomposition}, $\Pi_-\Pi_i$ belongs to $\mathcal C_-$. It follows, using $*$ and the centrality of $[p]$, that $x_{-i}\in\mathcal C_-$ and $x_{+i}\in\mathcal C_+$.\\

Moreover, anlogosouly to \cite[Lemma 3.3.34]{andrianov_book}, we have
\begin{equation}\label{inert Satake values}
 \Phi_n(\Pi_i)=x_0e_i(y),\qquad
 \Phi_n(C)=x_0^2e_n(y),\qquad \Pi_-\Pi_+=C,
\end{equation}
where $e_i$ denote the elementary symmetric polynomial of degree $i$ and $C:=p^{n^2}[p]$.
\subsection{Factorisation of the standard polynomial}\label{appendix inert factorisation}

We follow \cite[Proposition 3.5.19 and Theorem 3.5.23]{andrianov_book}. The outer factors $X_\pm$ are already given, so we first divide $R_p$ on the right by $X_+$. The relation $\Pi_-\Pi_+=C$ gives $X_\pm(0)=1$, and we can set $Y_-:=R_pX_+^{-1}$ in the formal power series ring. The anti-automorphism $*$ fixes $R_p$ and interchanges $X_-$ and $X_+$, so with $Y_+=Y_-^*$ we have
\begin{equation}\label{inert two factors}
 R_p(v)=Y_-(v)X_+(v)=X_-(v)Y_+(v).
\end{equation}
Thus the middle factor must be $B(v)=X_-^{-1}Y_-=Y_+X_+^{-1}$. By \eqref{inert signature}, the coefficient of degree $i$ in $X_\pm$ has signature $\pm2i$, and the same holds for $X_\pm^{-1}$. If $Y_-$ has only nonpositive signatures, then $Y_+$ has only nonnegative signatures, since $*$ reverses the signature. The two expressions for $B$ then show that it has signature zero. In $X_-^{-1}Y_-$ only the constant term of $X_-^{-1}$ can contribute to signature zero. Hence it suffices to compute $Y_-$ and show that it has only nonpositive signatures. Its signature-zero part will give $B_p$.\\

As $X_+^{-1}$ has coefficients in $\mathcal C_+$, the coefficients of $Y_-$ lie in $\mathcal V_+$. By \eqref{inert Satake values}, $\Phi_n(x_{+i})=e_i(y)$, so the defining Satake image of $R_p$ gives
\[
 \Phi_n(Y_-(v))=\prod_{j=1}^n(1-y_j^{-1}v).
\]
Injectivity on $\mathcal V_+$ shows that $Y_-$ is a polynomial and that its coefficient of $(-v)^i$ is determined by the image $e_i(y^{-1})$, where $y^{-1}=(y_1^{-1},\ldots,y_n^{-1})$. We seek these coefficients among sums of parabolic cosets whose signatures can be read off from the elementary divisors.\\

For $a,b\geq0$, $a+b\leq n$, and $0\leq r\leq a$, set
\begin{equation*}
 D_{a,b}:=\textup{diag}(1_{n-a-b},p1_a,p^21_b),\qquad H_r:=\textup{diag}(1_r,0_{a-r})
 \end{equation*}
 \begin{equation*}
 \Pi_{a,b}^{(r)}
 :=\left\langle\m{p^2D_{a,b}^{-1}&\textup{diag}(0_{n-a-b},H_r,0_b)\\
                      0&D_{a,b}}\right\rangle_P.
    \end{equation*}
Put also
\[
 S_{h,c}:=\sum_{a=n-h}^{n-c}\Pi_{a,c}^{(h+a-n)}
 \qquad(0\leq c\leq h\leq n).
\]
Following \cite[Proposition 3.5.19, equation (3.5.56)]{andrianov_book}, for $0\leq i\leq n$ we set
\[
 u_i^-:=p^{-i(2n-i)}[p]^{-1}\sum_{j=0}^i\alpha_{ij}(p)S_{i-j,0}.
\]
We will show that $u_i^-$ belongs to $\mathcal V_+$, has only nonpositive signatures, and satisfies $\Phi_n(u_i^-)=e_i(y^{-1})$. These properties identify $u_i^-$ with the coefficient of $(-v)^i$ in $Y_-$ and allow us to read off the middle factor.\\

By \eqref{inert signature}, the cosets $[p]^{-1}\Pi_{a,0}^{(r)}$ have signature $2(a-n)$, proving the signature assertion. To prove that $u_i^-\in\mathcal V_+$ it suffices to show that $S_{h,0}\in\mathcal V_+$, since $[p]^{-1}$ is central.\\

Let $g$ be a representative of the parabolic double coset $\Pi_{n-h,0}^{(0)}$. For the spherical double coset $\Gamma g\Gamma$, the analogue of \cite[equation (3.3.61) and Lemma 3.5.16]{andrianov_book} gives
\[
 \iota_p(\Gamma g\Gamma)
 =\sum_{c=0}^h S_{h,c}\in\mathcal V_+.
\]

As in \cite[Proposition 3.5.19, equations (3.5.59)--(3.5.61)]{andrianov_book}, we show $S_{h,c}\in\mathcal V_+$ by induction on $h-c$. For $h=c$ we use $\Pi_c\Pi_+=p^{(n-c)^2}S_{c,c}\in\mathcal C_+$, the analogue of \cite[Proposition 3.5.17, equation (3.5.34)]{andrianov_book}, where the power counts the Hermitian residue block of order $n-c$. At difference $h>0$, the spherical sum gives $S_{h,0}$, since the other summands have smaller difference. To obtain $S_{h+b,b}$ for $1\leq b\leq n-h$, we multiply $S_{h,0}$ by $\Pi_{n-b,b}^{(0)}$ on the right. When $h<n$, this completes the cases of difference $h$ and supplies $S_{h+1,1}$, which is needed to isolate $S_{h+1,0}$ in the next spherical sum.\\

For the product $\Pi_{a,0}^{(r)}\Pi_{n-b,b}^{(0)}$ we use the indices of \cite[Proposition 3.5.17, equation (3.5.35)]{andrianov_book}. Its possible terms are $[p]\Pi_{d,t}^{(R)}$, where $d=a+b-2t$ and $R=r-s$. The lower-right block counts use $\mathbb F_{p^2}$, and the unipotent residue blocks are Hermitian. Combining the lower-right count with the ratio $h_p(R,a-t)/h_p(R,d)$ from \eqref{inert rank count} leaves the factor
$\varphi_{d-R}(p^2)/\bigl(\varphi_{a-t-R}(p^2)\varphi_{b-t}(p^2)\bigr)$.
The unipotent indices give $p^{2bn-2t(d+t)}$. The entries between the nonsingular block of order $R$ and the $t$ additional rows and columns give $p^{2tR}$ choices. Clearing these entries by row and column operations leaves a Hermitian matrix of size $a-R$ and rank $s$ with a zero block of size $a-t-R$. Thus the coefficients depend on $a,r$ only through $a-r=n-h$, so the terms group into $[p]S_{h+2t-b-s,t}$ with $0\leq t\leq b$ and $s\geq0$. For $t=b$, $s=0$, only the power $p^{2bn+2br-2ab}=p^{2bh}$ remains. Thus this term is $p^{2bh}[p]S_{h+b,b}$, and all other terms have smaller index difference. Right multiplication by $\mathcal C_+$ preserves $\mathcal V_+$, and $[p]$ is central and invertible. This completes the induction.\\

To compute $\Phi_n(u_i^-)$, let $h_p(r,a)$ count the Hermitian matrices of size $a$ and rank $r$ over $\mathbb F_{p^2}$. In the analogue of \cite[Lemmas 3.5.20 and 3.5.21]{andrianov_book}, we choose the radical over $\mathbb F_{p^2}$. On the quotient of dimension $r$, the nondegenerate Hermitian forms constitute one congruence class with stabiliser $U_r(\mathbb F_p)$, so there are $|\textup{GL}_r(\mathbb F_{p^2})|/|U_r(\mathbb F_p)|$ such forms. This gives
\begin{equation}\label{inert rank count}
 h_p(r,a)=\frac{\varphi_a(p^2)}{\varphi_r(p^2)\varphi_{a-r}(p^2)}
        p^{r(r-1)/2}\prod_{\ell=1}^r(p^\ell+(-1)^\ell).
\end{equation}

Over each left coset of the lower-right block of type $D_{a,0}$ there are $h_p(r,a)$ unipotent choices. Together with the Satake image of the lower-right double coset this gives
\[
 \Phi_n(\Pi_{a,0}^{(r)})=p^{-a^2}x_0^2 h_p(r,a)e_a(y).
\]
Putting $k=n-a$ and dividing by $\Phi_n([p])$, we obtain
\[
 \Phi_n([p]^{-1}S_{h,0})
   =\sum_{k=0}^h p^{k(2n-k)}h_p(h-k,n-k)e_k(y^{-1}).
\]
Thus the coefficient of $e_k(y^{-1})$ in $\Phi_n(u_i^-)$ is $p^{k(2n-k)-i(2n-i)}$ times the sum below. The Hermitian analogue of \cite[equation (3.5.63)]{andrianov_book} is
\begin{equation}\label{inert triangular inversion}
 \begin{aligned}
 \sum_{j=0}^{i-k}\alpha_{ij}(p)h_p(i-j-k,n-k)
 &=\frac{\varphi_{n-k}(p^2)}{\varphi_{n-i}(p^2)\varphi_m(-p)}
    \prod_{\ell=0}^{m-1}(1-(-p)^\ell)\\
 &=\begin{cases}1&k=i,\\0&0\leq k<i,\end{cases}
 \end{aligned}
\end{equation}
where $m=i-k$. The first equality follows by substituting \eqref{inert rank count} and using the finite Gaussian binomial identity with parameter $-p$. This gives $\Phi_n(u_i^-)=e_i(y^{-1})$.\\

By injectivity, $Y_-(v)=\sum_{i=0}^n(-1)^iu_i^-v^i$. Its signature-zero part comes from $a=n$, where $[p]^{-1}\Pi_{n,0}^{(r)}=Z_r$. Hence the middle factor is
\[
 B(v)=\sum_{i=0}^n(-1)^ip^{-i(2n-i)}
    \sum_{j=0}^i\alpha_{ij}(p)Z_{i-j}v^i=B_p(v).
\]

\subsection{Inert inverse factors}\label{proof inert inverse factors}
This is the analogue of \cite[Lemma 3.5.24]{andrianov_book}, and can be shown in a similar way. Namely, we first observe that it suffices to show one equality, since the other follows by applying the $*$-anti-automorphism. So it is enough to show
\[
 X_{-}(v)^{-1} = \sum_{r=0}^{\infty} p^{-(2n-1)r}t^{-}(\mathfrak{p}^{r})v^{r}.
 \]
The operators $t^-(\mathfrak{p}^r)$ as well as $x_{-i}$ are all in $\mathcal{C}_-$ and hence it is enough to show the equality after applying the map $\Omega$, which we know is injective on $\mathcal{C}_-$. That is, we need to verify the identity
\[
 \left(\sum_{i=0}^{n}(-1)^i\Omega(x_{-i})v^i\right)^{-1} =\sum_{r=0}^{\infty} p^{-(2n-1)r}\Omega(t^{-}(\mathfrak{p}^{r}))v^{r}.
\]

As in Andrianov, this can be reduced to a statement about the Hecke algebra $\mathcal{H}(\textup{GL}_n(\mathcal{O}_{K_\mathfrak{p}}), \textup{GL}_n(K_{\mathfrak{p}}))$, and in particular
\[
\sum_{i=0}^\infty t(\mathfrak{p}^i) v^i = \left( \sum_{i=0}^{n}(-1)^i p^{(i-1)i}\pi_i(\mathfrak{p})v^i\right)^{-1},
\]
where $\pi_i(\mathfrak{p}) = \textup{GL}_n(\mathcal{O}_{\mathfrak{p}}) \textup{diag}(1_{n-i},p1_{i})\textup{GL}_n(\mathcal{O}_{\mathfrak{p}})$, and
$t(\mathfrak{p}^i) := \sum_{\textup{det}(g)\mathcal{O}_{\mathfrak{p}} = \mathfrak{p}^i} \textup{GL}_n(\mathcal{O}_{\mathfrak{p}}) g$.
This is the analogue of \cite[Proposition 3.2.22]{andrianov_book} for $\textup{GL}_n$ over the local field $K_{\mathfrak{p}}$, whose residue field has cardinality $p^2$. Write $j$ for inversion on double cosets. We have
\[
 \Omega_p(t^-(\mathfrak p^d))=j(t(\mathfrak p^d)),
 \qquad
 \Omega_p(x_{-i})=p^{-i(2n-i)}j(\pi_i(\mathfrak p)).
\]
The first equality follows from the diagonal representatives for $t^-$, since $\widehat D$ and $D^{-1}$ have the same elementary divisors. For the second, use $\Omega_p(\Pi_j)=p^{j^2}x_0\pi_j(\mathfrak p)$ and $\Omega_p([p])=x_0^2E_p(p1_n)E_p$ in the definition of $x_{-i}$. Applying $j$ gives the required equality.

\subsection{Action of the middle factor}\label{appendix inert middle}
As in the split case, we first reduce the action to a finite character sum. We recall that
\[
 \epsilon_{\mathbb A}((q,s);(\mc L,\tau))
 =\beta(q)\sum_{u\in\textup{GL}_n(K)}
   \lambda'(\xi^{-1}uq)\Phi(u,q,s,\tau).
\]
All Hermitian matrices of a fixed rank over $\mathbb F_{p^2}$ are congruent, so the operator $Z_a$ decomposes as
\[
 Z_a=\sum_{\substack{H\in\textup{Herm}_n(\mathbb F_{p^2})\\\textup{rk}\,H=a}}
       \Gamma_P n(p^{-1}\widetilde H),
\]
where $\widetilde H$ is any integral Hermitian lift of $H$. Then
\[
 (q,s)n(p^{-1}\widetilde H)^{-1}=(q,s-p^{-1}q\widetilde Hq^*)
\]
at $p$. The character in Definition \ref{hecke action} is trivial on these unipotent elements, so acting with this coset gives the additional factor
\[
 e_p(-p^{-1}\textup{tr}(T_u\widetilde H)),\qquad
 T_u=(q^*u^*\tau uq)_{\mathfrak p}.
\]
Let $A=(\xi^*\tau\xi)_{\mathfrak p}$, which is integral, and let $r$ be the rank of $\overline A$ over $\mathbb F_{p^2}$. For each nonzero $u$-term of $\epsilon_{\mathbb A}$, we have $h_u=(\xi^{-1}uq)_{\mathfrak p}\in E_p$ and $T_u=h_u^*Ah_u$. Thus $T_u$ is integral and its reduction has rank $r$ for every such summand.\\

Define the nontrivial character $\overline\psi_p(t)=e_p(-p^{-1}\widetilde t)$ of $\mathbb F_p$, where $\widetilde t$ is any lift to $\mathbb Z_p$. For a Hermitian matrix $T$ over $\mathbb F_{p^2}$ put
\[
 S_a^{(n)}(T)=
 \sum_{\substack{H\in\textup{Herm}_n(\mathbb F_{p^2})\\\textup{rk}\,H=a}}
 \overline\psi_p(\textup{tr}(TH)).
\]
The trace lies in $\mathbb F_p$. Since congruence of $T$ permutes the rank-$a$ matrices in this sum, the sum depends only on the rank of $T$. We write $S_a^{(n)}(r)$ for its value at rank $r$.\\

With this notation, the action of $Z_a$ multiplies each nonzero $u$-term by $S_a^{(n)}(r)$. Thus the action of $B_p(v)$ is reduced to computing
\[
 F_{n,r}(v)=\sum_{i=0}^n(-1)^iB_i^{(n)}(r)v^i,
 \qquad
 B_i^{(n)}(r)=p^{-i(2n-i)}\sum_{j=0}^i\alpha_{ij}(p)S_{i-j}^{(n)}(r),
\]
where $\alpha_{ij}(p)$ is as in Proposition \ref{standard polynomial inert}. When $r=0$, we have $S_a^{(n)}(0)=h_p(a,n)$, and \eqref{inert triangular inversion} with $k=0$ gives $F_{n,0}(v)=1$. We take $F_{0,0}(v)=1$.\\

To compute $F_{n,r}$ for $r>0$, we compare the character sums at ranks $r$ and $r-1$. This is the Hermitian analogue of the calculation in \cite[Proposition 4.2.24 and Theorem 4.2.27]{andrianov_book}. Take the indices $\textup{diag}(1,A')$ and $\textup{diag}(0,A')$, where $\textup{rk}\,A'=r-1$, and write $H=\m{c&w^*\\w&C}$. In the difference of the sums, the terms with $c=0$ cancel. For $c\ne0$, clearing the off-diagonal entries in the first row and column gives $\textup{diag}(c,C-c^{-1}ww^*)$. Thus, for $a\geq1$, summing over $C-c^{-1}ww^*$ of rank $a-1$ gives
\[
 S_a^{(n)}(r)-S_a^{(n)}(r-1)
 = S_{a-1}^{(n-1)}(r-1)
   \sum_{c\in\mathbb F_p^\times}(\overline\psi_p(c)-1)
   \sum_{w\in\mathbb F_{p^2}^{\,n-1}}
        \overline\psi_p(c^{-1}w^*A'w).
\]
After diagonalising $A'$, each of its $r-1$ nonzero coordinates contributes $-p$, since for $t\in\mathbb F_p^\times$ the norm map gives
\[
 \sum_{z\in\mathbb F_{p^2}}\overline\psi_p(tz\bar z)
 =1+(p+1)\sum_{x\in\mathbb F_p^\times}\overline\psi_p(tx)=-p.
\]
The other $n-r$ coordinates contribute $p^2$ each, and the sum of $\overline\psi_p(c)-1$ over $c\ne0$ is $-p$. Hence
\[
 S_a^{(n)}(r)-S_a^{(n)}(r-1)
 =(-p)^{2n-r}S_{a-1}^{(n-1)}(r-1).
\]

In the difference of the coefficients $B_i^{(n)}(r)$ and $B_i^{(n)}(r-1)$, the term $j=i$ cancels since $S_0^{(n)}(r)=1$. For $j\leq i-1$, the coefficients $\alpha_{ij}(p)$ depend on $n$ and $i$ through $n-i$, so they are the same in dimension $n-1$ at index $i-1$. Substituting the last identity and using
$p^{-(2n-1)}(-p)^{2n-r}=-(-p)^{1-r}$ gives
\[
 F_{n,r}(v)-F_{n,r-1}(v)=(-p)^{1-r}vF_{n-1,r-1}(v).
\]
Induction on $r$, for all $n\geq r$, now gives
\[
 F_{n,r}(v)=\prod_{j=0}^{r-1}\bigl(1+(-p)^{-j}v\bigr),
\]
as required.
\printbibliography
\end{document}